\documentclass{acta_2017shortcap}
\usepackage{harvard_acta}
\usepackage{amssymb}
\usepackage{amsmath}
\usepackage{graphicx}
\usepackage[dvipsnames]{xcolor}
\newcommand{\ron}[1]{{#1}} 
 
\newcommand{\guergana}[1]{{\color{red}{#1}}} 
\newcommand{\simon}[1]{{\color{blue}{#1}}} 
\newcommand{\boris}[1]{{\color{purple}{#1}}} 

\def\<{\langle}
\def\>{\rangle}

\def\e{\varepsilon}

\def\wt{\widetilde}

\def\Z{{\mathbb Z}}
\def\cS{{\cal S}}
\def\cL{{\cal L}}
\def\cA{{\cal A}}
\def\cC{{\cal C}}

\def\cN{{\cal N}}

\def\cV{{\cal V}}
\def\cM{{\cal M}}
\def\cF{{\cal F}}
\def\cP{{\cal P}}

\def\cD{{\cal D}}
\def\cP{{\cal P}}

\def\cX{{\cal X}}
\def\R{\mathbb{R}}
\def\N{\mathbb{N}}

\def\Relu{\mathrm{ReLU}}
	
\def\Up{\overline{\underline{\Upsilon}}}

\newcommand{\be}{\begin{equation}}
\newcommand{\ee}{\end{equation}}
\newcommand{\gives}{\rightarrow}
\newcommand{\set}[1]{\left\{#1\right\} }
\newcommand{\mN}{\mathcal N}
\newcommand{\E}[1]{\mathbb E\left[#1\right]}
\newcommand\eref[1]{{\rm (\ref{#1})}}
\def\rnew{\color{red} }
\def\bnew{\color{blue} }

\newcommand{\ignore}[1]{}
\newcommand{\lr}[1]{\left(#1\right)}
\newcommand{\inprod}[2]{\left\langle#1,#2\right\rangle}
\newtheorem{theorem}{Theorem}[section]
\newtheorem{lemma}[theorem]{Lemma}
\newtheorem{definition}[theorem]{Definition}
\newtheorem{proposition}[theorem]{Proposition} 
\newtheorem{remark}{Remark}[section] 

\newcommand{\abs}[1]{\left|#1\right|}
\newcommand{\norm}[1]{\left|\left|#1\right|\right|}

 \def\argmin{\mathop{\rm argmin}}

 \def\dist{\mathop{\rm dist}}

\title[]{
Neural Network Approximation}
\author[{Acta Numerica}]{
Ronald DeVore\\
{\tt rdevore@math.tamu.edu}
\and
Boris Hanin\\
{\tt bhanin@princeton.edu}
\and 
Guergana Petrova\\
{\tt gpetrova@math.tamu.edu}
}

\begin{document}

 \maketitle

\vskip 0.45in
\begin{abstract}
Neural Networks (NNs) are the method of  choice for building learning algorithms. They are now being 
 investigated for other numerical tasks such as solving high dimensional partial differential equations.
Their popularity stems from their empirical success on several challenging learning problems (computer chess/go, autonomous navigation, face recognition).  However, most scholars agree that a convincing  theoretical explanation for this success is still lacking. 
Since these applications revolve around approximating an unknown  function $f$ from
  data observations,   at least part of the answer to this success must involve the ability of NNs to produce accurate function approximations. 
  
  This article surveys  the known approximation properties of
the outputs of NNs with the aim of uncovering the properties that are not present
in the more traditional methods of approximation used in numerical analysis, such as approximations using polynomials, wavelets, rational function, and splines. Comparisons are made with traditional approximation methods from the viewpoint of rate distortion, i.e. error versus the number of parameters used to create the approximant. 
Another major component
in the  analysis of numerical approximation is  the computational time needed to construct the approximation  and  this in turn is intimately
connected with the stability of the approximation algorithm.   So the stability of numerical approximation using NNs is
a large part of the analysis put forward.

The survey, for the most part, is concerned with NNs using  the popular ReLU activation function.  In this case, the outputs
of the NNs are piecewise linear functions on rather complicated partitions of the domain of $f$ into cells that are convex polytopes.  When the architecture of the NN is fixed and the parameters are allowed to vary, the set of output
functions of the NN is a parameterized nonlinear manifold.  It is shown that this manifold has certain space filling
properties leading to an increased ability to approximate (better rate distortion) but at the expense of numerical
stability. The space filling  creates a challenge to the numerical method in finding best or good parameter choices when trying to approximate.

\end{abstract}

\tableofcontents 

\section{Introduction}
\label{Introduction}

Approximation using Neural Networks (NNs)   is  the method of choice for building numerical algorithms in Machine  Learning (ML)  and Artificial Intelligence (AI).  
It is now being looked at as a possible platform for computation in many other areas.  Although NNs have been around for over 70 years, starting with the work of Hebb in the late 1940's \cite{hebb1949organization} and  Rosenblatt in the 1950's \cite{rosenblatt1958perceptron}, it is only recently that  
  their popularity has surged as they have achieved state-of-the-art performance in a striking variety of machine learning problems. Examples of these are computer vision \cite{krizhevsky2012imagenet}, employed for instance in  self-driving cars, natural language processing \cite{wu2016google}, used in  Google Translate, or reinforcement learning, such as  superhuman performance at Go \cite{silver2016mastering,silver2017mastering}, to name a few. 
  
 Nevertheless, it is generally agreed upon that there  is  still a lack of solid mathematical analysis to explain the reasons behind these empirical successes. As a start, the understanding of the  approximation properties of NNs is of vital importance since approximation is one of the main components of any algorithmic design. A rigorous analysis of   what special properties  NNs hold as a method of approximation  could lead to both  significant practical improvements \cite{bronstein2017geometric,lecun2015deep} and a priori  performance guarantees for computational algorithms based on NNs.

At the heart of providing such a rigorous theory is understanding the benefits of  using NNs as an approximation tool  when compared with  other more classical methods of approximation such as polynomials, wavelets, splines, and sparse approximation from
bases, frames,  and dictionaries.  Indeed, most applications of NNs are built  on some form of function approximation.  This includes not only learning theory and statistical estimation,  but also the new forays of NNs into other application domains such
as numerical methods for solving partial differential equations (PDEs).

An often cited theoretical feature of neural networks is that they produce universal function approximants \cite{cybenko1989approximation,hornik1989multilayer} in the sense that, given any continuous target function $f$ and a target accuracy $\epsilon>0$, neural networks with enough judiciously chosen parameters produce an approximation to $f$ within an error of size $\epsilon$. This universal approximation capacity has been known since the $1980$'s.  But surely, this    cannot be the main reason why neural networks are so effective in practice. 
Indeed, all  families of functions used in numerical approximation such as polynomials, splines, wavelets, etc., produce  universal approximants.  What we need to understand is in what way  NNs are more effective than other methods as an approximation tool.

The purpose of the present article is to describe the approximation properties of NNs as we presently understand them,  and to compare their performance with other  methods of approximation.
To accomplish such a comparative  analysis, we introduce, starting in \S \ref{S:optimalperformance}, the tools by which various methods of approximation are evaluated.  These include approximation rates on model classes, $n$--widths, metric entropy,  and approximation classes.  Since NN approximation is a form of nonlinear manifold approximation, we make this particular form of approximation the focal point of our exposition. The ensuing sections of the paper examine the specific approximation properties of NNs.  After making some remarks that apply to general activation functions $\sigma$,  we turn our attention to the performance of    {\it  Rectified Linear Unit}  (ReLU) networks. These are the most heavily used in numerical settings and fortunately also the NNs most amenable to analysis. 

Since the output of a ReLU network is a continuous piecewise linear function (CPwL), it is important to understand what the class  of outputs  of a  ReLU network depending on $n$ parameters looks like in terms of their allowable partitions and the correlation between the linear pieces. This  topic is addressed in \S\ref{S:Relunetworks}.   This structure increases in complexity with the depth of the network.  It turns out that deeper NNs  give a richer  set of outputs  than shallow networks.   Therefore, much of  our analysis centers on deep ReLU networks. 

The key takeaways from this paper are the following. For a fixed value of $n$, the  outputs of  ReLU networks depending on $n$ parameters  form a rich  parametric family  of CPwL  functions.  This manifold exhibits certain space filling  properties (in the Banach space where we measure performance error), which are both  a  boon and a bottleneck. On one hand,  space filling  provides the possibility  to  approximate with relatively few parameters  larger classes of functions than the classes that are currently approximated by classical methods. On the other hand,  this flexibility comes at the expense of both the stability of the algorithm by which one selects the right parameters  and the a priori performance guarantees and uncertainty quantifications of performance when using NNs in numerical algorithms. This points to the need for  a comprehensive study of the trade-offs  between stability   of numerical algorithms based on NNs and their numerical efficiency.
 
 This exposition is far from providing a satisfactory theory for approximation by NNs, even when we restrict ourselves to ReLU networks.
 We highlight several fundamental questions that remain unanswered.  Their solution would not only lead to a   better understanding  of NN approximation but would most likely guarantee  better performance in numerical algorithms.  These issues include:
 \begin{itemize}
 \item matching upper and lower for the rate of approximation of standard model classes when using  ReLU networks
  \item how to precisely describe  the types of function classes that benefit from NN approximation.
 \item how to   numerically impose  stability in parameter selection; 
 \item how the imposition of stability limits the   performance of the network;  

\end{itemize}

\section{What is a Neural Network?}

This section begins by introducing feed-forward neural networks  and their elementary properties.  We begin with a general setting and  then specialize to the case of
fully connected networks. While the latter networks are generally not the architecture of choice in most targeted applications, their architecture provides the most convenient way to understand the trade-offs between approximation efficiency and the complexity of the network.  They also allow for a clearer picture
of the balance between width and depth in the assignment of parameters.

In its most general formulation, a  {\it feed-forward neural network} $\mathcal N$ is associated with  a directed acyclic graph (DAG),
\[\mathcal G = \lr{{\cal V},{\cal E}},\]
called the {\it architecture} of $\mathcal N$, determined by a finite set ${\cal V}$ of vertices and a finite set of directed edges ${\cal E}$, in which every vertex $v\in {\cal V}$ must belong to at least one edge $e\in {\cal E}$. The set  ${\cal V}$ consists of  three distinguished subsets. The first is the set ${\cal I}$ of input vertices. These vertices have no incoming edges and are placeholders for independent variables (i.e. network inputs). The second is the set ${\cal O}$ of output vertices. These vertices have no outgoing edges and will store, for given inputs, the corresponding value of the dependent variables (i.e. the network output). The third is the set of hidden vertices ${\cal H}= {\cal V}\backslash \set{{\cal I},{\cal O}}$. For a given  input, hidden vertices store certain intermediate values used to compute the corresponding output.
The vertices and edges also have the following adornments:
\begin{enumerate}
    \item[(1)] With every 
    $v\in {\cal V}\setminus{\cal I}$, there is  an associated function $\sigma_v:\R\gives \R$, called an {\it activation function}, and a scalar $b_v\in \R$, called a {\it }bias. 
    \item[(2)] For every $e\in {\cal E}$,   there is a scalar $w_e\in \R,$ called a {\it weight}. 
\end{enumerate}

In going forward, we often refer to the vertices as nodes. The weights and biases are referred to as the {\it trainable parameters} of ${\cal N}$. For a fixed network architecture, varying the values of these trainable parameters produces a family of output functions. The key to describing how these functions are constructed is that to each vertex $v\in {\cal V}\setminus{\cal I}$  we  associate a computational unit called a neuron. This unit takes as inputs the scalar outputs $x_{v'}$ from  vertices $v'\in {\cal V}\setminus {\cal O}$ with an edge
$e=(v',v)\in {\cal E}$ terminating at $v$, and outputs the scalar
\begin{equation}\label{E:neuron-comp}
    x_v:=\sigma_v\lr{b_v+\sum_{e=(v',v)\in {\cal E}} w_e x_{v'}}.
\end{equation}

 The word {\it neuron} comes from the fact that  \eqref{E:neuron-comp} can be viewed as a simple computational model for a single biological neuron. A neuron associated to a vertex 
 $v\in {\cal V}\setminus{\cal I}$ observes signals $x_{v'}$ computed by upstream neurons associated to $v'$, takes a superposition of these signals, mediated by synaptic weights $w_e$, $e=(v',v)$, and outputs  $x_v$  which is then seen by the downstream neurons.
For all neurons associated to vertices $v\in {\cal O}$, the  activation function $\sigma_v$ is the identity. The neuron associated to the $i$-th input vertex $v\in {\cal I}$, $i=1,\ldots,d$, where $d:=|{\cal I}|$, observes a scalar incoming (i.e. externally provided) signal $x_i$ and outputs $x_i$, which is then seen by the downstream neurons.

We view the network scalar inputs $x_i$, $i=1,\ldots,d$,   as an independent  variable  $x=(x_1,\ldots,x_d)\in \Omega\subset \R^d$, and define the output function $S_\cN$, $S_\cN:\Omega\gives\R^{d'}$ of the network $\cN$   by
\be
S_\cN(x):=\lr{x_v,~v\in {\cal O}}, \quad d':=|{\cal O}|.
\ee
Thus, $S_\cN$ is a function mapping 
$\Omega\subset\R^d$ into $\R^{d'}$, called the {\it output} of $\cN$. Note that for a fixed 
network architecture 
${\cal G} = ({\cal V},{\cal E})$,  the outputs $S_{\cal N}$ form a family of functions, determined by the trainable parameters 
$\{w_e,b_v\}$, $e\in {\cal E}$, $v\in {\cal V}\setminus{\cal I}$.

\subsection{Fully Connected Networks}
\label{SS:FCN}

The preceding is a very general definition of neural networks and encompasses virtually all   network architectures used in practice. In this article, however, we  restrict our study to rather special examples of such networks, the so-called {\it fully connected networks}. The architecture of such a network is given by  a directed acyclic graph whose vertices are organized into layers. 

Each vertex of every layer is connected via outgoing edges to all vertices from the next layer and to no other vertices from any other layer, see Figure \ref{F1_Ron}. The zero-th layer, called the {\it input layer},  consists of all $d:=n_0$ input vertices 
 ${\cal I}$, called inputs, where the $i$-th input,  receives scalar signal $x_i$ from outside the network.   The combined input
 $x:=(x_1,\ldots,x_d)$ forms the independent   variable of the function $S_\cN$.  The input layer is followed  by the hidden vertices ${\cal H}$, organized in $L$   {\it hidden layers}, with the $j$-th layer  consisting of  $n_j$ hidden vertices,  
 $j =1,\ldots,L$.  The integer $n_j$ is called the {\it width} of the   $j$-th  layer. Finally, the $(L+1)$-st layer, called   the
 {\it output layer}, 
 consists of all $d':=n_{L+1}$ output vertices 
 ${\cal O}$, called outputs.  The output vector of such a fully connected network is the value $S_{\cN}(x)\in \mathbb{R}^{d'}$ of the function $S_{\cN}$ for the input  $x$.

As is customary, we specify that there is a single  activation function $\sigma$ that is used at each hidden vertex $v\in {\cal H}$, i.e., $\sigma_v=\sigma$ for all  $v\in {\cal H}$. Recall that we always take the activation $\sigma_v$ at the output vertices $v\in {\cal O}$ to be the identity. In this way, each coordinate of $S_{\cN}(x)$ is a linear combination of the $x_v$'s  at layer $L$ plus a bias term, which is a constant. 
 
Thus, for a fully connected network $\cN$, the output function $S_\cN$ 
can be succinctly described by {\it weight matrices} and {\it bias vectors}
$$
W^{(\ell)}\in 
\R^{n_{\ell}\times n_{\ell-1}}, \quad b^{(\ell)}\in \R^{n_\ell}, \quad \ell=1,\ldots,L+1,
$$
associated to layer $\ell$ as follows.  
If $X^{(\ell)}\in \R^{n_\ell}$ is the vector of outputs  $x_v$ corresponding to nodes $v$ in  layer $\ell$, $\ell=0,\ldots,L+1$, then  the output $S_{\mN}$ is given by
\begin{equation}\label{E:output}
    S_\cN(x) = X^{(L+1)}=W^{(L+1)}X^{(L)}+b^{(L+1)},
\end{equation}
where the vectors  $X^{(\ell)}$ satisfy the recursion
\begin{eqnarray}
\label{E:FC-recursive}
X^{(\ell)}=\sigma(W^{(\ell)}X^{(\ell-1)}+b^{(\ell)}),\quad \ell =1,\ldots L,\qquad X^{(0)}=x.
\end{eqnarray}
Here and throughout this paper, we use the convention that the activation function $\sigma:\R \gives \R$ is defined to act on any vector $z=(z_1,\ldots,z_{n_\ell})\in\R^{n_\ell}$, $n_\ell\geq 1$ coordinatewise, that is
\[\sigma(z)= \sigma(z_1,\ldots, z_{n_\ell}):= (\sigma(z_1),\ldots, \sigma(z_{n_\ell})).\]

 \subsection{The Set $\Upsilon^{W,L}(\sigma;d,d')$}
 \label{SS:Upsilon}

We will almost always consider only fully connected feed-forward NNs whose hidden layer widths are all the same, namely,  $n_1=\dots=n_L=W$.  Note that  we can embed 
any fully connected feed-forward NN into a network with constant  width $W:=\max_{j=1,\ldots,L} n_j$ by inserting $(W-n_j)$ additional zero bias vertices to layer $j$ and adding new edges with weights set to $0$ between these vertices and those in the next layer
{if these vertices are in the first layer, we also add new edges with weights set to $0$ between them and the input vertices). We use this fact frequently in what follows, sometimes without mentioning it.

We refer to $W$
as the {\it width} of the network and to $L$ as its {\it depth}. In such networks, each vertex $v$ from a hidden layer can be associated with a pair of indices $(i,j)$, where $j$ is the layer index and $i$ is the row index of the location of $v$. We commonly refer to all vertices from a fixed row as a {\it channel}, and those from a fixed column as a {\it layer}.   It is useful to introduce for every vertex $v$ from the hidden layers  the function $z_v:=z_{i,j}$ which records how the value at this neuron depends on the original input   $x=(x_1,\dots,x_d)$ before the activation $\sigma$ is applied.
 It follows that,
\be
\label{interior}
\sigma(z_v(x_1,\dots,x_d)):=\sigma(z_{i,j}(x)) :=X^{(j)}_i,\quad  i=1,\dots,W; \ j=1,\dots,L,
\ee
which is the value of the 
$i$-th coordinate of the vector $X^{(j)}$ defined in \eqref{E:FC-recursive}.

For a fully connected feed-forward network $\cN$ with width $W$, depth $L$, activation function $\sigma$,    input dimension $d$, and output dimension $d'$, we define the set
\[\Upsilon^{W,L}:=\Upsilon^{W,L}(\sigma;d,d')\]
as  the  collection of all output functions $S_\cN$ that we obtain by varying the choice of the trainable parameters of $\cN$.   Recall that $S_\cN$
is a mapping from $\R^d$ (or $\Omega\subset \R^d$) to $\R^{d'}$. For notational simplicity, we often   omit the dependence of $\Upsilon^{W,L}$ on $\sigma, d$ and $d'$ when  these are understood from the context. 
Figure \ref{F1_Ron} shows the graph associated to a typical network that outputs functions from
$\Upsilon^{W,L}(\sigma;d,d')$ with $d=2,\, W=3,\, d'=1$.
\begin{figure}[h]
  \centering
\includegraphics[scale=.7]{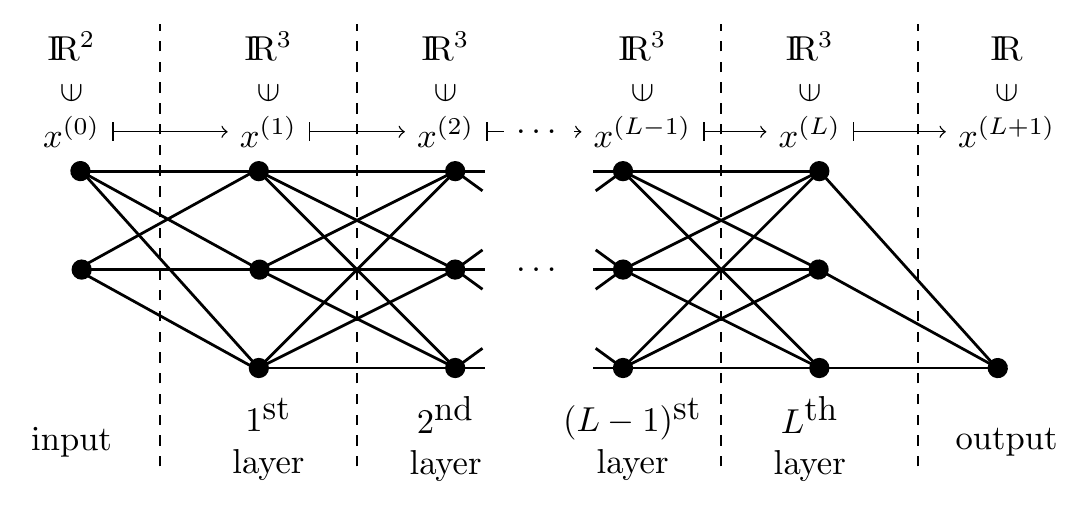}
\caption{The graph associated to the outputs $\Upsilon^{3,L}(\sigma;2,1)$ of a fully connected network with input dimension $2$, width $3$, $L$ hidden layers and output dimension $1$.}
\label{F1_Ron}
\end{figure}

Notice that $\Upsilon^{W,L}$ 
is closed under addition of weights and biases \textit{in the output layer}. This follows immediately from \eqref{E:output}.  However, it is \textit{not closed} under addition of functions because we can find two outputs $S_{\cN_1}, S_{\cN_2}$ from $\Upsilon^{W,L}$ with $S_{\cN_1}+S_{\cN_2} \not\in \Upsilon^{W,L}$. This will become apparent even in our discussion of one layer ReLU networks, see \S \ref{S:Relunetworks}. Therefore $\Upsilon^{W,L}$ is not a linear space. Each function  $S_\cN\in \Upsilon^{W,L}$ is determined by 
\be
\label{numberparameters}
n(W,L):= (d+1)W+W(W+1)(L-1)+ d'(W+1)
\ee
parameters consisting of the entries of  its weight matrices $W^{(1)},\ldots, W^{(L+1)}$ and bias vectors $b^{(1)},\ldots, b^{(L+1)}.$  We note in passing that it is possible that  distinct choices  of trainable parameters end up describing  the same   outputs.  

We  order the parameters of $\cN$ and organize them into a vector $\theta$, where $\theta=\theta_\cN \in \R^n$.  The output $S_\cN$ is then given by
\be
S_\cN= S(\cdot,\theta_\cN)=:M(\theta_\cN),
\ee
where $M:\R^n\to \Upsilon^{W,L}(\sigma;d,1)$ is the map from network weights and biases to output functions.  In this way, we view $\Upsilon^{W,L}$ as a parametric manifold.   Here, we are using the term `manifold' in a very loose sense because we are not
attributing any of the  topological or differential  properties  usually associated with this term.  

The mapping $M$ is completely determined once we have fixed the architecture and the activation function $\sigma$.  
Therefore, having made these choices, designing approximation methods for a target function $f$ boils down to choosing parameters $\theta=:a(f)$ when $f$ or information about $f$ is given. In this way,
any neural network based approximation  method consists of determining a mapping $a:f\mapsto a(f)$ that assigns to each potential target function $f$ a sequence of parameters $a(f)\in\R^{n}$.  The approximation to $f$ is then given by
\be
\label{appoperator}
A(f):= M(a(f)).
\ee

Our main focus in this paper is to understand the 
approximation power of NNs and thus we work under the assumption that we have full access to the target function  $f$.  However, in the last two  sections, we do make
forays into the more realistic (numerical) settings where we are only provided (partial) information about $f$ in terms of data observations, or we are only allowed to query $f$
to gain information. This separation between the approximation setting and the numerical setting is important since it may be that we could approximate $f$ well if we had unlimited access to $f$, but in reality we are limited by the information provided to us. 

Fully connected feed-forward NNs 
are an important approximation tool that is amenable to theoretical analysis. In practice, the most common choice of activation function 
$\sigma$ is the so-called rectified linear unit
\[\sigma(t)=\Relu(t):=t_+:=\max\set{0,t}.\]
This will constitute the main example of activation function studied in this article.

\subsection{ Fundamental Operations  with Neural Networks}
 \label{SS:fund}
In this section, we discuss  some   fundamental operations that one can implement with NNs.   Recall that $\Upsilon^{W,L}=\Upsilon^{W,L}( \sigma; d,d')$  is the set of functions that are outputs of a NN with the activation function $\sigma$, input dimension $d$, output dimension $d'$, and  $L$ hidden layers each  of fixed width $W$. 
 
Let us begin by pointing out that deep neural networks naturally allow for two fundamental operations -- parallelization and concatenation -- which we will often use. 
 \vskip .1in
  \noindent   
{\bf Parallelization:} {\it If the NNs $\cN_j$ 
 have width $W_j$, depth $L$, input dimension $d$, output dimension $d'$, and an activation function $\sigma_j$,  $j=1,\dots, m$, then the parallelization of these networks is a new network ${\rm PAR}(\cN_1, \ldots,\cN_m)$ with width $W=W_1+\ldots+W_m$, depth $L$, input dimension $d$ and output dimension $d'$. Its graph is obtained by placing the hidden layers of $\cN_j$  on the top of each other.   The parallelized  network can output any  linear combination $S=\sum_{j=1}^m\alpha_jS_j$, where $S_j\in \Upsilon^{W_j,L}(\sigma_j;d,d')$, $j=1,\dots,m$. }

 As described above, the network 
 ${\rm PAR}(S_{\cN_1},\ldots,S_{\cN_m})$ does not have full connectivity
 since the nodes of $\cN_j$ are not connected to the nodes of $\cN_i$, $i\neq j$. However, we can view the resulting network as a fully connected network by completing it, that is, by adding the missing edges and assigning to them zero weights.

To describe network concatenation, let us agree that, given $m$ functions $f_j:\R^{d_j}\to \R^{d_{j+1}}$, we will define their composition  $f_m\circ \cdots \circ f_1:\R^{d_1}\gives \R^{d_{m+1}}$ by \be
 \label{composition}
\lr{f_m\circ \cdots \circ f_1}(x):=f_m\lr{f_{m-1}\lr{\cdots f_1(x)}}.
 \ee
If $f:\R^d\to \R^d$, we also introduce the notation
 \be 
 \label{nfold} 
 f^{\circ m}:=f\circ f\circ \cdots \circ f,
 \ee 
where the composition is performed $m-1$ times. 
\vskip .1in
 \noindent
 {\bf Concatenation:} {\it If the NNs $\cN_j$ have width $W_0$, depth $L_j$, input dimension $d_j$, output dimension $d_{j+1}$, and activation functions $\sigma_j$, $j=1,\dots, m$, then the concatenation of these networks is a network ${\rm CONC}(\cN_1, \ldots,\cN_m)$ with width $W_0$, depth $L=\sum_{j=1}^mL_j$, input dimension $d_1$ and output dimension $d_{m+1}$. Its graph is obtained by placing the hidden layers of these networks  side by side  with full connectivity between the hidden layers of $\cN_j$ and $\cN_{j+1}$. The concatenated  NN can output any  composition  $S=S_m\circ S_{m-1}\circ \cdots \circ S_1$,  where the functions $S_j\in \Upsilon^{W_0,L_j}(\sigma_j;d_j,d_{j+1})$, $j=1,\ldots,m$.  It does this by assigning   weights and biases, associated to edges connecting the last hidden layer of
 an $\cN_j$ to a node of   the first hidden layer of the neighbor $\cN_{j+1}$, using the output weights and biases of $\cN_j$ and input weights and biases of $\cN_{j+1}$.}
 
 
 Parallelization and concatenation of neural networks allow us to perform the following operations between their outputs.
  \vskip .1in
  \noindent
 {\bf Addition by increasing width:} It   follows from 
 {\bf Parallelization} that for any $L\ge 1$ and 
 $S_j\in \Upsilon^{W_j,L}(\sigma;d,d')$, $j=1,\dots,m$, the linear combination 
 $$ 
 \sum_{j=1}^m \alpha_j S_j\in\Upsilon^{W,L}(\sigma;d,d'), \quad W:=W_1+\cdots +W_m.
 $$
 
 \vskip .1in
 \noindent
 {\bf Composition:} It follows from {\bf Concatenation} that for any $W\ge 1$ and $S\in\Upsilon^{W,L_1}(\sigma;m,d')$, $T\in \Upsilon^{W,L_2}(\sigma;d,m)$, the
 composition
 \begin{equation}
 \nonumber
     S\circ T\in \Upsilon^{W,L}(\sigma; d,d'),\quad  L:=L_1+L_2.
 \end{equation}
 \vskip .1in
 
  \begin{remark}
  \label{R:furtherconcat}
  In particular, 
 it follows from {\bf Parallelization}  and 
 {\bf Concatenation} that given the outputs 
 $T_j\in \Upsilon^{W_j,L_1}(\sigma; d,1)$,
 $j=1,\ldots,m$, and   $S\in\Upsilon^{W,L_2}(\sigma;m,d')$, 
 with $W=\sum_{j=1}^mW_j$,
 then 
 the function 
 $$
 S(T_1,\ldots,T_m)\in\Upsilon^{W,L}(\sigma; d,d'), \quad  L:=L_1+L_2.
 $$
 \end{remark}

 \vskip .1in
 \noindent
 {\bf Shifted dilates:}  {\it If $S\in\Upsilon^{W,L}(\sigma; d,d')$, then for any $a\in\R$ and $c\in\R^d$, the shifted dilate 
 $$
 T(x):=S(ax+c)\in \Upsilon^{W,L}(\sigma; d,d').
 $$
 }
 \vskip .1in
 \noindent 
 To prove this, let $\cN$ be the NN which outputs the function $S$.  To output the function $T$, it is enough to alter the 
 weights and biases of the first hidden layer of $\cN$.  Namely, if a neuron
 from this layer computes $\sigma (w \cdot x+b)$, $w\in\R^d$, $b\in\R$,  we replace it with 
 $\sigma (aw\cdot x+ 
 w\cdot c+b)$. Here,  $x\cdot x'$ denotes the inner product of two vectors $x,x'$ of the same dimension.   The remaining layers stand the same.

 \subsection{One Layer Neural Networks}
\label{S:one-layer-meas}
The  function $S_\cN\in \Upsilon^{W,1}(\sigma; d,1)$ produced by a single hidden  layer,   fully connected feed-forward neural network $\mathcal N$ with activation function $\sigma$,    $d$ inputs and one output  has the representation
\be 
\label{onelayerrep}
S_\cN(x)= b_0+\sum_{j=1}^{W} a_j\sigma(  w_j\cdot x +b_j),\qquad a_j,b_j\in \R,\quad w_j\in \R^{d},
\ee
where $W$ is the width of the first (and only) hidden layer and 
 $b_0$ is the bias of the  output node.  
The above  can equivalently be written as
\begin{eqnarray}
\nonumber
S_\cN(x) = b_0+\int_{\R^{d+1}} \sigma(w\cdot x+b)d\mu_\cN(w,b),  \quad \mu_\cN(w,b):=\sum_{j=1}^{W} a_j \delta_{(w_j,b_j)}, 
\end{eqnarray}
where $\delta_z$ denotes the mass one atomic measure at the point $z$. The correspondence between purely atomic Borel measures and  $\Upsilon^{W,1}(\sigma;d,1)$ is useful in addressing various structural properties of this set via functional analytic arguments.
As an example of this, we discuss briefly the density question of whether for each continuous function $f$, defined on a compact set $\Omega\subset \R^d$, we have
\be
\label{density}
{\rm dist}(f,\Upsilon^{W,1}(\sigma;d,1))_{C(\Omega)}\to 0, 
\quad W\to\infty,
\ee
where for the current discussion  the distance is measured in the uniform norm $\|f\|_{C(\Omega)}:=\sup_{x\in\Omega}|f(x)|$.
This question is discussed in detail in \cite{pinkus1999approximation}, see also \cite{cybenko1989approximation}, \cite{petersennotes}.  Here we only point out some key results.

Note that \eref{density} does not hold for every  activation 
function $\sigma$. For example, if $\sigma=P$ is a univariate polynomial of degree $m$, then $\sigma(w\cdot x+b)$ with $w\in\R^d$, $ b\in\R$ is a multivariate polynomial in $x=(x_1,\dots,x_d)$
of  total  degree $m$   and hence   
$\Upsilon^{W,1}(\sigma;d,1)\subset X$, where $X$ is a linear space of fixed finite  dimension.  
Thus, \eref{density} does not hold.

\vskip .1in
\noindent
{\bf A sufficient condition:}  If $\sigma$ is a continuous function on $\R$ such that for each  finite, signed regular
Borel measure $\mu\neq 0$ on $\Omega$ the function
\be
\nonumber
F_\mu(w,b):=\int_\Omega\sigma(w\cdot x+b)\, d\mu(x),\quad w\in\R^d,\, b\in\R,
\ee
is not identically zero, then the density condition \eref{density} holds. This condition can be used to prove the following examples of activation functions for which  \eref{density} holds.
\begin{itemize}
\item {\bf Sigmoidal activation function:} A   function $\sigma$, defined and continuous on $\R$,  is called a sigmoidal function if 
$$
\lim_{t\to \infty}\sigma(t)=1, \quad \hbox{and} \quad 
\lim_{t\to -\infty}\sigma(t)=0.
$$
For each such $\sigma$ the density statement \eref{density} holds, see \cite{cybenko1989approximation} for one of the first proofs in this case.
\item {\bf ReLU activation function:}  If $\sigma(t)=t_+$, $t\in\R$, then the density condition \eref{density} holds.
\end{itemize}

\section{ReLU Networks}
\label{S:Relunetworks}
In this section, we summarize what is known about the outputs of 
NNs with ReLU activation  (ReLU networks).  We begin by making general remarks that hold for any ReLU network
and then turn to special cases, especially those that form our main interest of study in this paper.

Perhaps the most important structural property of ReLU networks is that any output of such a network is a continuous piecewise linear  function. To describe this precisely, we start with the following definitions.
\begin{definition}\label{D:cpwl1}
A \textit{polytope partition} of $\R^d$ is a finite collection $\cP=\set{P_j}$ of convex closed  $d$ dimensional polytopes (not necessarily bounded) which are exhaustive and have disjoint interiors $P_j^o$, that is,
\[\bigcup_{j} P_j = \R^d,\qquad P_j^{\circ}\cap P_k^{\circ}=\emptyset,~~~\forall j\neq k.\]
Each such convex polytope is the intersection of a finite number of closed half spaces. We refer to the polytopes $P_j$ of such a partition as {\it cells}.
\end{definition}
\begin{definition}\label{D:cpwl}
A function $S:\R^d\gives \R$ is a continuous piecewise linear function (CPwL) if $S$ is globally continuous and there is a  polytope partition $\cP=\set{P_j}$ on which $S$ is  locally affine, that is,
\[S|_{P_j}\quad \text{ is affine }~~~~\hbox{for \ all}\,\,\  j.\]
We then say that $S$ is subordinate to the polytope partition $\cP.$
\end{definition}
\vspace{.2cm}We denote by 
$$
\Sigma_{n,d}:=\Sigma_{n,d}({\rm CPwL})
$$
the collection of all CPwL functions $S:\R^d\gives \R$ that are subordinate to some polytope  partition with at most $n$ cells. This collection is
a nonlinear set. For example, if $S_1$ and $S_2$ are subordinate to different partitions of size $n$ then the sum $S_1+S_2$ is typically  not in $\Sigma_{n,d}$.
Going forward in this paper, we do not study $\Sigma_{n,d}$ but only use it for comparison purposes. Note that if a CPwL function $S$ is subordinate to $\cP$ then it is also subordinate to any refinement of $\cP.$

In the special case when $d=1$, polytope partitions of $\R$ are simply decompositions of $\R$ into  intervals with disjoint interiors, and thus
$\Sigma_{n,1}$
is in fact the set of all univariate continuous linear free-knot splines with at most $n-1$ break points.

\begin{theorem}
\label{T:CPwL}
Let $\cN$ be a ReLU network with $d$ inputs,   one output node, and $m$ hidden neurons.  Then, the output $S_\cN$ of $\cN$ is a CPwL
function subordinate to  a partition $\cP_\cN$ with at most $3^m$ cells, i.e., $\# \cP_\cN\le 3^m$.
\end{theorem}
\noindent
{\bf Proof:} 
 Let us denote by 
$z_1(x),\ldots, z_{m}(x)$ the pre-activations of the network's neurons, that is the values stored at the neurons for input $x$ before ReLU is applied. For every \textit{activation pattern}
\[\nu=\lr{\nu_1,\ldots, \nu_{m}}\in \set{-1,0,1}^{m},\]
we define
\begin{equation}
\label{E:omega-def}
\Omega_\nu :=\set{x\in \R^d~:~ \text{sgn} (z_j(x))=\nu_j,\ j=1,\dots,m},    
\end{equation}
where for the purpose of this formula $\text{sgn} (0):=0$. By construction, each $\Omega_\nu$ is the (possibly empty) collection of all inputs $x\in\R^d$ at which the network  neurons have a given pattern of being on, off, or zero, prescribed by $\nu.$ A simple inductive argument, see Lemma 7 in \cite{hanin2019deep},  shows that $\Omega_\nu$ is a convex polytope. Moreover, defining $P_\nu$ to be the closure of $\Omega_\nu$, we see that the collection
\[\cP_\cN:= \set{P_\nu~|~ P_\nu^o\neq\emptyset}\]
is a polytope partition of $\R^d$ and that $S_{\mathcal N}$ is a CPwL function subordinate to this partition.\hfill $\Box$

Having established that each output of a ReLU network is a CPwL function, it is of interest to give bounds for the number of cells in such a partition.
The above theorem gives a bound $3^m$. However, many of the cells $\Omega_\nu$, defined in \eref{E:omega-def}, are either empty or have dimension smaller than $d$.
We shall see as we proceed in this section, that this bound can be improved in the cases of interest to us. At this stage, let us just mention the following almost trivial result. 
  \vskip .1in
  \noindent   
{\bf Claim.} {\it Consider a fixed architecture for neural networks with $m$ neurons as above. Let $S(\cdot;\theta)$
be the outputs of a ReLU network with  parameters $\theta\in \R^n$. Then, outside a set of measure zero in $\R^n$, any selection of parameters results in an $S(\cdot,\theta)$ which is subordinate to a partition with at most $2^m$ cells.}
\vskip .1in

 This claim is proved by showing that outside of a set of measure zero in parameter space $\R^n$, all cells $\Omega_\nu$ defined in \eref{E:omega-def}
 are empty or have dimension $<d$ whenever one of the components $\nu_j$ of $\nu$ is zero.

Our purpose in the remainder of this section is to explore the properties of both the polytope partitions created by ReLU networks and the complexity of the CPwL functions that they output. We  start in \S \ref{S:d=1} by studying in detail ReLU networks with input and output dimension $1$, postponing a discussion of higher input dimensions to \S \ref{S:d>1}.

\subsection{Univariate ReLU Networks}
\label{S:d=1}
In this section, we  consider ReLU networks with input and output dimensions both equal to one, that is, $d=d'=1$.  In this case, the polytope partitions of $\R$ are simply decompositions of $\R$ into a finite collection of intervals with disjoint interiors, and the CPwL functions subordinate to such partitions are customarily referred to as continuous linear free-knot splines. 

\subsubsection{Single Layer Univariate ReLU Networks} 
\label{SSS:shallowd1}
For the  set $\Upsilon^{W,1}:=\Upsilon^{W,1}({\rm ReLU};1,1)$,
 we have the simple inclusion, see \cite{daubechies2019nonlinear},
\begin{equation}
    \label{eq}
\Sigma_{W,1}\subsetneq\Upsilon^{W,1}\subsetneq \Sigma_{W+1,1}, 
\end{equation}
where we recall our notation  $\Sigma_{W,1}:=\Sigma_{W,1}({\rm CPwL})$ for the set of CPwL functions subordinate to a partition of $\R$ into   $W$ intervals.
This shows that $\Upsilon^{W,1}$ and  the set of linear free-knot splines, determined by comparable number of parameters, essentially have the same approximation power.
Recall that, according to \eref{numberparameters}, 
$\Upsilon^{W,1}$ is deterined by $3W+1$
parameters, while $\Sigma_{W,1}$ is determined by $2W$ parameters.

We point out a particularly important  family of functions generated by ReLU networks,  namely,  the hat functions
$H_{\bf p}$, where ${\bf p}=(p_1,p_2,p_3)\in \R^3$, $p_1<p_2<p_3$, defined as
\begin{eqnarray}
\label{hatfunction}
H_{\bf p}(t)=\begin{cases}
0, \quad \quad t\notin [p_1,p_3],\\\\
\frac{t-p_1}{p_2-p_1}, \quad t\in [p_1,p_2],\\\\
-\frac{t-p_3}{p_3-p_2}, \quad t\in [p_2,p_3].
\end{cases}
\end{eqnarray}
$H_{\bf p}$ is a CPwL function  that takes the value one at $p_2$, zero at $p_1$ and $p_3$, is linear on $[p_1,p_2]$ and $[p_2,p_3]$, and vanishes outside of $[p_1,p_3]$.  
Note that since $H_{\bf p}\equiv 0$ outside $[p_1,p_3]$, we have
\begin{eqnarray}
\nonumber
H_{\bf p}(t)=\frac{1}{p_2-p_1}(t-p_1)_+
-\frac{p_3-p_1}{(p_3-p_2)(p_2-p_1)}(t-p_2)_+
+\frac{1}{p_3-p_2}(t-p_3)_+,
\end{eqnarray}
and hence $H_{\bf p}\in \Upsilon^{3,1}({\rm ReLU};1,1)$. In particular, the  hat function $H$, defined as $H:=H_{(0,1/2,1)}$ and viewed as a function on $[0,1]$ has the representation
\begin{equation}
\label{hat}
H(t)=2(t-0)_+-4(t-\frac{1}{2})_+.
\end{equation}
Thus,  $H\in\Upsilon^{2,1}({\rm ReLU};1,1)$ when considered only on $[0,1]$.

\subsubsection{Deep Univariate ReLU Networks} 
\label{SS:deepunirelu}
According to the discussion at the start of \S \ref{S:Relunetworks}, any function  from the set  $\Upsilon^{W,L}:=$  
$\Upsilon^{W,L}(\Relu;1,1)$ is a CPwL function
on $\R$.  It is of interest to understand  exactly which  CPwL functions   are  in this set.  We shall see that such a characterization is rather straightforward  when $L=1$, but the situation gets more complicated as $L$ gets larger.

When $L=1$, any  selection  of weights and biases produces as output a CPwL function $S$ with at most $W$  breakpoints.  Indeed, 
$S$ can be expressed as $S=b_0+\sum_{j=1}^W a_j \eta_j(t)$, where the functions  $\eta_j(t)=(\pm t+b_j)_+$. Obviously the bound $W$ cannot be improved.  Although $\Upsilon^{W,1}$ does not contain all of $\Sigma_{W+1,1}$, it does contain all of $\Sigma_{W,1}$, see \eref{eq}.

When $L>1$, the situation gets much more complicated. Even though there  is no precise characterization of the set of outputs, we can provide some important insight.     When $L$ grows, two important things happen: 
\begin{itemize}
    \item [(i)]
the number breakpoints of functions from  $\Upsilon^{W,L}$ can be exponential in $L$;  
 \item[(ii)] not every CPwL function with this large number of breakpoints is in $\Upsilon^{W,L}$, in fact,  far from it.
 \end{itemize}
We first address (i).  Fix $W$ and let $S\in \Upsilon^{W,L}=\Upsilon^{W,L}(\Relu;1,1)$.  We define $m(L)$ as the
  largest number of breakpoints that any $S\in \Upsilon^{W,L}$
  can have.  We know that $m(1)=W$. Moreover, once the parameters are chosen for the first layer, any output $S$ has breakpoints in a fixed set $\Lambda$ of cardinality at most $W$.   We can bound $m(L+1)$ in terms of 
  $m(L)$ as follows.  Each $S\in \Upsilon^{W,L+1}$ can be expressed as
\be 
\label{written}
S=\sum_{k=1}^W a_k [S_k]_+  +b,
\quad a_k,b\in \R, \quad S_k\in 
\Upsilon^{W,L}, \quad k=1,\ldots,W.
\ee 
There is a  set $\Lambda$, 
$\#(\Lambda)\leq m(L)$ 
such that each of the $S_k$ have their breakpoints in $\Lambda$.
Fix $k$ and consider the  function $[S_k]_+$.  It has two types of breakpoints. One  are those it inherited from $\Lambda$ and the second
is the set $\Lambda_k'$ of new breakpoints that arose after the application of ReLU.  We have $\#(\Lambda_k')\le \#(\Lambda)+1\le m(L)+1$, $k=1,\dots,W$.
Hence, $S$ has at most $m(L)+W(m(L)+1)$ breakpoints.It follows that
\be 
\label{countbps}
m(L+1)\le (W+1)m(L)+W,\quad L\ge 1.
\ee
This recursion with the starting value $m(1)=W$ gives the  bound
\be 
\label{countbps1}
m(L)\le (W+1)^L,\quad L=1,2,\dots.
\ee
This bound can be improved somewhat at the expense of a more involved argument.

This potential exponential growth of the breakpoints as a function of the number of neurons can in fact be attained. A simple example, first noted by Telgarsky \cite{telgarsky2015representation}, is to compose the hat function $H_{(0,1/2,1)}$ on $[0,1]$, see \eref{hat}, with itself $(L-1)$ times.  The resulting function $S:= H^{\circ L}$ is the  saw tooth function with $2^{L-1}$ teeth, see Figure \ref{F16}. Since $H\in \Upsilon^{2,1}(\Relu;1,1)$, it follows from {\bf Composition} that 
$S\in \Upsilon^{2,L}(\Relu;1,1)$.
\begin{figure}[h]
  \centering
\includegraphics[scale=.5]{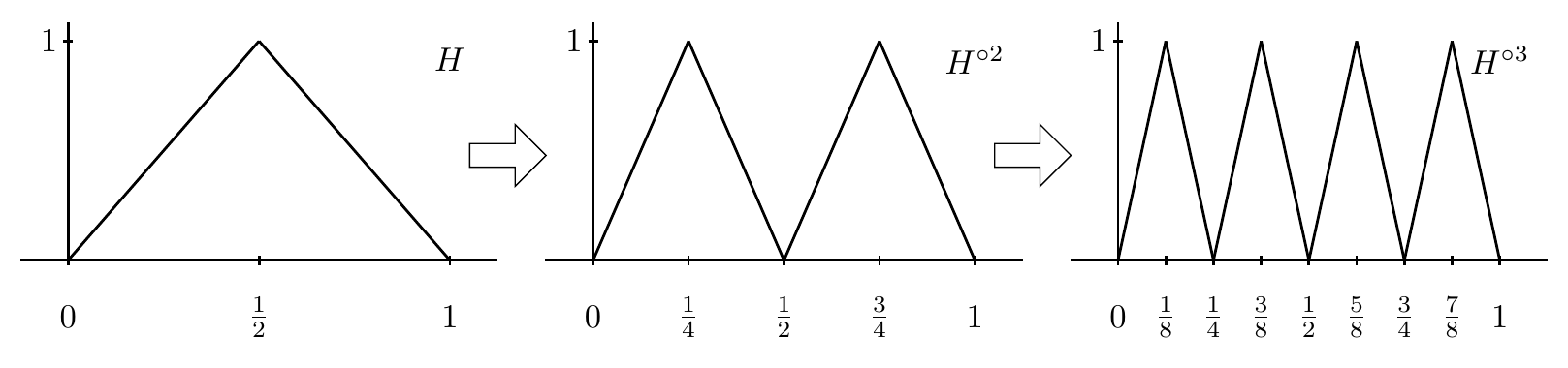}
\caption{The Sawtooth Functions $H^{\circ L}$.}
\label{F16}
\end{figure}

While a function in $\Upsilon^{W,L}(\Relu;1,1)$ can have an exponential (in $L$) number $N_L$ of breakpoints, one should not get disillusioned into thinking
that this set of functions is anywhere close to $\Sigma_{N_L,1}$, 
which was the point in  (ii).  The reason for this is that
 there are linear dependencies between the linear pieces in the case of a large number of breakpoints.  
 There is not yet a  good understanding   of exactly which CPwL functions are in $\Upsilon^{W,L}(\Relu;1,1)$ when $L$ is large.  A possible starting point to unravel this is  to consider special cases such as breakpoints in $[0,1]$ at the  dyadic integers $j2^{-n}$, $j=0,1,\dots,2^n$, with $n>1$.   
 \vskip .1in
 \noindent
 {\bf Problem 1:} {\it For $L>n $, characterize the CPwL functions in $\Upsilon^{W,L}(\Relu;d,1)$ which have breakpoints only at the
 dyadic integers $j2^{-n}$, $j=0,\dots,2^n$.}

\subsection{Multivariate ReLU Networks}
\label{S:d>1}
We now turn to studying the properties of ReLU networks with input dimension $d>1$, starting with those networks that have one hidden layer. Deeper multivariate ReLU networks are discussed in \S \ref{S:d>1-L>1}.

\subsubsection{Multivariate ReLU Networks with one hidden layer}
\label{S:d>1-L=1}
A ReLU network $\mN$ with input dimension $d>1$, output dimension $1,$ and one hidden layer of width $W$ outputs a function of the form
\begin{equation}\label{onelayernet2}
    S_\cN(x) = b_0+\sum_{j=1}^{W} a_j\eta_j(x),\quad \eta_j(x):=(z_j(x))_+,\qquad  x\in \R^d,
\end{equation}
where $b_0,a_j\in\R$, $j=1,\dots,W$,  and   $z_j$, $j=1,\dots,W$, is the function computed by the $j^{th}$ neuron before applying ReLU,
\begin{equation}
\nonumber
    z_j(x)=z_j(x;w_j,b_j)=w_j\cdot x+b_j,\qquad w_j\in \R^d,\, b_j\in \R.
\end{equation}
 Let us record the following useful fact.
 \vskip .1in
 \noindent
 {\bf Observation:} {\it Any function $S\in \Upsilon^{W,1}(\Relu;d,1)$ has a representation \eref{onelayernet2} with the $w_j$, $j=1,\dots,W$, unit norm vectors.} 
 \vskip .1in
 \noindent
This follows by removing the zero weight vectors and normalizing the remaining $w_j$'s by  adjusting the constants $b_j$ and $a_j$. We use this representation of functions in $\Upsilon^{W,1}(\Relu;d,1)$ in going forward.   
Given the collection of $W$ unit vectors $w_j\in\R^d$ and biases $b_j\in\R$
from \eref{onelayernet2}, we define the hyperplanes
\[H_j:=\set{x\in \R^{d}: w_j\cdot x+b_j=0},\quad j=1,\dots,W,\]
and the  collection ${\cal H}:=\set{H_1,\ldots, H_W}$, associated to the network $\cN$.  This collection is an example of a \textit{hyperplane arrangement}, a classical subject in combinatorics \cite{stanley2004introduction}. 
\begin{figure}[h]
  \centering
\includegraphics[scale=.25]{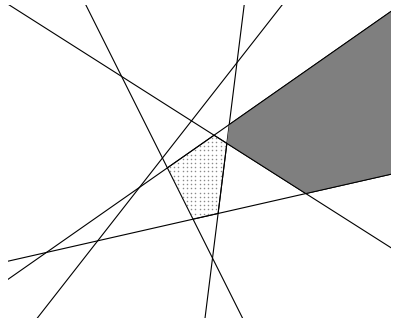}
\caption{A hyperplane arrangement in $\R^2$ with two of its cells shaded.}
\label{hyperplanes}
\end{figure}

We now describe how the hyperplane arrangement $\cal H$ associated to $\cN$ determines the polytope partition of $\R^d$ to which $S_{\mN}$ is subordinate in the sense of Definition \ref{D:cpwl}. Because of our assumption that each $w_j$ has unit norm, only activation patterns with entries $\nu_j\in \set{\pm 1}$ lead to cells with nonempty interiors. For any such activation pattern $\nu=\lr{\nu_1,\ldots, \nu_W}$, we may write in the notation of \eqref{E:omega-def},
\[\Omega_\nu = \bigcap_{j=1}^W H_j^{\nu_j},\qquad H_j^{\nu_j}=\set{x\in \R^d~:~\text{sgn}(z_j(x))=\nu_j},\]
as an intersection of half-spaces.  Thus, in the special case of $\Upsilon^{W,1}$ the cells in partitions $\cP$ of the output functions have a simple global description as the connected components of $\R^d$ when the hyperplanes are removed. The partition $\cP$ associated to $\cal H$ is the collection of the closures $P_\nu:=\overline{ \Omega}_\nu$ for which   $\Omega_\nu\neq \emptyset$.    Each cell $P_\nu$ is a convex closed polytope. 
By a special case of a classical result of Zaslavsky \cite{Zaslavsky1975}, the number of cells $\# \cP$ in $\cP$ satisfies
\begin{equation} 
\label{countcells}
    \#\cP\leq \sum_{j=0}^d \binom{W}{j}.
\end{equation}
In fact, Zaslavsky's theorem shows that away from a co-dimension $1$ set of weights and biases (i.e. when the hyperplanes are in general position), this upper bound is attained. 

In summary, any function $S\in\Upsilon^{W,1}(\Relu;d,1)$ is a CPwL function subordinate to a partition $\cP$,  generated by an  arrangement of $W$ hyperplanes determined by the weights and biases of the network $\cN$. However, it is imporant to note that, unlike the case of one hidden layer with input dimension $1,$ not every CPwL functions subordinate to a polytope partition arising from a hyperplane arrangement is the output of a    one layer ReLU network.   There are several ways to see this that we now discuss.

First of all, let us show that  $\Upsilon^{W,1}(\Relu;d,1)$ does not contain any compactly supported function on $\R^d$  once $d>1$. To see this, consider a function $S$ of the form \eref{onelayernet2} and suppose that $S$ has a compact support on $\R^d$. For each $j=1,\ldots,W$, there is a ball $B_j\subset\R^d$ outside of the support of $S$ that intersects the hyperplane $H_j$ but none of the other hyperplanes $H_i$, $i\neq j$.  We have, 
$$
0=S(x)=L(x)+ a_j\eta_j(x), \quad x\in B_j,
$$
where $L$ is an affine function. Since $L$ and $\eta_j$ are linearly independent on $B_j$, this  implies $a_j=0$.  Hence, all $a_j$ are zero and $S$ is a constant. Since $S$ was assumed to have compact support, this constant is zero.

Another way to see that $\Upsilon^{W,1}(\Relu;d,1)$ does not contain all CPwL functions subordinate to a given hyperplane arrangement is the following. Once ${\cal H}=\{H_1,\dots,H_W\}$ is chosen, thereby determining the partition $\cP$, the outputs of $\Upsilon^{W,1}$ that are subordinate to $\cP$ are all contained in a linear space of dimension $2W+1$. This follows from the representation \eref{onelayernet2}. Indeed, since we've assumed that $\norm{w_j}=1$ for each $j$, there are  only two choices of $(w_j,b_j)$ for the functions $z_j(x)=w_j\cdot x + b_j$, 
$j=1,\ldots,W$, such that
\[H_j=\set{x\in \R^d~:~z_j(x)=0}.\]
However, when $d>1$, Zaslavsky's theorem  shows that the number of cells in $\cP$ can grow 
as fast as $CW^d$ when $W\geq d$. Hence, in general, the set of all CPwL functions subordinate to $\cP$ is a linear space with dimension much larger  than $2W+1$.

The following lemma gives a simple way of checking when a CPwL function subordinate to a partition generated by an arrangement of $W$ hyperplanes is in   $\Upsilon^{W,1}(\Relu;d,1)$. 
Before formulating the lemma, let us  note that if $T$ is a CPwl function subordinate to $\cP$, then  on any cell $\cP_\nu$ of $\cP$, the gradient $\nabla T$   is a constant vector. 
It follows that
$\nabla T$ is a piecewise constant vector valued function subordinate to $\cP$.

\begin{lemma}
\label{L:one-layer-char}
Let $\cP$ be a partition of $\R^d$ generated by  a  hyperplane arrangement 
${\cal H}=\set{H_1,\ldots, H_W}$,  where $H_j:=\{x\in \R^d\,\,:\,\,w_j\cdot x+b_j=0\}$, and $w_j$ is a unit vector, $j=1,\ldots,W$. Let $T$ be CPwL function that is subordinate to $\cP$. Then $T$ has the representation 
$$
T=S+L, \quad S\in \Upsilon^{W,1}(\Relu;d,1), \quad L\,\,\hbox{- globally affine},
$$
if and only if  the following condition holds:
 \begin{enumerate}
    \item[(\textbf{A})] For each $j=1,\ldots,W$, there is a real number $a_j$ such that for  every $x\in \R^d$ on the   hyperplane $H_j$,
 and on no other hyperplane,
    the jump  in $ \nabla T$   across $H_j$ at $x$,    equals $a_jw_j$.
\end{enumerate}
\end{lemma}
\noindent
{\bf Proof:} First, let  $T=S+L$ with $S$  and $L$ as above.  We know that the function 
$S\in \Upsilon^{W,1}(\Relu;d,1)$ has the representation \eref{onelayernet2} with the $w_j$'s being unit vectors. Given $x\in \R^d$ that belongs to $H_j$ and to no other hyperplane, the jump in 
$\nabla T$ at $x$ is the same as that of $a_j\nabla \eta_j$ 
at $x$, which is $a_jw_j$. This shows that $T$ satisfies ({\bf A}). 

For the converse, suppose that $T$ is any CPwL that is subordinate to $\cP$ and  that $T$ satisfies condition ({\bf A}).   We define $S:=\sum_{j=1}^W a_j\eta_j\in \Upsilon^{W,1}(\Relu;d,1)$,
where the $a_j$ are given by {\bf (A)} and $\eta_j(x):=(w_j\cdot x+b_j)_+$. Consider the function $(T-S)$ which is piecewise linear subordinate to the partition $\cP$.  We claim that this function is a globally affine function. Indeed,  otherwise there would be two adjacent cells which share a $d-1$ dimensional boundary (which is part of some $H_j$) and the jump of $\nabla(T-S)$ across this boundary is not zero.  But both $T$ and $S$ have the same jump $a_jw_j$ of their gradient across this boundary. This is a contradiction and proves the lemma.
\hfill $\Box$

\subsubsection{Deep Multivariate ReLU Networks}
\label{S:d>1-L>1}
The discussion at the beginning of \S \ref{S:Relunetworks} showed that any  $S\in \Upsilon^{W,L}(\Relu;d,1)$  is a CPwL function subordinate to a partition $\cP$  of $\R^d$ into convex  polytopes. The partition we produced to show this was not  determined by a  hyperplane arrangement.  It turns out, as we shall see in
\S\ref{SS:GenCPwL}, that $S$ is always subordinate to some partition given by a  hyperplane arrangement.  However, the latter partition is not
a minimal partition to which $S$ is subordinate.  In other words, unlike the case of $L=1$, the minimal polytope partition of $\R^d$ to which $S$ is subordinate is not simply given by the cells of a hyperplane arrangement.

Given $S\in \Upsilon^{W,L}(\Relu;d,1)$, we do not know the best bound for the number of cells in a minimal convex polytope partition to which $S$ is subordinate, but we can give some bounds.  Recall that from   \eqref{E:omega-def} we have the bound $3^{WL}$.   We also stated that in the generic case this bound can be improved to $2^{WL}$. Indeed, in the generic case,  this partition consists of the closures of those convex sets
\[
\Omega_\nu=\set{x\in \R^d~:~ \text{sgn}(z_j(x))=\nu_j},\qquad \nu = \lr{\nu_1,\ldots, \nu_{WL}}\in \set{\pm 1}^{WL},
\]
which have a nonempty interior.
We continue to write $z_j(x)$ for the CPwL function computed by the $j^{th}$ neuron in $\mN$ before ReLU is applied, and we have assumed for simplicity that for every neuron $z_j$ the sets
\[H_{j}=\set{x\in \R^d~:~ z_j(x)=0}\]
have co-dimension at least $1$.  It is important to note that the $H_j$'s are no longer hyperplanes since the functions $x\mapsto z_j(x)$ are not affine. Instead, $H_j$ is the zero level set of $z_j$ and, following the language in \cite{hanin2019universal}, we refer to the $H_j$ as {\it bent hyperplanes} and
$${\cal H}=:\set{H_1,\ldots, H_{WL}}$$
as a bent hyperplane arrangement. We can now describe the cells in the partition $\cP$,
\[\cP= \set{\overline{\Omega}_\nu,\, \nu\in \set{\pm 1}^{WL},\, \dim(\Omega_\nu)=d },\]
or equivalently the cells that are  the connected components of $\R^d\backslash {\cal H}$.

To understand this setting more clearly, let us consider a neuron $z$ in the second hidden layer of $\mN$. Note that the function $x\mapsto z(x)$ is the output of an element of $\Upsilon^{W,1}$. Hence, it is CPwL subordinate to the partition defined by the hyperplane arrangement \[{\cal H}^{(1)}=\set{H_{1},\ldots, H_{W}}\]
created by the neurons $z_1,\ldots, z_{W}$ in the first hidden layer of $\mN$. On each cell $\cC$ of the arrangement ${\cal H}^{(1)}$, the function $x\mapsto z(x)$ is affine.  Let $H_z$ denote the bent hyperplane associated with this neuron $z$ from the second layer. We  see that $H_z\cap \cC$ is given by the (possibly empty) intersection of a single hyperplane with $ \cC$. However, because $x\mapsto z(x)$ is a different affine function on different cells, its zero set $H_z$ may ``bend'' at the boundary between two cells and is not given globally by a single hyperplane. More is true: while in every cell $H_z$ coincides with a single hyperplane, globally, it may have several connected components. 

Just as in the case of univariate ReLU networks considered in \S \ref{S:d=1}, the number of cells defined by the bent hyperplane arrangement $H$ in a deep ReLU network with any input dimension can grow exponentially with depth. In fact, see \cite[Theorem~5]{montufar2014number},  there are ReLU networks of depth $L$ and width $W \ge d$ giving rise to partitions with   at least
\be 
\nonumber
\left\lfloor \frac{W}{d} \right\rfloor^{d(L-1)}
\sum_{j=0}^d \binom{W}{j}
\ee
cells. Similarly to the univariate case, the exponential growth in the number of pieces of the CPwL functions produced by deep networks is a consequence of composition. 
 
Let us summarize what we know about the sets $\Upsilon^{W,L}(\Relu;d,1)$.  Each $S\in \Upsilon^{W,L}(\Relu;d,1)$ is a CPwL function on a finite partition of $\R^d$ into convex polytopes.  The number of cells in the partition for $S$ can be very large compared to the number of parameters $n(W,L)$. For example, when $L=1$ the number of cells  can be of order $W^d$, and as $L$ grows the number of cells can grow exponentially with respect to $L$. Although the number of cells is large, not every CPwL function subordinate to such a partition is in $\Upsilon^{W,L}(\Relu;d,1)$ since there is linear dependency imposed
 on the affine pieces.  However, every CPwL function subordinate to a partition into convex polytopes is eventually in the $\Upsilon^{W,L}(\Relu;d,1)$
 spaces, provided we take $L$ and $W$ large enough, see {\bf CPwL1} and 
 {\bf CPwL2} in \S\ref{SS:GenCPwL}.  Let us also repeat the fact that every convex polytope is the intersection of a finite number of half spaces given by a suitable hyperplane arrangement.  Finally, by refining partitions, we have that every $S$ that is in one of the spaces $\Upsilon^{W,L}(\Relu;d,1)$ is a CPwL function on a partition given by a hyperplane arrangement but the number of these hyperplanes
 may be huge.

\subsection{Properties of deep ReLU networks}
\label{SS:deepReLUNN}
Deep ReLU networks have a variety of remarkable properties that make their outputs  a powerful approximation tool.   We describe some of these properties in the present
section.  We study the sets $\Upsilon^{W,L}:=\Upsilon^{W,L}(\Relu;d,1)$, where $W$ is generally fixed and $L$ is allowed to vary. 
We begin by introducing a special class of ReLU networks that are effective in the construction of numerical approximation methods.

\subsubsection{Special networks and their set of outputs $\Up^{W,L}$}
\label{SS:specialnetworks}
We describe in this section a set of NNs that we call {\it special networks}, following \cite{daubechies2019nonlinear}, which designate certain channels for specific tasks.  We introduce the notation
$\eta_{i,j}$ 
 for the function 
 $x\mapsto \eta_{i,j}(x)$ of the initial input $x$ at the $(i,j)^{th}$ node. Usually, $\eta_{i,j}=[z_{i,j}]_+$, 
 that is, $\eta_{i,j}$ is the CPwL function computed by the $(i,j)^{th}$ neuron after ReLU is applied, but in special networks we sometimes do not apply the activation
 $\Relu$ at certain nodes.  However, as we shall see, the outputs of a special network, when restricted to a bounded domain,  are still  
 functions in  $\Upsilon^{W,L}$.

 In a special network, we reserve the top $d$ channels to simply push forward the input values of $x$. Namely, channel $i\in\{1,\dots,d\}$ has 
\[\eta_{i,j}(x):=x_i,\qquad j=1,\dots,L,\] where $x=(x_1,\dots,x_d)$ is the initial input. This allows us to use $x$ as an input to the computation performed at any later layer of the network. We refer to such channels as {\it source channels} (SC).

 In a special network, we also designate some channels,  called {\it  collation channels} (CC), to simply aggregate the value of certain intermediate computations.  In a collation channel the ReLU activation may or may not be applied. 
 The key point here is that nodes from  both   collation  and  source channels may be  ReLU free. Therefore, such networks are not true ReLU networks. 

We denote the set of functions $S$ which are the outputs of a special network of width $W$ and depth $L$ by 
$\Up^{W,L}$.
 A useful observation made in \cite{daubechies2019nonlinear}  is that the functions that are outputs of a special network, when restricted to a bounded domain,  are in $\Upsilon^{W,L}$, that is,
 \be 
 \label{contain} 
 \Up^{W,L}\subset \Upsilon^{W,L}, \quad L\ge 1. 
 \ee
This is proved using the following observations.
 \begin{itemize}
     \item Given any configuration of weights and biases in any collation channel  of a special network and any fixed compact set $K$ of inputs, we may choose a sufficiently large value $b_{i,j}$ associated to the $(i,j)^{th}$ node  so that 
 $z_{i,j}(x)+b_{i,j}>0$ for all $x\in K$. Then, we construct the true ReLU network by assigning to this node the function $\eta'_{i,j}$, given by 
 $\eta'_{i,j}(x)=[z_{i,j}(x)+b_{i,j}]_+=z_{i,j}(x)+b_{i,j}$. The effect of $b_{i,j}$ on any subsequent computation is then eliminated  by adding an extra bias (to the bias present from the special network) for any  neuron from the next layer to which the output $\eta'_{i,j}(x)$ is passed to. We perform this procedure to every ReLU free node from the collation channels.
 \item Similar treatment as above is done for all nodes in all source channels.
 \item The ReLU network that has been constructed has the same output as that of the special network we started with.
 \end{itemize}
 This specific trick works only when $K$ is compact. Alternatively, at the expense of increasing the width, we can create a true ReLU network of width $W=W_0+2d+2k$, where $W_0+d+k$ is the width of the special network with $d$ source  and $k$ collation channels by using the identity $t=t_+-(-t)_+$. This approach works for arbitrary inputs but at the expense of increasing the width of the network.

 In what follows, we use  extensively  special networks to derive  some important properties of deep networks   since they facilitate many constructions.

\subsubsection{Some important properties of deep ReLU networks}
 \label{SS:highlighted}
 
As noted in {\bf Observation} from \S \ref{S:d>1} in the case of one layer networks, 
the vector $w$  consisting of all
incoming weights into any hidden node of a ReLU network, if nonzero, can be taken to be of Euclidean norm $\|w\|_2=1$. Indeed, this follows from the equality
$$
(w\cdot x+b)_+=\|w\|_2\left(\frac{w}{\|w\|_2}\cdot x+\frac{b}{\|w\|_2}\right)_+,
$$
and the fact that the factor $\|w\|_2$ can be absorbed 
by the outgoing weights.

Next, we return to the Addition Property.  Earlier we have shown that we can add functions in $\Upsilon^{W,L}(\sigma;d,d')$ by increasing the width  of the network
using the method of parallelization. Here, we want to observe that 
addition can also be performed by increasing depth and not significantly enlarging width. Here is a statement to that effect.

\vskip .1in
 \noindent{\bf Addition by increasing depth:}   {\it If $S_j\in \Upsilon^{W,L_j}$, $j=1,\dots,m$, then for any 
 $\alpha_j\in \R$ we have
 $$
 S:=\sum_{j=1}^m \alpha_j S_j\in \Upsilon^{W+d+1,L},
 $$
 where  $L:=L_1+\cdots +L_m$.  In this statement  all $S_j$'s are viewed as functions on $[0,1]^d$ (or any bounded rectangle ${\cal R}\subset \R^d$)}.

 \vskip .1in
 Indeed, if $\cN_1,\dots,\cN_m$, are the $\Relu$  networks  that produce the $S_j$'s, then we create from these the following special network.
 First, we augment each of the $\cN_j$'s by adding $d$ source channels, 
 and one collation channel.  
We denote this new augmented network by $\cN_j'$.   Next, we place the hidden layers of the augmented networks side by side,
connect the source channels of the $\cN_j'$'s and place (with appropriate weights) the outputs of $\cN_j'$, $j=1,\ldots, m-1$, in the collation channel. Finally, the desired sum is the  output of the concatenated network  
(with appropriate weights).
 As a result, we obtain a special network with width $W+d+1$ and depth $L=L_1+\ldots+L_m$. The result follows from the containment \eref{contain}.

 
 \vskip .1in
Another operation on the output functions of $\Relu$ networks that can be easily performed with increasing depth 
 is to take their minimum or  maximum.  For this, let us first observe that given  $t,t'\in \R$, we have
\begin{eqnarray}
\label{max1}
\max\{t,t'\}=(t-t')_++(t')_+-(-t')_+,\\ \min\{t,t'\}=(t')_+-(-t')_+-(t'-t)_+.
\nonumber
\end{eqnarray}
Hence, $\min\{t,t'\}$, $\max \{t,t'\}\in \Upsilon^{3,1}(\Relu;2,1)$. 
We can extend the above minimization to an arbitrary number of inputs.

\vskip .1in
\noindent 
{\bf Minimization/Maximization  1 (MM1):} {\it 
Let  
$$
z_j(x):=w_j\cdot x+b_j, \quad w_j\in\R^d, \quad b_j\in\R, \quad 
j=1,\ldots,m,\quad x\in \R^d,
$$
be $m$ affine functions on $\R^d$. Then,  for $W=3\cdot 2^{\lceil \log_2 m\rceil-1}$ and $L=\lceil\log_2 m\rceil$
we have
\be 
\label{mindz} 
 \min \{z_1(x),\dots,z_m(x)\}, \max \{z_1(x),\dots,z_m(x)\} \in \Upsilon^{W,L}(\Relu;d,1).
\ee
In particular, if $x_i\in \R$,
 $i=1,\ldots,m$, we have
\be 
\label{mind} 
 \min \{x_1,\dots,x_m\}, \max \{x_1,\dots,x_m\} \in \Upsilon^{W,L}(\Relu;m,1).
\ee
Moreover, 
$\Relu\left(\min \{z_1(x),\dots,z_m(x)\}\right)$ and
$\Relu\left(\max \{z_1(x),\dots,z_m(x)\}\right)$ are elements of $\Upsilon^{W,L+1}(\Relu;d,1)$.
}

\vskip .1in
 We discuss the case of minimum only, since the case of maximum is almost the same.   We start with proving \eref{mind}. In our construction, we will use the fact that
$$
\min \{x_1,\dots,x_{2^k}\}=\min_{1\leq j<2^k,\,j \,{\rm odd}}\min\{x_j,x_{j+1}\}.
$$
We first use {\bf Parallelization} to  construct for each $k\geq 1$ a neural network $\cN_k$  with $2^k$ inputs and $2^{k-1}$ outputs that 
creates a vector in $\R^{2^{k-1}}$ with components $\min\{x_j,x_{j+1}\}$, $j=1,\ldots,2^k-1$, $j$-odd, by stacking on the top of each other the networks that produce $\min\{x_j,x_{j+1}\}$.  Since each of the networks in the stack has width $3$, we end up with a network with width 
$W_k=3\cdot 2^{k-1}$ and depth $L_k=1$. We concatenate the networks 
$\cN_k, \ldots,\cN_1$ in this order, by feeding the output of $\cN_j$ as input to $\cN_{j-1}$. It is easy to see that the concatenated network $\cN^k$ outputs $\min \{x_1,\dots,x_{2^k}\}$, has depth
$L=k$ and varying widths. We can augment the network by adding extra nodes and edges to each layer so that we end up with a network of width $W=3\cdot 2^{k-1}$.

For general $m$, we let $k:= \lceil \log_2 m\rceil$ and define $\hat x_j=x_j$, $1\le j\le m$ and $\hat x_j:=x_m$, $m<j\le 2^k$.  Applying the above to this new sequence gives the result \eref{mind}. To show \eref{mindz}, we feed $z_j(x)$  into the first hidden layer of $\cN_k$  by assigning appropriate input weights and node biases.

At the end, if we want to output the ReLU of min/max, we just add another hidden layer to perform the 
ReLU.

Another way to compute the above min/max is via increasing the depth and keeping the width relatively small by 
utilizing a recursive formula, first used in \cite{hanin2019universal}.
\vskip .1in
\noindent 
{\bf Minimization/Maximization  2 (MM2):} {\it 
Let 
$$
z_j(x):=w_j\cdot x+b_j, \quad w_j\in\R^d, \quad b_j\in\R, \quad j=1,\ldots,m,\quad x\in \R^d,
$$
be $m\ge 2$ affine functions on $\R^d$. Then,  
we have
\be 
\nonumber
 \min \{z_1(x),\dots,z_m(x)\}, \max \{z_1(x),\dots,z_m(x)\} \in \Upsilon^{d+1,m-1}(\Relu;d,1).
\ee
In addition, we have that both functions 
$\Relu\left(\min \{z_1(x),\dots,z_m(x)\}\right)$ and 
 $\Relu\left(\max \{z_1(x),\dots,z_m(x)\}\right)$ are elements of  $\Upsilon^{d+1,m}(\Relu;d,1)$.
 In this statement all $z_j$'s are viewed as functions on $[0,1]^d$ (or any bounded rectangle 
 ${\cal R}\subset \R^d$).}

\vskip .1in
We discuss the case of maximum only (the case of minimum is treated likewise). 
 Let $\mu_1(x):=z_1(x)$ and $\mu_k(x):= \max\{z_1(x),\dots,z_k(x)\}$, $k\ge 2$.  We use the recursion formula
$$
\mu_k(x)= (\mu_{k-1}(x)- z_k(x))_++ z_k(x),\quad 2\le k\le m,
$$
and discuss the case ${\cal R}=[0,1]^d$. For the  case of a 
general rectangle ${\cal R}$ one 
needs to add appropriate biases.
Our construction is the following: 
\begin{itemize}
    \item the first $d$ channels of the network push forward the variables $x_1,\ldots,x_d$.
Their nodes can be viewed as ReLU nodes since $t_+=t$ for $t\geq 0$.
\item the  $(d+1)^{st}$ channel computes in its first node $(z_1(x)-z_2(x))_+$. Note that if we wanted to,  we could stop and output   $\mu_2(x)$ at this stage.  
The  
$j^{th}$ node of this channel, $j=2,\ldots,m-2$, computes
$(\mu_j(x)-z_{j+1}(x))_+$, 
which is then given as an input to the $(j+1)^{st}$ node.  The final layer $L=m-1$  will hold $\mu_{m-1}(x)$  and hence can output $\mu_m(x)$. 
\end{itemize}
To show the last statement, we 
add a hidden layer after the last hidden layer of the NN from the construction above to perform the ReLU of the 
max/min. Of course, we could augment the resulting network by adding nodes and connections so that we have a fully connected feed-forward NN.

\vskip .1in

Note that if we want to compute the min/max of $m$ linear functions $z_j$, $j=1,\ldots,m$, viewed as functions on the whole $\R^d$, we  can do this if we take  $W=2d+1$, since any channel $i$, 
$1\leq i\leq d$,   in the above construction doubles in order to be able to forward the input $t=t_+-(-t)_+$.

More general statements hold when instead  of computing the min/max of affine functions  we have to find the min/max of outputs of neural networks. 

\vskip .1in
\noindent
{\bf Minimimization/Maximization 3 (MM3):}  {\it  Let $m\ge 2$ and  let the functions $S_j\in \Upsilon^{W_j,L_0}(\Relu;d,1)$,  $j=1,\ldots, m$.   Then, 
\be
\label{minf}
S:=\min\{S_1,\dots,S_m\}\in \Upsilon^{W,L}(\Relu;d,1),
\ee
where $W:=\max\{W_1+W_2+\cdots+W_m, 3\cdot 2^{\lceil \log_2 m\rceil-1}\}$, $L=L_0+ \lceil \log_2m\rceil $. If $W_j\geq 3$, $j=1,\ldots,m$, we have $W=W_1+W_2+\cdots+W_m$.
The same statement holds for 
$\max\{S_1,\dots,S_m\}$.
}
\vskip .1in
\noindent 
We use {\bf Parallelization} to construct the first $L_0$ hidden layers of the network $\cN$ that outputs $S$ Then, from the $L_0$-th  layer
we can output any of the $S_j$, $j=1,\dots,m$.   We concatenate this with the network in {\bf MM1} which has $\lceil\log_2 m\rceil$ hidden layers to complete the construction of $\cN$. 
Clearly, the resulting network has varying width,
where the first $L_0$ layers are with width $W_1+W_2+\cdots+W_m$, while the last 
$\lceil\log_2 m\rceil$ layers have width $3\cdot 2^{\lceil \log_2 m\rceil-1}$. We 
augment this network by adding extra nodes and edges. At the end, our network has width
$W=\max\{W_1+W_2+\cdots+W_m, 3\cdot 2^{\lceil \log_2 m\rceil-1}\}.$
In the case of $W_j\geq 3$, 
$W_1+\ldots+W_m\geq 
3\cdot2^{\lceil\log_2 m\rceil-1}$, which gives that $W=W_1+\ldots+W_m$.

It is also possible to do the minimization by increasing the  depth of the network while keeping the width relatively the same.
\vskip .1in
\noindent 
{\bf Minimimization/Maximization 4 (MM4):}  {\it  Let $m\ge 2$ and   let the functions $S_j\in \Upsilon^{W_0,L_j}  (\Relu;d,1)$,    $j=1,\dots,m$.    Then,  
$S:=\min\{S_1,\dots,S_m\}$ belongs to $\Upsilon^{W,L}(\Relu;m,1)$,
where the width $W:=\max\{W_0,3\}+d+1$, and the depth $L=\sum_{j=1}^m L_j +m-1 $.
The same statement holds for the $\max\{S_1,\dots,S_m\}$.
Here, all $S_j$'s are viewed as functions on $[0,1]^d$ (or any bounded rectangle 
 ${\cal R}\subset \R^d$).}

\vskip .1in

\noindent 
In order to construct a neural network $\cN$ which shows that $S\in \Upsilon^{W,L}$, we utilize {\bf Concatenation} in place of {\bf Parallelization}  and we use special networks.  Let $\cN_j$
be a network of width $W_0$ and depth $L_j$ which outputs $S_j$, $j=1,\dots,m$.  To each of the networks $\cN_j$, we add $d$ source channels to push forward the original inputs $x_1,\dots,x_d$ and  one collation channel that we will use to update computations towards outputting $S$.  Let us denote these special
networks by $\cN_j^\prime$.   We now explain how to construct $\cN$. The first $L_1$ hidden layers of $\cN$ consist of those of $\cN_1'$. The collation
channel simply pushes forward zero for these layers.  We  concatenate $\cN_1'$ with  $\cN_2'$ by  placing $S_1$ in the collation  channel of $\cN_2'$  and then pushing it forward, and by placing the outputs of the source channels of
$\cN_1'$, multiplied by appropriate weights (those that enter the first layer of $\cN_2$), into the first hidden layer of 
$\cN_2'$.   If $m=2$, we can complete the construction by placing a last hidden layer which takes $S_1$ from the collation channel
and $S_2$ as an output from $\cN_2'$ and computes $(S_1-S_2)_+$, $(S_2)_+$, and $(-S_2)_+$. We augment the resulting network with additional nodes, if necessary, so that we have a special network with width $W$. This  network outputs $S$, has depth $L=L_1+L_2+1$, and width $W$. If $m>2$, we continue by concatenating with
$\cN_3'$. The collation channel is now occupied by $T_2:=\min\{S_1,S_2\}$. If $m=3$, then we complete as before by adding a layer to compute $(T_2-S_3)_+$, $(S_3)_+$, and $(-S_3)_+$.
Continuing  this way we  obtain the desired network.

\subsubsection{General CPwL functions}
\label{SS:GenCPwL}

We have  observed  earlier that if $\cP$ is a partition into a finite number of cells obtained from a hyperplane arrangement, then not every CPwL function subordinate to this partition is
in  the set $\Upsilon^{W,1}(\Relu;d,1)$. In particular, the latter set does not include CPwL functions with compact support.  This can be remedied by slightly increasing the depth of the network. 
More precisely, let $\Delta$ be any simplex in $\R^d$ and $x^*$ be any point in its interior. Since $\Delta$ is a convex polytope with $d+1$ facets, there exist $d+1$ affine functions $z_j:\R^d\gives \R$ such that $z_j(x^*)=1$ and
\[\Delta = \set{x\in \R^d~:~z_j(x)\geq 0,\ j=1,\dots,d+1}.\]
Thus, the {\it tent function}
\be 
\label{tent}
T:=T_{\Delta,x^*}:= \Relu\left(\min \{z_1,\dots,z_{d+1}\}\right)
\ee 
vanishes outside of $\Delta$ and satisfies $T(x^*)=1$. For example, when $d=1$ this is the hat function on $\R$. The construction {\bf MM1} ensures  the following.

\vskip .1in
\noindent{\bf Tent functions:} {\it  For each $d$ dimensional simplex $\Delta$ and each $x^*$ in its interior, the   tent function $T_{\Delta,x^*}$ is in 
$\Upsilon^{W, L}(\Relu;d,1)$
with $W=3\cdot2^{\lceil\log_2(d+1)\rceil-1}$ and $L=1+\lceil \log_2(d+1)\rceil$.}
\vskip .1in

With these remarks in hand, let us now turn to the question of whether every   CPwL function is in one of the spaces $\Upsilon^{W,L}(\Relu;d,1)$. 
 We make the following observations.
\vskip .1in
\noindent 
 {\bf CPwL1:} {\it If $S$ is a CPwL function 
 on $\R^d$,  then    
 $$
 S\in \Upsilon^{W',L}(\Relu;d,1), \quad \hbox{with} \,\,L= \lceil \log_2 (d+1)\rceil, \,\, W'-\,\,\hbox {sufficiently large}.
 $$
 }
 \noindent
 This is proved in \cite{arora2016understanding} by  using the fact 
 that any CPwL function $S$  can be written  as a linear combination of piecewise linear convex functions, each with at  at most $(d+1)$ affine pieces, that is,
 \begin{eqnarray} 
 \label{sumS} 
 S
 =\sum_{j=1}^{p} \e_j\left(\max_{i\in S_j}z_i\right),\quad \e_j\in \{-1,1\},
\quad S_j\subset \{1,2,\ldots,k\}, \quad 
 \end{eqnarray}
with $s_j:=\#(S_j)\leq d+1,$ for some affine functions $z_1,\ldots,z_k$. 

To show that $S\in \Upsilon^{W',L}(\Relu;d,1)$, we use
{\bf  MM1} to show that for every $j=1,\ldots,p$, the function $\max_{i\in S_j}z_i$ is in $\Upsilon^{W,L}({\rm ReLU};d,1)$ with
 $$
  W=3\cdot 2^{\lceil\log_2 s_j\rceil-1},\quad L=\lceil\log_2 s_j\rceil.
 $$
  For the proof, we can assume that $s_j=d+1$ for all $j$ by artificially writing  an index already in $S_j$ several times, so that we end up with networks with the same depth $L=\lceil\log_2 (d+1)\rceil$.  Using {\bf Parallelization}, we then  stack these networks to produce 
 $S\in \Upsilon^{W',L}({\rm ReLU};d,1)$ with $L=\lceil\log_2 (d+1)\rceil$ and 
 $W'=3p2^{\lceil\log_2 (d+1)\rceil-1}$.

 At the other extreme, one may wish to keep the width $W$ of the ReLU network as small as possible at the expense of letting the depth of the network grow. 
 In this direction, we have the following result.
 \vskip .1in
\noindent 
 {\bf CPwL2:} {\it If $S:{\cal R}\to\R$ is a CPwL function
 defined on a rectangle ${\cal R}\subset \R^d$, then   $S\in \Upsilon^{d+2,L}({\rm ReLU};d,1)$ for $L$
   suitably large.
   }
   
  \vskip .1in 
  To show this, we use the representation \eref{sumS} for $S$
  and utilize   {\bf MM2} to view each of the functions $\max_{i\in S_j}z_i$ as an output of a  network $\cN_j$.  
  that produces $\Upsilon^{d+1,s_j-1}(\Relu;d,1)$. 
  We concatenate the networks $\cN_j$, $j=1,\ldots,p$, by placing them next to each other and connecting their source channels. We add a collation channel where we store the
 consecutive outputs  $\e_j\left(\max_{i\in S_j}z_i\right)$
 from $\cN_j$. The resulting network computes $S$,  has
 width $W=d+2$ and depth at most $L=pd$. 
 
 A result along these lines was proven in \cite{hanin2019universal}, where it was shown that the output of any ReLU network can be generated by a sufficiently deep ReLU network with fixed width $W=d+3$.
 \vskip .1in
 \noindent

\subsubsection{Finite element spaces}
\label{SS:FEM} 

One of the most popular methods of approximation used in numerical analysis is the Finite Element Method (FEM).  This method employs certain linear spaces of  piecewise polynomials. In its simplest case, one partitions a given polyhedral domain $\Omega\subset \R^d$ into simplicial cells with some requirements on the cells to avoid hanging nodes and small angles.   The linear space $X(\cP)$  of CPwL functions subordinate to such a partition $\cP$ is used for the
approximation.  
A natural question is if and how we can use ReLU NNs  in place of 
$X(\cP)$. 

Rather than address this question in its full generality, we consider only a very  special setting that will be sufficient for our discussion of NN approximation given in later
sections of this paper.  The reader  can consult \cite{Xu+} and \cite{opschoor2019deep} for a more far reaching exposition of the relation between FEM and NNs.

For our special example,
we begin with $\Omega=[0,1]^d$, and  given  $n\ge 1$,   we consider  the uniform partition ${\cal Q}_n$ of $\Omega$ into $n^d$  cubes with sidelength $1/n$.  We denote by $V_n$ the set of vertices of the cubes in ${\cal Q}_n$.  There are $(n+1)^d$ such vertices.  Each
cube $Q\in { \cal Q}_n$ can in turn be partitioned into $d!$ simplices using the so-called Kuhn triangulation with northwest diagonal.  This gives  a partition ${\cal K}$ of $\Omega$ into
$n^d d!$ simplices.  Let $X({\cal K})$ be the space of all CPwL functions defined on $\Omega$ and subordinate to ${\cal K}$.   This is a linear space of dimension $N= (n+1)^d$.   A basis for $X({\cal K})$
is given by the nodal functions 
$\{\phi_v, \,v\in V_n\}$, which are the CPwL functions  defined on $\Omega$, subordinate to ${\cal K}$, and satisfy
\be 
\label{nodal} 
\phi_v(v')=\delta(v,v'), \quad  v,v'\in V_n,
\ee 
where $\delta$ is the usual Kronecker delta function.  Each $S\in X({\cal K})$ has the representation
\be 
\label{nodalrep}
S=\sum_{v\in V_n}S(v)\phi_v.
\ee

\vskip .1in
\noindent
{\bf FEM spaces:}    Let  $X({\cal K})$  be the finite element space  in $d$ dimensions described above, 
and let $d^*:=(d+1)!$. Then the following holds:
$$
X({\cal K})\subset \Upsilon^{W,L}(\Relu;d,1), \quad W=3(n+1)^d2^{\lceil\log_2 d^*\rceil-1}, \,\, L:=1+\lceil\log_2d^*\rceil,
$$
$$
X({\cal K})\subset
\Upsilon^{d+2,L'}(\Relu;d,1), \,\, 
L'=(n+1)^dd^*.
$$
\vskip .1in

 To prove these  statements,  we first observe that each nodal basis function $\phi_v$ can be expressed as
 \be 
 \label{expressphi} 
 \phi_v=\Relu\left(\min \{ z_\Delta: \Delta\in D_v\}\right ),
 \ee 
 where $D_v$ is the set of simplices in ${\cal K}$  that have $v$ as one of their vertices.  There are $(d+1)!$ such simplices when $v$ is an internal vertex and less than that for vertices on the boundary of $\Omega$.
    The function  $z_\Delta$ is the linear function which is one at $v$ and vanishes on the facet of $\Delta$ opposite to $v$. If $|D_\nu|<(d+1)!$, we add artificially some of the functions $z_\Delta$ that are already in $D_\nu$ so that we end up with $(d+1)!$ not necessarily different functions, since later we will do parallelization that requires the depth of certain networks to be the same. 

 It follows from {\bf MM1} that 
$\phi_v\in \Upsilon^{\widetilde W,L}(\Relu;d,1)$, 
$\widetilde W:=3\cdot2^{\lceil\log_2 d^*\rceil-1}$, $L=1+\lceil\log_2d^*\rceil$.    
Using {\bf Parallelization, we stack} the networks $\cN_v$ that output the $\phi_v$'s to obtain a  network with 
width $W=(n+1)^d \widetilde W$ that can output any linear combination of the $\phi_v$'s.  Hence, $X({\cal K})$ is contained in $\Upsilon^{W,L}(\Relu;d,1)$ for the advertised values of $W$ and $L$.  Note that since $X({\cal K})$ is generated by parallelization of 
$(n+1)^d$ networks of width $\widetilde W$, the number of parameters used in the resulting network is 
${\cal O}((n+1)^d)$.

 To prove the second containment, we 
first observe from {\bf MM2} that each of the nodal basis functions $\phi_v\in \Upsilon^{d+1,|D_v|}(\Relu;d,1)$. We then \bf Concatenation} of the networks
$\cN_v$ used to produce $\phi_v$  and add a collation channel 
for the computation of  $S$, see \eref{nodalrep} (as it is done for the {\bf CPwL2} construction). The resulting network has 
width $W=d+2$ and depth at most $(n+1)^dd^*$.

 Several remarks are in order concerning this result.  First, note that the number of parameters used in both NNs is comparable (up to a factor depending on $d$
 )
 to the dimension $(n+1)^d$ of $X({\cal K})$.  The most important point to stress is that when using the set $\Upsilon^{W,L}$ in place of a 
 piecewise linear FEM space, we are using a much larger nonlinear family as an approximation tool.
 Indeed, the set $\Upsilon^{W,L}$ not only contains the FEM space 
 $X({\cal K})$ based on the initial choice of partitioning, but  it also contains an infinite number of  such FEM  spaces corresponding to an infinite number of possible ways to partition $\Omega$.  In fact, the NN approach is even more than a simple  generalization of the  Adaptive Finite Element Method (AFEM), where one is allowed to adaptively choose partitions (from a restricted family of partitions). It will be shown in \S \ref{SS:super} that these NNs provide a provably better approximation rate to various Sobolev and Besov classes than that provided by the FEM spaces.   While this seems like a tremendous advantage for NNs over FEMs,
 one must address (in the specific problem setting) how one (near) optimally chooses the parameters of these NNs.

 In the case when  FEMs are used to numerically solve  linear elliptic PDEs, one can employ  the Galerkin method which finds  a (near) best approximation
 to the solution to the PDE by projecting onto $X({\cal K})$.  
This  is well understood and quantified in both theory and practice through
 theorems that bound error and establish stable numerical implementation.

 When we eventually discuss   quantitative theorems for  NN approximation, we shall see that the known   results point to a tremendous
 potential increase in approximation efficiency (error versus number of parameters needed) when using NNs for the numerical solution of elliptic problems.  Whether this advantage
 can be maintained in concrete stable numerical implementation is less clear.

\subsection{Width versus depth}
 \label{SS:wd}
 An underlying issue when choosing a NN  architecture to be used in a numerical setting is whether to increase the width or the depth of the NN when one is willing to allocate more parameters to improve accuracy.
 Suppose, we fix a bound $n$ on the number of parameters to be used and ask which of the sets $\Upsilon^{W,L}$  depending on at most $n$ parameters should we employ in designing a numerical algorithm. All other issues being the same, the general consensus is that in practice deeper networks are preferable. We make some comments to explain this preference from the point of view of the enhanced approximation capacities of deeper networks.
 
First, we have shown that addition of the output functions of a NN can be implemented by either increasing width (parallelization) or depth (concatenation)  with a controlled increase in the number of parameters. However, certain operations like composition and forming minimums can only be implemented by increasing depth.  So, for  example, if we  fix a width $W=W_0$ sufficiently large to accommodate $d$ source channels and a couple of collation  channels, then we can seemingly implement as outputs from $\Upsilon^{W_0,L}$ all functions that occur as outputs of shallower networks with a comparable number of parameters.  The only rigorous statement given to this effect was for ReLU networks with $d=1$.  In this case, it was proved in    \cite{daubechies2019nonlinear}   that for any fixed $W_0\ge 4$,
 we have
 \be
 \label{depthwidth}
 \Upsilon^{n,1}(\Relu;1,1)\subset \Upsilon^{W_0,L}(\Relu;1,1),
 \ee
 provided the elements in these sets are viewed as functions on $[0,1]$(or any finite interval $[r,e]$), and $L\asymp n/W_0^2$, where the constants in $\asymp$ are absolute constants.  
 Note that 
 the number of parameters $n(W_0,L)$ determining the set $\Upsilon^{W_0,L}$ is $n(W_0,L)\asymp W_0^2L\asymp n$, and therefore comparable to the number of parameters in $\Upsilon^{n,1}$. However, $\Upsilon^{W_0,L}$ is richer, since it contains, for example compositions. So,  in this special case depth beats width.  This leads us to formulate the following general question.
 \vskip .1in
 \noindent
 {\bf Problem 2:} Are shallow networks always contained in deep networks of fixed width $W_0$ with the same number of parameters?  More precisely, is it true that if we fix the  depth $L$ and the width $W_0$, we have the inclusion
 $\Upsilon^{W,L}(\Relu;d,1)\subset \Upsilon^{W_0,L_0}(\Relu;d,1)$ whenever $L_0$ and $W$ satisfy 
 $n(W_0,L_0)\asymp n(W,L)$, where the constants in $\asymp$
depend  at most on $d$?
 \vskip .1in
 
The results given above show that {\bf Problem 2} has a positive answer if we are not concerned  about the number of parameters.  Indeed, each function in
 $\Upsilon^{W,L}(\Relu;d,1)$ is a CPwL function, 
 and therefore, according to {\bf CPwL2}, is in $\Upsilon^{d+2,L}(\Relu;d,1)$ for $L$ suitably large.  So, the key issue in {\bf Problem 2} is the control on the number of parameters.

 \subsection{Interpolation by neural network  outputs}
 \label{SS:interpolation}
 
 A common strategy for approximating a given target function $f$ is to interpolate some of its point values.  Although this is often not a good method for approximation, it is
 important to understand when we can interpolate a given set of data, and how stable is this process.  A satisfactory understanding of interpolation 
 using NNs is  far from complete.
 The purpose of this section is to frame the interpolation  problem and point out what is known.  We begin by considering interpolation by NNs with an arbitrary activation function $\sigma$
 and later specialize to ReLU activations.
 
 Let $\Upsilon ^{W,L}(\sigma)=\Upsilon^{W,L}(\sigma;d,1)$ be the set of outputs of NNs with activation $\sigma$, input dimension $d$,  output dimension $1$,  width $W$and depth $L$.  Given a finite set of data points $(x^{(i)},y_i)$, with $x^{(i)}\in\R^d$, $y_i\in\R$, $i=1,\dots,D$, a natural question is whether there is an $S\in \Upsilon^{W,L}(\sigma)$ which interpolates the given data
 in the sense that 
 \be
 \label{interpolate}
 S(x^{(i)})=y_i,\quad i=1,\dots,D.
 \ee
 This is the {\it existence question} for data interpolation.  In the case interpolants exist, let us denote by $\cS_I:=\cS_I(W,L;\sigma,d)$ the set of functions $S\in \Upsilon^{W,L}(\sigma;d,1)$ which satisfy the interpolation conditions \eref{interpolate}.
 
  Given  that we are going to use the set $\Upsilon^{W,L}(\sigma;d,1)$ for interpolation, the first question to ask is:
 \vskip .1in
 \noindent{\bf Question}: {\it Determine the  largest value $D^*:=D^*(W,L;\sigma,d)$ such that the interpolation problem    has a solution from $\Upsilon^{W,L}(\sigma;d,1)$  for all data sets of size 
 $D^*$}.
 \vskip .1in
 \noindent 
 One expects that $D^*$ should be closely related to the number of parameters 
 used to describe $\Upsilon^{W,L}(\sigma;d,1)$.  
 
 There seems to be only one general theorem addressing the interpolation problem for general activation functions $\sigma$.  It applies to the case of single hidden layer networks, that is, $L=1$,
  and is discussed  in detail  in the  survey article \cite{pinkus1999approximation}, see Theorem 5.1 in  that paper.

 \vskip .1in
 \noindent
 {\bf Interpolation from $\Upsilon^{W,1}(\sigma;d,1)$:}  If  
 $\sigma\in C(\R)$ is not a polynomial, then 
 \be 
 \label{interpsigma}
 D^*(W,1;\sigma,d) \ge W,\quad W\ge 1.
 \ee 
  
  The following sections discuss the interpolation problem for ReLU activation where more results are known.

 \subsubsection{Interpolation for $\Upsilon^{W,1}(\Relu;1,1)$}
 \label{SSS:interpolation1}
 The interpolation  question is easiest to answer for ReLU networks with  $d=L=1$.  In this case, 
 we know that the set $\Upsilon^{W,1}(\Relu;1,1)$ is almost the same as the space 
 $\Sigma_{W,1}=\Sigma_{W,1}({\rm CPwL})$ of CPwL functions subordinate to a partition of $\R$ into  $W$ intervals
 ($W-1$ breakpoints), see \eref{eq}. 
  Interpolation by functions in 
  $\Sigma_{W,1}$ is well understood, see \cite{deBoor}.  
  In the case of  $\Upsilon^{W,1}(\Relu;1,1)$, we claim  that  
  \be 
\label{Dstar1} 
D^*(W,1;\Relu,1)=W+1.
\ee  
To show this, we will use the representation \eref{onelayernet2} for functions $S$ from this set.
Consider data points $\{(t^{(i)},y_i)\}$, $i=1,\ldots,D$, with $t^{(1)}<\cdots<t^{(D)}$. We first show that when  $D= W+1$ 
 interpolation is not only possible, but there are infinitely many $S\in \Upsilon^{W,1}(\Relu;1,1)$ for which
\[S(t^{(i)})=y_i,\quad i=1,\dots,W+1.\]
We take any points $\xi_j$ that satisfy the interlacing property
  $$
 t^{(1)}<\xi_1<t^{(2)}<\cdots <\xi_W<t^{(W+1)},
  $$
and consider the function  
 \be 
 \label{genS}
 S(t):= c+\sum_{j=1}^Wa_j(t-\xi_j)_+
 \in \Upsilon^{W,1}(\Relu;1,1).
 \ee 
  We establish that interpolation is possible by induction on $W$.  When $W=1$, we choose $c:=y_1$ and $a_1$ so that $c+a_1(t^{(2)}-\xi_1)=y_2$.  For the induction step, let 
 $S_0(t):= c+\sum_{j=1}^{W-1}a_j(t-\xi_j)_+$ satisfy the first $W$ interpolation conditions.  We define $a_W$ so that we have  
 $$
 a_W(t^{(W+1)}-\xi_W) +S_0(t^{(W+1)}) =y_{W+1}.
 $$
 Then, 
 $S(t):=S_0(t)+a_W(t -\xi_{W})_+\in\Upsilon^{W,1}(\Relu;1,1)$ and satisfies all of the interpolation conditions.  This shows that the interpolation conditions can always be satisfied and that the set $\cS_I$ is infinite since we have infinitely many choices for the $\xi_i$'s.
 
 Finally, we want to  see that interpolation at $W+2$ points is generally  not  possible. For this, we use the following proposition which will also be useful when we discuss the
Vapnik–Chervonenkis (VC) dimension of NNs.
 \begin{proposition} 
 \label{P:interpolation}
 Let $n\ge 3$ and let $t^{(1)}<t^{(2)}\cdots<t^{(n)}$ be $n$ distinct arbitrary points. Let $y_1,\dots,y_n$ be such that $y_jy_{j+1}<0$,
 $j=1,\dots,n-1$.  Then, there is no $S\in\Upsilon^{n-2,1}(\Relu;1,1)$ such that $S(t^{(j)})=y_j$, $j=1,\dots,n$.
 \end{proposition}
 \noindent{\bf Proof:}  We can without loss of generality assume that $y_1>0$.
 We prove the proposition by induction on $n$.
 We first consider the case $n=3$.
If $S\in\Upsilon^{1,1}(\Relu;1,1)$, then  we have the representation 
 $S(t)=c+a(t-\xi)_+$ or $S(t)=c+a(\xi-t)_+$.
 We assume the first representation  since the second one is  treated in a similar way. In this case, we note that:
 \begin{itemize}
     \item if  $\xi\le t^{(1)}$, then $S$ is linear on $[t^{(1)},\infty)$, and therefore cannot satisfy the three interpolation conditions.  
     \item if $\xi\in (t^{(1)},t^{(2)})$, then in order for $S$ to satisfy the first two interpolation conditions, we would need $c>0$ and $a<0$.
 So, the function $S$ is then a non-increasing function of $t$ and thus $S(t^{(3)})\le S(t^{(2)})$, which shows that  $S$ cannot satisfy the third interpolation condition.
 \item if $\xi\ge t^{(2)}$ then $S$ cannot satisfy the first two interpolation conditions, since $S$ is constant on $(-\infty,\xi]$.
 \end{itemize}
 Now, we consider the induction step.  Suppose that we have proved the proposition for a value of $n\ge 3$ and consider $n+1$  interpolation points.  If 
 $S$ is any function in $\Upsilon^{n-1,1}(\Relu;1,1)$, then
 $$
 S(t)=S_0(t)+a(t-\xi_{n-1})_+, \quad \hbox{or}\quad 
 S(t)=S_0(t)+a(\xi_{n-1}-t)_+,
 $$
 where $S_0\in \Upsilon^{n-2,1}(\Relu;1,1)$, and has break points $
 \xi_1<\ldots<\xi_{n-2}$, with $\xi_{n-2}<\xi_{n-1}$. We again consider only the first possibility, since the other is handled similarly. We show that $S$ cannot satisfy the interpolation conditions by considering the following two cases:
 \begin{itemize}
     \item if   $ t^{(n)}<\xi_{n-1}$, then $S_0(t^{(j)})=S(t^{(j)})=y_j$, 
     $j=1,\ldots,n,$ and $S_0$  would contradict the induction hypothesis. Hence this case is not possible.
     \item if $\xi_{n-1}<t^{(n)}<t^{(n+1)}$, then  $S$ is a linear function on $[\xi_{n-1},\infty)$. Note that   
      $y^*:=S_0(\xi_{n-1})=S(\xi_{n-1})$ and because of $y_ny_{n+1}<0$, ${\rm sign}(S(\xi_{n-1}))={\rm sign}(y_n)$. 
    Note that  $S$ cannot satisfy the last three interpolation conditions corresponding to $t^{(n-1)},t^{(n)},t^{(n+1)}$ unless we have $t^{(n-1)}<\xi_{n-1}$. Thus, 
     $S_0(t^{(j)})=S(t^{(j)})=y_j$,  for $j=1,\ldots,n-1,$ and 
     $S_0(\xi_{n-1})=y^*$, where the sign of $y^*$ is the same as the sign of $y_n$. Therefore, according to the induction hypothesis, such $S_0$ cannot be an output of $\Upsilon^{n-2,1}(\Relu;1,1)$. 
     \end{itemize}
 This completes the proof of the proposition. \hfill $\Box$

 \subsubsection{Interpolation for $\Upsilon^{W_0,L}(\Relu;1,1)$}
 \label{SSS:interpolation2}

It is also possible to produce an interpolant to   given data by using deep networks with a fixed width.  This of course follows from \eref{depthwidth} together with what we have just proved.  However, we wish to give a direct construction because it will be used later in this paper.
\begin{proposition}
\label{int11} 
Given $D$ points $0\leq t^{(1)}<t^{(2)}<\cdots<t^{(D)}\le 1$ and values $y_j\in\R$, $j=1,\dots,D$, there is an 
$S\in \Upsilon^{3,D-1}(\Relu;1,1)$ which interpolates this data, namely,
\be 
\label{int110}
S(t^{(j)})=y_j,\quad j=1,\dots,D.
\ee 
 Therefore, $D^*(3,L;\Relu,1)\ge L+1 $, where $L\ge 1$.
\end{proposition}
\noindent{\bf Proof:}  
   We have shown above  that there is an $S$ of the form  \eref{genS} with $W:=D-1$, that satisfies the interpolation conditions \eref{int110}. We view $S$ as a function on $[0,1]$ and 
construct a special network that outputs any such $S$.   The  first channel of this special network is a source channel that pushes forward the input $t$ and the last channel is a collation  channel. This last channel is  initialized  with $0$ at layer $1$ and then
successively collects the sums $\sum_{i=1}^{j-1}a_i(t-\xi_i)_+$ 
 at layers $2,\ldots,D-1$, respectively, while  the middle 
 channel successively produces the terms $a_j(t-\xi_j)_+$, at layers $j=1,\dots,D-1$, using the inputs $t$ from the source channel. We can then output $S$ from layer $D-1$.  \hfill $\Box$

   \subsubsection{Interpolation from $\Upsilon^{W,L}(\Relu;d,1)$}
   \label{SS:interpolationRelu}
 We turn now to results that hold for  general $d\ge 1$.  There is a simple way to derive interpolation results for arbitrary $d>1$ from those for $d=1$.  Let
\[
\mathcal X:=\set{x^{(j)},\, j=1,\ldots, D}
\subseteq  \R^d
\] 
be any finite collection of data sites.
A simple measure theoretic argument shows that there exists a unit vector $v\in\R^d$ for which the points $t^{(j)}\in \R$, given by
$t^{(j)}:= v\cdot x^{(j)}$, $j=1,\dots,D$, are all distinct.  If $g$ is any univariate function which satisfies

\be 
\nonumber
g(t^{(j)})= y_j,\quad j=1,\dots,d,
\ee 
then the ridge function $f(x):=g(v\cdot x)$ satisfies $f(x^{(j)})=y^j$, $j=1,\dots,D$. 
We utilize this observation
to prove the following.
\begin{proposition}
\label{P:interpolation1} 
For any $W,L\ge 1$, we have
\be 
\nonumber
D^*(W,L;\Relu,d)\ge D^*(W,L;\Relu,1).
\ee 
\end{proposition}
\noindent 
{\bf Proof:} Given a  data set $\cX\subset \R^d$ of size $D$,  we choose $v$  as above to arrive at the points $t^{(j)}\in \R$, $j=1,\dots,D$.  If $D\le D^*(W,L;\Relu,1)$, then there is  an $S\in \Upsilon^{W,L}(\Relu;1,1)$ which satisfies $S(t^{(j)})=y_j$, $j=1,\dots,D$.  Then, the function 
$T(x):=S(v\cdot x)\in \Upsilon^{W,L}(\Relu;d,1)$ and interpolates the multidimensional data set $\cX$. \hfill $\Box$

While the above proposition is of theoretical interest, it is not used in practice because the ridge function interpolant does not reflect the local flavor of the data.  A more common scenario is to construct via ReLU networks a dual basis $\{\phi_j\}$, $j=1,\dots,D$, for the data sites, that is, a basis that satisfies the conditions
\be  
\nonumber
\phi_i(x^{(j)})=\delta_{i,j},\quad 1\le i,j\le D.
\ee 
The goal is to construct a locally supported dual basis.  In that case, the interpolation operator
\be 
\nonumber
P_{\cal X}(f):=\sum_{j=1}^Df(x^{(j)})\phi_j,
\ee 
is a bounded 
projection onto ${\rm span}\{\phi_j\}$
whenever the data sites are in $\Omega$.  The norm of this projector,
\be 
\nonumber
\|P_{\cal X}\|_{C(\Omega)\to C(\Omega)}= \max_{x\in\Omega}\sum_{j=1}^D|\phi_j(x)|,
\ee 
to a large extent determines the approximation properties of interpolation at these sites.

We have already discussed such dual bases in the context of FEM, where for the Kuhn simplicial decomposition 
${\cal K}$ of $\Omega=[0,1]^d$,
we showed that the space $X({\cal K})$ spanned by the nodal basis $\{\phi_v\}$, 
see \eref{nodal}, is contained in $\Upsilon^{W,L}(\Relu;d,1)$ for certain choices of $W$ and $L$  with the number of parameters $n(W,L)$ comparable to the dimension of $X({\cal K})$.  In this case, the nodal basis form a partition of unity $\sum_{v}\phi_v\equiv 1$ on $\Omega$ and the projection operator 
$P_{\cal X}$ is of norm one.  It follows therefore that
\begin{eqnarray} 
\nonumber
{\rm dist}(f,\Upsilon^{W,L}(\Relu;d,1))_{C(\Omega)}&\le& \|f-P_{\cal X} f\|_{C(\Omega)}\nonumber\\
&\le& \|f-S\|_{C(\Omega)}+\|P_{\cal D}(f-S)\|_{C(\Omega)}\nonumber\\
&\le&  2\, {\rm dist}(f,X({\cal K}))_{C(\Omega)},
\nonumber
\end{eqnarray} 
where we insert the best approximation $S$ to $f$ from $X({\cal K})$ to obtain the last inequality.
This allows one to deduce estimates for NN approximation from those known in FEM and also to exhibit simple linear operators which
achieve these bounds.  

\subsection{VC  dimension of ReLU 
outputs}
\label{SS:VC}

An important ingredient in understanding the approximation power of ReLU networks is the Vapnik-Chervonenkis (VC)   dimension  of the   sets of their outputs  $\Upsilon^{W,L}(\Relu;d,1)$.   
This topic is by now well studied, see \cite{bartlett2019nearly} for a summary of the most recent results.     Here, we shall only discuss the results on VC dimension that are important  for approximation.

Let  $\cF$ be a collection  of real valued  functions defined on $\Omega\subset\R^d$.  We say that a  set  $\{x^{(1)} ,\dots, x^{(n)} \}\subset \Omega$   is   {\it shattered} by $\cF$ if for each subset  $\Lambda\subset
\{1,\dots,n \}$, 
  there is a function $f=f_\Lambda \in\cF$ such that 
\be
\nonumber
f(x^{(i)})>0, \quad {\rm iff}\quad i\in\Lambda.
\ee
The maximum value of $n$ for which there exists such a collection of $n$  points that are shattered by ${\cal F}$   is  called the     Vapnik-Chervonenkis (VC) dimension of $\cF$ and is denoted by ${\rm VC}(\cF)$,  see  \cite{Vapnik}. 

In the case that $\cF$ is one of the sets $\Upsilon^{W,L}(\Relu;d,1)$, then the VC dimension of $\cF$ is the largest value of $n$ for which there exist
$n$ points such that for any assignment of signs $\e_i\in\{-1,+1\}$, there is an $S\in\cF$ such that
\be 
\label{shatter1}
\e_iS(x^{(i)})>0,\quad i=1,\dots,n.
\ee 
This follows from the fact that  whenever $S\in\cF$, then $S+c$, $c\in\R$, is also in $\cF$.
We sometimes use   property \eref{shatter1}  instead of the original definition of shattering for the output of NNs  in going forward.

Let us note that the definition of VC dimension of $\cF$ only requires the existence of one set of points where shattering takes place.  When  proving
upper  bounds on the error of approximation,  it is useful to know precisely which collections of points can be shattered.  The reader will see how this issue arises
when we use VC dimension in proving approximation results.

We are interested in describing  the VC dimension of the set $\cF$ of   outputs of ReLU networks  in terms of the number $n(W,L)$  of their parameters.
Let us now consider what is known in the special cases of interest to us.

\subsubsection{The VC dimension of   $\Upsilon^{W,1}(\Relu;d,1)$}
\label{SSS:VCshallow}
We first consider the space $\Upsilon^{W,1}$  of function of $d$ variables which is described by $n(W,1)=(d+2)W+1$ parameters.  

\begin{lemma}
\label{L:pseudo1} We have the following upper and lower bounds for the VC dimension of $\Upsilon^{W,1}(\Relu;d,1)$, $W\geq 1$:
\begin{itemize}
    \item[(i)] If $d=1$, then $ VC(\Upsilon^{W,1}(\Relu;1,1))=W+1$;
    \item[(ii)] If $d\ge 2$, then  $VC(\Upsilon^{W,1}(\Relu;d,1))\le C_0W\log_2 W$, where $C_0$ depends only on $d$;
    \item [(iii)] If $d\ge 4$, then  $VC(\Upsilon^{W,1}(\Relu;d,1))\ge c_0W\log_2 W$, where $c_0$ depends only on $d$;      
    \item[(iv)] If $d=2,3$, then   $VC(\Upsilon^{W,1}(\Relu;d,1))\ge 
W+1$.
\end{itemize}
  \end{lemma}
\noindent
  {\bf Proof:}
  (i)  The VC dimension in this case is at least $W+1$ because we can interpolate any $W+1$ data by an element from $\Upsilon^{W,1}(\Relu;1,1)$, see 
  \eref{Dstar1}.  On the other hand, the VC dimension is at most $W+1$ because of Proposition \ref{P:interpolation}. 
  (ii)  This  upper bound can be found in \cite{bartlett2019nearly}.
  (iii) The lower bounds in the case $d\geq 4$   can be derived from known lower  bounds  for the VC dimension of the  collection $\cC=\{R\}$ of sets $R$ which are the union of $W$ closed  half spaces.     Indeed,    whenever  points   $P_1,\dots,P_m$     are shattered by $\cC$, then they  are shattered by $\Upsilon^{W,1}$.  To see this, suppose that $\Lambda$ is any subset of these points. We can construct $S\in\Upsilon^{W,1}$ that is positive on this set and zero  on the remaining points as follows.  Let $R_j$, $j=1,\dots,W$, be the closed  half spaces whose union contains only the points from $\Lambda$ and none  of the rest of the $P_j$'s, $j=1,\ldots,m$.  Each of these half spaces can be represented as $w_j\cdot x+b_j\geq 0$ for some 
  $w_j\in \R^d$, $b_j\in \R$.  Then, if $\e>0$, the function 
  $$S:= \sum_{j=1}^W(w_j\cdot x+b_j+\e)_+\in \Upsilon^{W,1}(\Relu;d,1)$$
  will be
  positive on the points in $\Lambda$
  and zero on the rest of the  points $P_j$, provided we take $\e$ small enough.  It follows that 
  $$
  VC(\Upsilon^{W,1}(\Relu;d,1))\ge VC(\cC),
  $$
  and hence the lower bounds stated in (iii), follow from  the lower bounds on VC dimension of $\cC$ given in
  \cite{CKM}.
  (iv) The lower bounds in this case follow from  the fact that we can interpolate any data at any $(W+1)$ data sites, see Propostion \ref{P:interpolation1} and \eref{Dstar1}.
 \hfill $\Box$
  
  It appears that the VC dimension of $\Upsilon^{W,1}(\Relu;d,1)$ when $d=2,3$ is not completely determined because of the discrepancy between
  the upper and lower bounds in the above lemma.

\subsubsection{VC dimension of $\Upsilon^{W_0,L}(\Relu;d,1)$}
\label{SSS:VCdeep}
Next, we  consider the case where $W_0$ is fixed but sufficiently large, depending only on $d$, and $L$ is allowed to vary. Note that in this case the number of parameters 
of the network $n(W_0,L)\asymp W_0^2L$.  The following theorem gives  bounds on the VC dimension of such  networks.

\begin{theorem}
\label{L:VCdeep}
   Let  $W_0$ be fixed, and sufficiently large depending only on $d$. There are fixed constants $c_1,C_1$, depending only on $d$, such that 
\be
\label{VCUp}
 c_1 L^2    \le {\rm VC}(\Upsilon^{W_0,L}(\Relu;d,1))\le C_1 L^2.
 \ee
 \end{theorem}
   
  The upper bound in this theorem follows from Theorem 8 in    \cite{bartlett2019nearly}.   The remainder of this section will provide a proof of the lower bound in  a form which will be used later in this paper to prove certain approximation results.
  Related lower bounds are stated in Theorem 3 of \cite{bartlett2019nearly}.
  
 \subsubsection{Bit extraction using ReLU networks}
 \label{SSS:bitextraction}

 We discuss in detail  a very specific way to prove the lower bound in Theorem \ref{L:VCdeep}.  This particular construction, called {\it bit extraction},  is useful in 
 proving upper bounds for approximation using  deep neural networks.
 For  a fixed $W_0$ and $C$, depending only on $d$, the set $\Upsilon^{W_0, Cn}(\Relu;d,1)$  not only  shatters $N=n^2$ equally spaced points $x^{(1)},\dots,x^{(N)}\in \Omega:= [0,1]^d$,  but for certain bit data $y_j$, it contains
 an $S$ such that $S(x^{(j)})=y_j$, $j=1,\dots,N$.
 
 In order to avoid certain technicalities, we present this result only in the case $d=1$.  The full implementation for $d\ge 2$ can be found in
 \cite{yarotsky2018optimal} and \cite{Shencomp}.
 
 \begin{theorem}
 \label{T: bitextraction}   
Let  $N:=n^2$ with $n\ge 4$ be an even integer.  Define   $t_i:=i/N$, $ i=0,1,\dots,N$, and consider any  data   $y_i$,
$i=0,\dots,N$, with the properties:
\begin{itemize}
    \item [(i)]
 $y_{jn}=0,\quad  j=0,\dots,n$;
\item [(ii)] $y_{i+1}=y_i+\e_i$, with $\e_i\in\{-1,1\}$ for all $i=0,\dots,N-1$.
\end{itemize}
Then,   there is an
 $S\in \Upsilon^{11,15n+2}(\Relu;1,1)$,  such that 
 \be
\nonumber
S(t_i)=y_i, \quad  i=0,\dots,N.
\ee
Moreover, we have
 \be 
 \label{additionally1}
  |S(t)-S(t_i)| \le 1,\quad t\in [t_i,t_{i+1}],\quad i=0,1,\dots,N-1.
  \ee 
  \end{theorem}
\bigskip

The novelty in this theorem is that while the number of parameters in the  NN is $ Cn$, the
number of  data  points is $N+1=n^2+1$.
The theorem provides the lower bound  in \eqref{VCUp} for the VC  dimension. Indeed, at any point $t_{2i}$ not of the form $t_{jn}$,  we can assign any data $y_{2i}\in \{0,2\}$ because of property (ii).   Thus, we can shatter these points  using  the set $\Upsilon^{11,15n+2}(\Relu;1,1)$. Since there are at least $cn^2$ such points, with $c$ an absolute constant,  we have the lower bound in \eref{VCUp} in the case $d=1$ with $W_0=11$.  The general case of $d> 1$ also easily follows from this by using the method of proof used in Proposition \ref{P:interpolation1}.

Before we present  the proof of  Theorem \ref{T: bitextraction}, which is a bit laborious, we introduce some notation, make several observations, and present the general idea of the proof.  

First, note that for each $i=0,1,\dots,N-1,$ there is a unique representation 
\be 
\nonumber
t_i=\frac{i}{N}= 
 \frac{j(i)}{n}+\frac{k(i)}{N}, \quad j(i),k(i)\in \{0,1,\dots,n-1\}.
\ee
Next,  recall that  any $t\in [-1,1]$,  can be represented  as
\be 
\nonumber
t=\sum_{k=1}^\infty B_k(t)2^{-k},
\ee
where the bits  $B_k(t)\in \{-1,1\}$ of $t$ are  found using the familiar quantizer function
 \be 
\nonumber
 Q:=-\chi_{[-1,0]}+\chi_{(0,1]},
 \ee 
with  $\chi_I$ denoting the characteristic function of a set $I$.
 The first bit  of $t$, $B_1(t)=Q(t)$ and has the residual $R_1(t):=2t-B_1(t)\in [-1,1]$.   We find the later bits and residuals  recursively
 as 
\begin{equation}
    \label{bits}
 B_j(t)=Q(R_{j-1}(t)), \quad 
 R_j(t):= 2R_{j-1}(t)-B_j(t), \quad j=2,3, \ldots.
\end{equation}

Given our assigned bit sequence 
$\{\e_i\}$, $i=0,\ldots,N-1$, available to us from the values $y_i$, $i=0,1,\dots,N$, we define the numbers
\be
\label{definey}
Y_j:=\sum_{k=0}^{n-1}\e_{jn+k}2^{-k-1},\quad j=0,\dots,n-1.
\ee
Note that $Y_j\in [-1+2^{-n},-2^{-n}]\cup [2^{-n},1-2^{-n}]\subset [-1,1]$, and
the bits $B_{\nu}(Y_j)=\varepsilon_{jn+\nu-1}$, $\nu=1,\dots,n$.

The idea of proving
Theorem \ref{T: bitextraction} 
is to produce a function $S$ from the set  $\Upsilon^{11,15n+2}(\Relu;1,1)$, such that for each $i=1,\ldots N$, 
$i=j(i)n+k(i)$, 
\begin{equation}
\label{E:S-bitsum}
    S(t_{i})=
    y_{i}=\sum_{\nu=1}^
    {k(i)}
    \varepsilon_{j(i)n+\nu-1}=
    \sum_{\nu=1}^{k(i)} B_{\nu}(Y_{j(i)}), \quad k(i)=1,\ldots,n-1, 
\end{equation}
$$
S(t_{jn})=
    y_{jn}=0, \quad j=0,\ldots,n,
$$
that in addition satisfies \eqref{additionally1}. We construct $S$ by showing  that each of the functions 
\[t_i=\frac{i}{N}= 
 \frac{j(i)}{n}+\frac{k(i)}{N}\mapsto j(i), k(i), Y_{j(i)}, \chi_{\set{\nu:\,\, \nu\leq k(i)}},\qquad Y_j\mapsto B_\nu(Y_j),
 \]
are each  outputs of ReLU networks of an appropriate size. 

To do this, let $\delta=2^{-N}$ and define:
\begin{itemize}
    \item the CPwL function $J=J_N$ which  has breakpoints  at each  of the  points 
\be 
\label{breakpointsJ}
\xi_j:=j/n, \ j=1,\dots,n-1,\quad   \xi_j':=(j+1)/n-\delta, \ j=0,\dots,n-1,
\ee 
and no other breakpoints, and takes the value $j$ on the  interval $[\xi_j,\xi_j']$. 
We also require $J(1)=n$. Note that $J$ has  the property
$$
J(t_i)=j(i), \quad  i=0,1,\dots,N-1;
$$
\item the CPwL function $K(t)=K_N(t):=J(nt-J(t))$. Observe that the key property of $K$ is
\be 
\label{properties} 
K(t_i)=k(i), \quad  i=0,1,\dots,N-1;
\ee 
\end{itemize}

Next, we would like to implement quantization by a neural network.  However, the function  $Q$ is not continuous, 
and so we cannot exactly reproduce $Q$.  Instead, we use a surrogate 
\be 
\label{surrogatequantizer}
\hat Q(t)=-1+(\frac{1}{\delta}t+1)_+-(\frac{1}{\delta}t-1)_+,
\ee 
where $\delta:=2^{-N}$.  
The surrogate  $\hat Q$ is in $\Upsilon^{2,1}(\Relu;1,1)$ and   coincides with $Q$ on   $[-1,1]\setminus [-\delta,\delta]$. 

We define the surrogate bits $\hat B_\nu(t)$ for $t\in [-1,1]$   by using $\hat Q$ in place of $Q$ in the recursive definition of $B_\nu$, described in \eref{bits}.
Because of  the choice of $\delta$, 
$\hat B_\nu$ can be used in place of $B_\nu$ to compute the bits of $t$ whenever  $t$ has the representation
\be 
\label{Bcomputes}
t=\sum_{\nu =1}^k B_\nu(t)
2^{-\nu},
\quad \hbox{with}\quad k\leq N-1.
\ee 
For  such a $t$, we have $\hat B_\nu(t)=B_\nu(t)$, $\nu=1,\dots, N-1$.

Finally, we introduce 
\begin{itemize}
\item  the CPwL function $Y$ which  has exactly the same breakpoints as $J$,  see \eref{breakpointsJ},  and satisfies 
\be 
\label{defY}
Y(\xi_j)=Y(\xi_j')= Y_j, \quad j=0,1,\dots,n-1,
\ee 
with $Y_j$ defined in \eref{definey} and $Y(1)=0$.  
\end{itemize}
\noindent
The function $S$ will be the output of a special neural network  of width $W=11$
and depth $L=15n+2$, which is a concatenation of four special networks that we describe below.
  The top channel of each of these networks  is a source channel which simply passes forward the input $t$.   Some of the other channels are collation channels and are occupied by zeros in their first layers so that they  can be used later for  passing forward certain function values.   
  
  We want to point out  that our  construction is probably not optimal in the sense that it does not provide a NN with the best possible minimal width and depth that outputs $S$. In addition,  some  of the channels in our  NN are ReLU free.  We have discussed earlier how we can construct a true ReLU network with the same outputs as a network that has ReLU free nodes. 
  
 In going further, we note that any 
 CPwL function $T$ with $k$ break points  is in
$\Upsilon^{3,k}(\Relu;1,1)$, where the network used to output $T$ has  one source channel, one computational channel, and one collation channel that  collects the successive terms we have computed, see the construction in Proposition \ref{int11}.

{\bf Proof of Theorem \ref{T: bitextraction}:} We can now give  the proof of Theorem \ref{T: bitextraction}.  The network $\cN$ which outputs
the function $S$ of the theorem is a concatenation of five special networks $\cN_1,\cN_2,\cN_3,\cN_4,\cN_5$. Each of them has a source channel as its first channel.  It  pushes forward the input $t\in[0,1]$.  The first of these networks outputs $K(t)$, the second outputs $Y(t)$, the third takes input $K(t)$ and $Y(t)$ and outputs a CPwL function  $\tilde S$  which almost satisfies the theorem.  Namely, it satisfies the interpolation conditions and it also satisfies
\eref{additionally1} except for a small set of $t$ values.  
The last two networks make a technical correction to $\tilde S$ to obtain the desired $S$ which satisfies \eref{additionally1} for all $t\in[0,1]$.  We now describe these five networks.  All of them have width at most $11$ and we make the width exactly $11$ by adding zero channels.  The depth of each network is also controlled so that the final network has depth $L=15n+2$.

\vskip .1in

\noindent
{\bf First NN:}      This network, which we denote by $\cN_1$, has depth 
$4n-2$ and for any input $t\in[0,1]$   outputs the function value $K(t)$.
From our remarks on interpolation,  see Proposition \ref{int11}, we know that $J(t)$ is the output of a special ReLU   network $\cN_0$  of width $W=3$ and depth $2n-1$, where channel three is a collation channel.  The CPwL function $K$ is the output of a ReLU network $\cN_1$ of width $W=3$ and depth $4n-2$, which is obtained by concatenating the network $\cN_0$ for $J$ with itself and using $nt-J(t)$ as the input to the second of these networks. 
The third channel is  a collation channel, used first to build $J(t)$. Once $J(t)$ is computed, it sends this value  as an input to the 
$2n$-th layer. Then, it is zeroed out  by  assigning  a weight $0$, and subsequently used as a collation   channel to  build  $K(t)$. It follows from \eref{properties} that the output of this network is $k(i)$ when the input is $t_i$.  
We add eight other  channels   with zero parameters. These channels will be used later. 
 
\vskip .1in
\noindent
{\bf Second  NN:}  
The second  network    $\cN_2$  takes the input $t$ from channel one and  outputs the  CPwL function 
$Y(t)$ which belongs to $\Upsilon^{3,2n-1}(\Relu;1,1)$.    This network has depth $2n-1$ and only needs three channels, 
 but we augment it with eight more  channels.  
After a concatenation with  the existing network $\cN_1$, it uses  Channel 2 to compute the terms involved in $Y(t)$,  while channels 3 and 4  push forward the values  $K(t)$ and the terms involved in $Y(t)$, respectively. Channels $5$ to $11$ have all parameters  zero.   
Note that after this concatenation, we have available to us as outputs $t$, coming from the source channel, $K(t)$, kept in channel 3, and $Y(t)$ kept in channel 4.
\vskip .1in
\noindent
{\bf Third  NN:} 
This network takes as inputs  $K(t),Y(t)$ and outputs a  function $\tilde S$ which satisfies the interpolation conditions and coincides with the desired  $S$ except for a small subset of $[0,1]$.  To describe this network, we shall use
the CPwL function $T$ with break points $-1,1,2$, defined as
\begin{equation}
\nonumber
T(t)=\begin{cases}
-1, \quad \quad t\leq -1,\\
t,\quad \quad -1\leq t\leq 1,\\
2-t, \quad 1\leq t\leq 2,\\
0, \quad \quad \quad t\geq 2.
\end{cases}    
\end{equation}
Since
$T(t)=-1+(t+1)_+-2(t-1)_++(t-2)_+$, it belongs to
$\Upsilon^{3,1}(\Relu;1,1)$.
Note that $T$ is the identity on the interval $[-1,1]$ and has the important property that for $t\in[-1,1]$ and for each  $1\le k<n$, we have
\be 
\label{importantT}
T\lr{B_\nu (t)+ 3(\nu -k)_+} =  B_\nu(t),\quad 1\le \nu \le k,
\ee 
 and is zero otherwise, since  $3(\nu -k)_+=0$ when $\nu\le k$ 
 and $3(\nu -k)_+\geq 3$
 when $\nu>k$. It follows from \eref{importantT} that 
\be 
\nonumber
\sum_{\nu=1}^n T\lr{B_\nu(t)+ 3(\nu-k)_+} =  \sum_{\nu=1}^{k} B_\nu(t).
\ee 
Now, for $i=0,\dots,N-1$, consider one of our points $t_i$   which is not a multiple of $n$, that is, $k(i)\neq 0$. Then $Y(t_i)= Y_{j(i)}$, $K(t_i)=k(i)$, and
\be 
\nonumber
\sum_{\nu=1}^n T(B_\nu(Y(t_i))+ 3(\nu-K(t_i))_+) = \sum_{\nu=1}^{k(i)} B_\nu(Y_{j(i)})=\sum_{\nu=1}^{k(i)} \e_{j(i)n+\nu-1}= y_i.
\ee 
Since we cannot produce $B_\nu$ with a ReLU network, we use the surrogate $\hat B_\nu$ in its place. This leads us to define the following function
\be 
\label{defineS} 
\tilde S(t):= 
\sum_{\nu=1}^n T(\hat B_\nu(Y(t))+ 3(\nu-K(t))_+), \quad t\in [0,1].
\ee 
This function satisfies the interpolation conditions \eref{E:S-bitsum} since   the bits $\hat B_\nu(Y(t))=B_\nu(Y(t))$, $\nu=1,\dots,n$, whenever $t$ is one of the points $t_i$, $i=0,\dots,N$, where interpolation is to take place.
In addition,  since 
for $j=0,\ldots,n-1$,
$K(t_{jn})=0$  and $K(1)=J(0)=0$, we have
$$
\tilde S(t_{jn})=
\sum_{\nu=1}^n T(\hat B_\nu(Y(t_{jn}))+ 3\nu)=0,\quad j=0,\ldots,n,
$$
because of the definition of $T$.

Next, we describe how  
$\tilde S$ is an output of a ReLU network $\cN_3$ with inputs $Y(t),K(t)$.  The network $\cN_3$  is organised as follows.
Channel 1 is left to be a source channel that forwards the value of $t$.
Channels 2 and 3 are occupied with the values of $K(t)$ and $Y(t)$, forwarded to the next layers. Channel 4 computes 
$(\nu-K(t))_+$ in layer $\nu$, for  $\nu=1,\ldots,n$.    Channel 5 computes the residual 
$R_{\nu-1}(Y(t))$ in layer $\nu$, 
$\nu=1,\ldots,n$ (this is a ReLU free channel), where
$R_0(Y(t))=Y(t)$, 
and  $R_\nu(Y(t))=2R_{\nu-1}(Y(t))-\hat B_\nu(Y(t))$,
$\nu=1,\ldots,n$.
Channels 6 and 7 
implement the network for $\hat Q$ and compute consecutively
$\hat B_1(Y(t))=\hat Q(Y(t))$, 
$\hat B_\nu=\hat Q(2R_{\nu-2}(Y(t))-\hat B_{\nu-1}(Y(t)))$, for $\nu=2,\ldots,n$. Channels 
8, 9, and 10  implement  $T$.  The 11-th channel successively  adds  the $T$ values   in the sum \eref{defineS}, and therefore $\cN_3$ outputs $\tilde S$.
 In total, the entire  network $\cN_3$ has $(n+1)$ layers and  width 11. 
 
 We have already observed that $\tilde S(t_i)=y_i$ and so the interpolation conditions are satisfied.
 The reader can imagine that we can take a max and min with an upper and lower CPwL to obtain the control
 \eref{additionally1}.   The network $\cN_4$ will do precisely that.  So, the remainder of the proof is to give one 
 such construction.
 
  We claim that the output  $\tilde S$ of $\cN_3$ already satisfies the inequalities
 \be 
 \label{additionally}
  |\tilde S(t)-\tilde S(t_i)| \le 1,\quad t\in[t_i,t_{i+1})\cap  \Omega_N,\quad i=0,1,\dots,N-1,
  \ee 
  where
 \be 
 \nonumber
 \Omega_N:= [0,1]\setminus \bigcup_{j=1}^{n}(t_{jn} -\delta,   t_{jn}), \quad \delta:=2^{-N}.
 \ee  
We verify this property when 
$t\in [0,1/n)\cap\Omega_N=[0,1/n-\delta]$ since the verification on the intervals $[j/n,(j+1)/n)\cap \Omega_N$, $j=1,\dots,n-1$, is the same. For 
$t\in [0,1/n-\delta]=[0,t_n-\delta]$ we have, see \eref{defY},  $Y(t)=Y_0$, and therefore for $\nu =1,\dots,n$, 
$\hat B_\nu(Y(t))=\hat B_\nu(Y_0)=B_\nu(Y_0)=\e_{\nu -1}$. Thus, see \eref{defineS}, we have
\begin{equation}
    \nonumber
\tilde S(t)= 
\sum_{\nu=1}^n T(\e_{\nu-1}+ 3(\nu-K(t))_+), \quad t\in [0,t_n-\delta]. 
\end{equation}
We consider the following  cases:
\begin{itemize} 
\item if $t\in[0,t_{1}-\delta]$, then $K(t)=0$ and  
\begin{equation}
    \label{plo}
    \tilde S(t)=\sum_{\nu=1}^{n}
T(\e_{\nu-1}+3\nu)=0,
\end{equation}
and thus \eref{additionally}
is satisfied for these $t$.
\item if $t\in[t_k,t_{k+1}-\delta]\subset [0,t_n-\delta]$, with $1\le k<n$,  then $K(t)=k$ and  
\begin{equation}
    \label{pli}
\tilde S(t)=\sum_{\nu=1}^{k}\e_{\nu-1}=y_k,
\end{equation}
and thus  \eref{additionally}
is satisfied again.
\item if $t\in (t_{k}-\delta,t_{k})$, $1\leq k\le n-1$, then $k-1\leq K(t)<k$ and
\begin{eqnarray}
\nonumber
3(\nu-K(t))_+ \ {\rm is}\ 
\begin{cases}
=0, \quad \nu\le k-1,\\
\leq 3, \quad \nu=k,\\
>3,\quad \nu\ge k+1.
\end{cases}
\end{eqnarray} 
It follows  from the definition of $T$ that
\begin{eqnarray}
\nonumber
\!\!\!\!\!\!\!\!\!\tilde S(t)\!\! &=&\!\!\!\sum_{\nu=1}^n T(\e_{\nu-1}+ 3(\nu-K(t))_+)\\
\nonumber
 \!\!&=&\!\!\!\sum_{\nu=1}^{k-1} T(\e_{\nu-1})+
 T(\e_{k-1}+3(k-K(t))_+)=y_{k-1}+ \eta,
 \end{eqnarray}
 with $|\eta|\leq 1$, and therefore \eref{additionally} is satisfied in this case as well. 
\end{itemize}

In summary, the function $\tilde S$ satisfies  the properties we want except for control on the small intervals that make up the complement of  $\Omega_N$.
We do not have a bound for $\tilde S$ on these small intervals.
Our last construction will be to take care of these intervals while leaving $\tilde S$ unchanged outside of them.

To see how to do this, 
we concentrate on the interval $[t_{(j-1)n},t_{jn}]$, for $j=1,\dots,n$, and let $I_j:=(t_{jn}-\delta,t_{jn} {]}$.  We know from our analysis that $\tilde S(t)=y_{jn-1}=\eta_j$ for  $t \in [t_{jn-1},t_{jn}-\delta ]$,  where $\eta_j=\pm 1$.  Assume for now that  $\eta_j=1$. Also, recall that $\tilde S(t)=0$
for $t \in [t_{(j-1)n},t_{(j-1)n+1}-\delta ]$.  If $M:=\|\tilde S\|_{C(\Omega_N)}\geq 1$, we consider the CPwL function $U_j$ whose graph  passes through the points  
 $(t_{(j-1)n},0)$,  $(t_{(j-1)n+1}-\delta,M)$, $(t_{jn-1},M)$,  $(t_{jn}-\delta,1)$,  and  $(t_{jn},0)$,
and is otherwise linear between these points.  Then, $U_j(t)\ge \tilde S(t)$ on $[t_{(j-1)n},t_{jn}-\delta ]$, and  thus on $[t_{(j-1)n},t_{jn}]\setminus I_j$ the function $\min\{\tilde S,U_j\}=\tilde S$ 
. In addition, $\min\{\tilde S,U_j\}$  will have values between $0$ and $1$ on $I_j$.  This is the correction we want on $I_j$.
We then define $U:=\sum_{j\in\Lambda_+} U_j\chi_{[t_{(j-1)n},t_{jn}]}$ where $\Lambda_+$ is the set of $j$'s for which $\eta_j=+1$.
In a similar way we define a lower envelope $\hat U$ for the $j$'s such that $\eta_j=-1$.    We can then take
\be  
\nonumber
S:= \max\{\min\{\tilde S,U\},\hat U\}.
\ee 
The function $S$ satisfies the conclusions of the theorem and we only have to see how it is outputted by a suitable neural network.
Each of the functions $U,\hat U$ have at most $4n+1$ breakpoints and hence are in $\Upsilon^{3,4n+1}(\Relu;1,1)$.  

\vskip .1in
\noindent
{\bf Fourth and Fifth NNs:}   
These are the networks
$\cN_4$ and $\cN_5$ outputting  $U$ and $\hat U$. We augment them  with collation channels so that they have   width $11$. Since they already have a source channel (channel 1), there is no need to add such a channel.

\vskip .1in
\noindent
{\bf The network $\cN$:}   
We use a construction similar to the one in {\bf MM4} to output $S$
by concatenating the networks for 
$\tilde S$, $U$, and $\hat U$. 
Following the construction in {\bf MM4}, we end up with a network with width W=11 (the same as the one for 
$\tilde S$) and depth $L=15n+2$  where we added the depth 
$7n-2$ of the network for $\tilde S$, the depths of the network for $U$ and $\hat U$, each of which is $4n+1$, and two more layers to perform the $\min$ and $\max$. 
 This completes the proof of the
 theorem.
 \hfill $\Box$

\begin{remark}
\label{R:VCconstruction}
We make some final remarks on the above construction.  We have used the fact that the $t_i$'s are $N=n^2+1$ equally spaced points.  It would be interesting
and useful to clarify for which other patterns of univariate points the construction can be done.    We know that we cant increase the number of points  significantly because of the upper bound in \eref{VCUp}.  Note that in the case $d>1$, we
can construct a similar interpolant for points $x^{(i)}\in \R^d$ if there is a $v\in \R^d$ such that $v\cdot x^{(i)}=t_i$, $i=1,\dots,N$.  Again, it would be
useful to know exactly when we can interpolate certain patterns of values like the  $y_i$'s with
$Cn^2$ points from $\R^d$.  The constructions in \cite{yarotsky2018optimal} and \cite{Shencomp} show that this is possible on
equally spaced tensor product grids.
\end{remark}

\section{Classical model classes: smoothness spaces}
\label{S:modelclasses}

In order to prove anything quantitative about the rate of approximation of a given target function $f$, one obviously needs to assume something about $f$.
Such assumptions are referred to as  {\it model class assumptions}.   We say that a set $K$ in a Banach space $X$ 
is a {\it model class} of $X$ if $K$ is compact in $X$.
The classical model classes for  multivariate functions are the unit balls of smoothness spaces such as Lipschitz, H\"older, Sobolev, and Besov spaces.
We give a brief (mostly heurestic) review of these spaces in this section.  A detailed development of these spaces can be found in the  standard references, see e.g. \cite{adams2003sobolev,peetre1976new,stein1970singular,devore1993besov,devore1998nonlinear}.

We consider these spaces on the domain $\Omega:=[0,1]^d$.  All definitions and properties extend to more general domains  such as Lipschitz domains in $\R^d$.  We use standard multivariate notation.

\subsection{$L_p$ spaces}
\label{SS:Lp}
As a starting point, we recall that the $L_p(\Omega)$ spaces consist of all Lebesgue measurable functions $f$ for which $|f|^p$ is integrable.
  We define
\be
\nonumber
\|f\|_{L_p(\Omega)}:=\left(\int_\Omega |f(x)|^p\, dx\right )^{1/p},\quad 0< p< \infty.
\ee
This is a norm when $1\le p<\infty$ and a quasi-norm when $0<p<1$.  
When $p=\infty$, one usually takes $X=C(\Omega)$, the space of continuous functions on $\Omega$ with the uniform norm
\be
\nonumber
\|f\|_{C(\Omega)}:=\sup_{x\in\Omega}|f(x)|.
\ee
However, on occasion we,  need the space $L_\infty(\Omega)$ consisting of all functions that are essentially bounded on $\Omega$ with
\be
\nonumber
\|f\|_{L_\infty(\Omega)}:=  {\rm ess} \sup_{x\in\Omega}|f(x)|.
\ee
We assume throughout that the reader is familiar with the standard properties of these spaces.

\subsection{Sobolev spaces}
\label{SS:Sobolev}
We begin by defining smoothess spaces of continuous functions. 
If $r$ is a positive integer then $C^r:=C^r(\Omega)$, $\Omega=[0,1]^d$,  is the set of all continuous functions $f$ defined on $\Omega$, which have classical derivatives
$D^\alpha f$ for all $\alpha$ with $|\alpha|=r$, where  $|\alpha|:=\sum_{j=1}^d|\alpha_j|=r$.  We equip this space with the semi-norm
\be
\nonumber
|f|_{C^r}:= |f|_{C^r(\Omega)}:=\max_{|\alpha|=r}\|D^\alpha f\|_{C(\Omega)}.
\ee
A norm on this space is given by $\|f\|_{C^r(\Omega)}:=|f|_{C^r(\Omega)}+\|f\|_{C(\Omega)}$.

The Sobolev spaces (of integer order) generalize the spaces $C^r$ by imposing weaker assumptions on the derivatives $D^\alpha f$.  First,
 the notion of weak  (or distributional) derivatives $D^\alpha f$ is introduced in place of classical derivatives.  Then, for any $1\le p\le \infty$,
the Sobolev space $W^r(L_p(\Omega))$ is defined as the set of all $f\in L_p(\Omega)$ such that $D^\alpha f\in L_p(\Omega)$ for all
$|\alpha|=r$. We  equip this space with the semi-norm
\be
\nonumber
|f|_{W^r(L_p(\Omega))}:=\max_{|\alpha|=r}\|D^\alpha f\|_{L_p(\Omega)},
\ee
and obtain a norm on this space by 
$
\|f\|_{W^r(L_p(\Omega))}:= |f|_{W^r(L_p(\Omega))} +\|f\|_{L_p(\Omega)}.
$

\subsection{Besov  spaces}
\label{SS:Lipschitz}

The Sobolev  spaces above  are not sufficient   because they  only classify smoothness  for integer values $r$.    There is a long history
of introducing smoothness spaces for any order $s>0$.  This began with Lipschitz and H\"older spaces and culminated with the Besov spaces
that we define in this section.  

Given a function $f\in L_p(\Omega)$, $0<p\le\infty$, and any integer $r$, we define its modulus of smoothness of order $r$ as
\be
\nonumber
\omega_r(f,t)_p:=\sup_{0<|h|\le t}\|\Delta_h^r(f,\cdot)\|_{L_p(\Omega)},\quad t>0,
\ee
where $h\in \R^d$ and $|h|$ is it Euclidean norm.  Here,  $\Delta_h^r$, is the $r$-th difference operator, defined by
\be
\nonumber
\Delta_h^r(f,x) := \sum_{k=0}^r (-1)^{r-k}\binom{r}{k} f(x+kh),\quad x\in \Omega\subset \R^d,
\ee
where this difference is  set to  zero whenever one of the points $x+kh$ is not in $\Omega$.  
It is easy to see that for any $f\in L_p(\Omega)$, we have $\omega_r(f,t)_p\to 0$, when $t\to 0$.   How fast this
modulus tends to zero with $t$ measures the $L_p$ smoothness of $f$. 

  For example, the Lipschitz space ${\rm Lip}(\alpha,p)$ for $0<\alpha\le 1$
and $0<p\le \infty$ consist of those functions $f\in L_p(\Omega)$ for which
\be
\nonumber
\omega_1(f,t)_p\le Mt^\alpha,\quad t>0,
\ee
and the smallest $M$ for which this holds is the semi-norm $|f|_{{\rm Lip}(\alpha,p)}$.   Again, we obtain a norm on this space by simply
adding $\|f\|_{L_p(\Omega)}$ to the semi-norm.

The Besov spaces  generalize the measure of smoothness in two ways.   They allow for $r$ to be replaced by any $s>0$ and 
they introduce a finer way to measure decay of the modulus as $t$ tends to zero.  This finer decay is controlled by a new parameter $0<q\le\infty$. 

If $f\in L_p(\Omega)$, $0<p,q\le\infty$ and $s>0$, the space $B_q^s(L_p(\Omega))$ is defined as the set of functions 
$f$ for which
\be
\label{Besovseminorm}
|f|_{B_q^s(L_p(\Omega))}:= \| t^{-s}\omega_r(f,t)_p\|_{L_q((0,\infty),dt/t)}<\infty,\quad \hbox{where}\,\,r:=\lfloor s\rfloor +1.
\ee
Notice here that the $L_q$ norm is taken with respect to the Haar measure $dt/t$.  The case $q=\infty$ is simply the supremum norm
over $t>0$.  The norm on this space is $\|f\|_{B_q^s(L_p(\Omega))}:=|f|_{B_q^s(L_p(\Omega))} +\|f\|_{L_p(\Omega)}$.

The Besov spaces are now a standard way of measuring smoothness.  Functions in this space are said to have smoothness of order $s$
in $L_p$ with $q$ giving a finer gradation of this smoothness.   We mention without a proof a few of the properties of these spaces that are frequently used in analysis.

First, notice that when $s\in (0,1)$ and $q=\infty$, these spaces are the ${\rm Lip}(s,p)$ spaces.  However, the space $B_\infty^1(L_p(\Omega))$ is not ${\rm Lip}(1,p)$ since $\omega_2$ is used in place of $\omega_1$ in the definition \eqref{Besovseminorm}, thereby resulting in a slightly larger space.  A second useful remark is that in  \eqref{Besovseminorm} we could have used any $r>s$ and obtained the same space and an equivalent
norm.  When we insert $q$ into the picture, the requirement for $f$ to be in the space $B_q^s(L_p(\Omega))$ gets stronger as $q$ gets smaller,
namely, we have the following embeddings: 
\vskip .1in
 
 \noindent
 {\bf BE1:} Let  $0< p\le\infty$.  If  $s>s'$ and $0<q,q'\le \infty$ or $s=s'$ and $q\le q'$, we have $|f|_{B_{q'}^{s'}(L_p(\Omega))}\le C|f|_{B_{q}^{s} (L_p(\Omega))}$ with the constant $C$ independent of $f$.
 
 \vskip .1in
 \noindent
 {\bf BE2:} If  $0< p<p'\le\infty$ and $0<q,q'\le \infty$  then  $|f|_{B_q^s(L_p(\Omega))}\le |f|_{B_{q'}^{s} (L_{p'}(\Omega))}$.

\vskip .1in
We also have the well known Sobolev embeddings for Besov spaces.
\vskip .1in
\noindent
{\bf BE3} Let $0<p\le\infty$.  For any $s>0$ and $0<q\le \infty$, we have that the unit ball $ U(B_q^s(L_\tau(\Omega)))$, $0<q\le \infty$, is a compact subset of
$L_p(\Omega)$ 
whenever $s>\frac{d}{\tau}-\frac{d}{p}$.
\vskip .1in

\begin{figure}[h]
  \centering
\includegraphics[scale=.5]{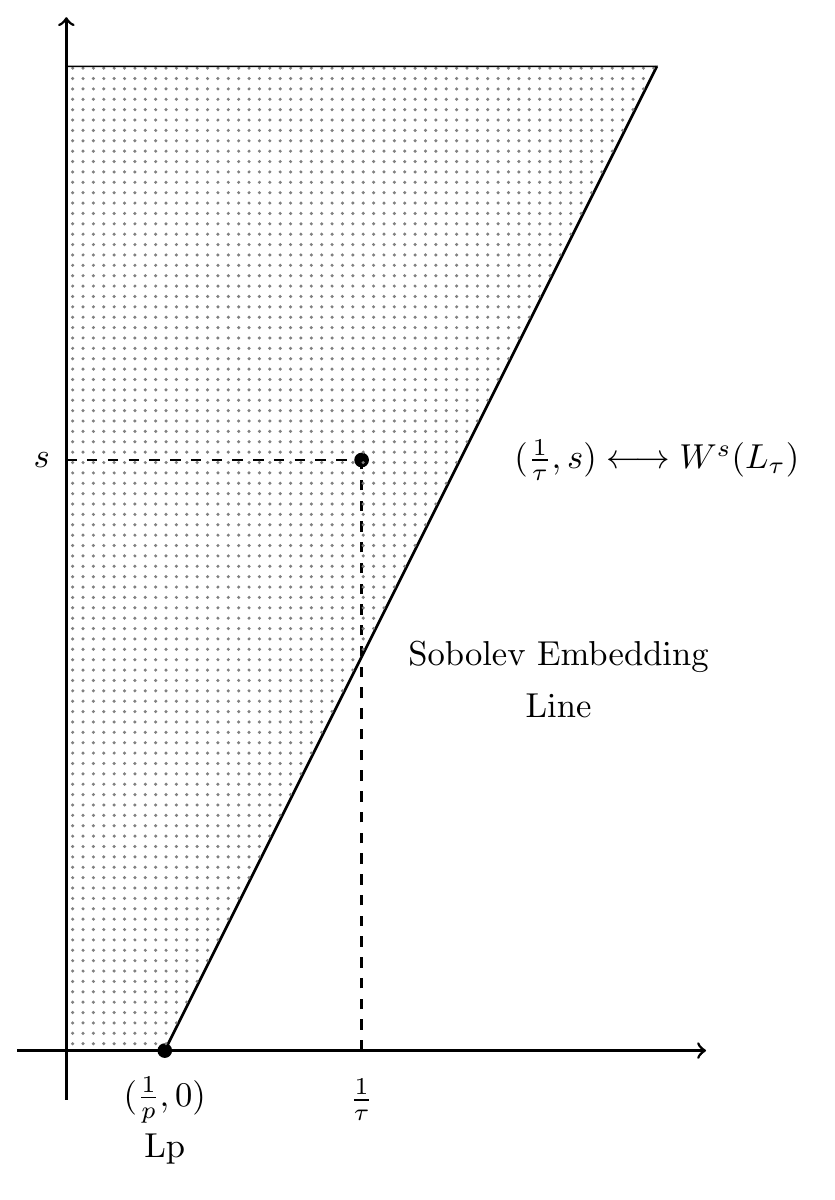}
\caption{The Sobolev embedding theorem.}
\label{fig2.pdf}
\end{figure}

There is a simple graphical way to describe these embeddings that we shall refer to  in this paper.  We use the upper right
quadrant of $\R^2$ to graphically represent smoothness spaces.  We can write any point in this quadrant as $(1/p,s)$ with $0<p\le \infty$ and $s\ge 0$.
We think of any such point as corresponding to a smoothness space with smoothness of order $s$ measured in $L_p$.  For example, the space
Lip$(\alpha,L_p)$ can be thought of as corresponding to the point $(1/p,\alpha)$, and all Besov spaces $B_q^s(L_p(\Omega))$, 
$0<q\leq \infty$, are identified with the same point $(1/p,s)$.
In terms of this graphical description, given an $L_p(\Omega)$ space, the smoothness spaces embedded into $L_p(\Omega)$ are the ones that correspond to points $(1/\tau,s)$ that lie on or above
the line with equation $s=d(1/\tau-1/p)$.  Those corresponding to points strictly 
above this line are compactly embedded.   These embedding results are summarized in Figure \ref{fig2.pdf}.

\subsubsection{Atomic decompositions}
\label{SS:atomic}

An   often used fact about Besov spaces is that functions in these spaces can be described by certain so-called {\it atomic decompositions}.  Historically, this began with the Littlewood-Paley decompositions, see \cite{frazier1991littlewood}.  In the case of  approximation by ReLU networks, the two most relevant decompositions are
those using spline functions or wavelets.  We discuss  the case of spline decompositions.  Details and proofs can be found,  for example, in \cite{DP}.

 Let $r\geq 1$ be a positive integer and consider the univariate cardinal B-spline $N_r$ of order $r$ (degree $r-1$), which is defined by
   \be
   \label{Bspline}
   N(t):= N_r(t):=
   \frac{r}{r!}\sum_{k=0}^r(-1)^{r-k}\binom{r}{k}
   (k-t)_+^{r-1},\quad t\in \R.
   \ee 
      The function   $N_r$ is a piecewise polynomial of degree $r-1$, is  in $C^{r-2}(\R)$, and is  supported on $[0,r]$. With this normalization of the
      B-spline we have $\|N_r\|_{C(\R)}\le 1$.
   
   The multivariate cardinal B-splines are defined as tensor products
   \be 
   \label{Bmulti}
  \nonumber
   N(x):=N_r(x_1,\dots,x_d):=N_r(x_1)\cdots N_r(x_d),\quad x=(x_1,\dots,x_d)\in\R^d.
   \ee 
   We do not indicate the dependence on $r$ when it is known from context.
   Let $\cD$ denote the collection of dyadic cubes in $\R^d$ and let $\cD_k$ denote the dyadic cubes of side length $2^{-k}$.  We also use the notation
   $\cD_+:=\bigcup_{k\ge 0} \cD_k$.  If $I\in\cD_k$ has smallest vertex $2^{-k}j$ with $j\in\Z^d$, we let
   \be 
   \label{dyadicB}
   N_I(x):=N_{I,r}(x):=N(2^kx-j), \quad x\in\R^d.
   \ee

   The splines $N_I$ provide an atomic decomposition for many function spaces and, in particular, the $L_p$, Sobolev, and Besov spaces.  Consider, for example,
     $\Omega=[0,1]^d$ and  denote by $\cD_k(\Omega)$ the set of those $I\in \cD_k$ for which the support of
     $N_I$ nontrivially intersects $\Omega$.  Then each $f\in L_1(\Omega)$ has a representation
   \be 
   \label{Brep}
   f=\sum_{I\in\cD_+(\Omega)} c_I(f)N_I,
   \ee 
   where  the $c_I$'s are linear functionals on $L_1$, and  $\cD_+(\Omega)=\bigcup_{k\ge 0} \cD_k(\Omega)$.  The representation \eref{Brep} is not unique since the
   $N_I$'s are not linearly independent.  However, we can fix the $c_I$'s so that all properties stated below in this section are valid.
   
  We can characterize membership of $f$ in a Besov space $B_q^s(L_p(\Omega))$ in terms of the decomposition \eref{Brep}, see Corollary 5.3 in \cite{DP}.   Namely, $f\in B^{s}_q(L_p(\Omega))$,   $0<s<\min\{r,r-1+1/p\}$, and $0<q,p\le \infty$ if and only if $f$ has the representation \eref{Brep} with coefficients $c_I(f)$ satisfying
   \be 
   \label{Bnorm}
   \|f\|'_{B_q^s(L_p(\Omega))}:=  
   \left\{ \sum_{k=0}^\infty 2^{skq}\left(\sum_{I\in\cD_k(\Omega)}|c_I(f)|^p|I|\right)^{q/p}            \right\}^{1/q}<\infty,
   \ee 
   for $0<q,p<\infty$, with the obvious modifications when either $p$ or $q$ is infinity.  Moreover, $\|\cdot\|'$ is equivalent to the usual Besov norm.
   This fact is the starting point for proving many approximation theorems for functions in Besov spaces.

\subsection{Interpolation of operators} 
\label{SS:interpolatiooperators} 
Next, we  mention how from known upper bounds for approximation error on a model class, we can derive new upper bounds on
a spectrum of new model classes by using results from the theory of interpolation of operators.  We assume the reader is familiar
with the rudiments of the theory of interpolation spaces via the real method of interpolation, see either \cite{BL1} or \cite{BS}.  

Given two Banach spaces $X,Y$ with (for convenience)  $Y$ continuously embedded in $X$, the real method of interpolation generates a family of new Banach spaces $(X,Y)_{\theta,q}$, $0<\theta < 1$, $0<q\le\infty$, which interpolate
between them.  These spaces are defined via what is called the $K$ functional for the pair
\be 
\nonumber
K(f,t):=K(f,t;X,Y):=\inf_{g\in Y} \|f-g\|_X+t|g|_Y,\quad t>0,
\ee 
where $\|\cdot\|_X$ is the norm on $X$ and $|\cdot|_Y$ is a semi-norm on $Y$ \footnote{When $Y$ is not continuously embedded in $X$, we use $\|\cdot\|_Y$ in the definition of $K$.}.  The space $(X,Y)_{\theta,q}$,
$0<\theta< 1$, $0<q\le \infty$, then consists of all $f\in X$, such that
\be 
\label{intnorm} 
\|f\|_{(X,Y)_{\theta,q}}:= \| t^{-\theta}K(f,t)\|_{L_q((0,\infty),dt/t)}<\infty,
\ee 
where the $L_q$ norm is taken with respect to the Haar measure $dt/t$. The important fact for us is that for classical pairs $(X,Y)$ of spaces such as $L_p$ and Besov/Sobolev spaces, the interpolation spaces are known and can be used to easily extend known error estimates for approximation.  We mention two typical approximation results.  By $U(Y)$ we mean the unit ball of the space $Y$.
\vskip .1in
\noindent 
{\bf Extend 1:}  If $\Sigma_n\subset X$ is a set that provides the approximation error 
$$
\inf_{S\in \Sigma_n}\sup_{f\in U(Y)}\|f-S\|_X=:E(U(Y),\Sigma_n)_X=\e_n,
$$
then for the space $Z=(X,Y)_{\theta,q}$, $0<\theta< 1$ and $\ 0<q\le \infty$, we have
\be 
\nonumber
E(U(Z),\Sigma_n)_X\le \e_n^\theta.
\ee

\vskip .1in
\noindent 
{\bf Extend 2:}  If for the Banach spaces $Y_0,Y_1$ continuously embedded in $X$, and the set $\Sigma_n\subset X$, we know that %
\be 
\nonumber
E(U(Y_0),\Sigma_{n})_X\le \e_n, \quad E(U(Y_1),\Sigma_n)_X\le \tilde  \e_n,
\ee 
then it follows that for $Z:=(Y_0,Y_1)_{\theta,q}$,
$0<\theta< 1$ and $\ 0<q\le \infty$, 
we have
\be 
\nonumber
E(U(Z), \bar \Sigma_n)_X\le C\e_n^{1-\theta}\tilde \e_n^\theta,  
\ee 
where $\bar \Sigma_n:=\{aS+bT:\ a,b\in\R;\ S,T\in\Sigma_n\}$, and $C$ depends only on $\theta$.
\vskip .1in
We prove only {\bf Extend 1} since the proof of {\bf Extend 2} is similar.  We can take $q=\infty$ because this is the largest space $Z$ for the given $\theta$. If $f\in U(Z)$, for $t=\e_n$,  there is a  $g\in Y$
 (if the infimum is not achieved, the proof follows from some limiting arguments) 
 that satisfies
\be 
\label{K1}
\|f-g\|_X+\e_n \|g\|_Y\le K(f,\e_n;X,Y) \le  \e_n^\theta.
\ee 
We know that there is an $S\in  \Sigma_n$ which approximates $g$ to accuracy $\e_n\|g\|_Y$.  For this $S$, we have
\be 
\label{K2}
\|f-S\|_X \le \|f-g\|_X+\|g-S\|_X\le K(f,\e_n;X,Y)\le \e_n^\theta.
\ee

Here is a simple but typical example of {\bf Extend 1}.  If we establish a bound $\e_n$ for approximation of functions in $U({\rm Lip}\ 1)$ with error measured in $X=C(\Omega)$, then we automatically get the bound $\epsilon_n^\alpha$ for approximating functions from $U({\rm Lip}\ \alpha)$, $0<\alpha<1$,
because Lip $\alpha=(C(\Omega),{\rm Lip}\ 1)_{\alpha,\infty}$. 
%

\section{Evaluation of  nonlinear methods of approximation}
\label{S:optimalperformance}

Before embarking on an analysis of the approximation performance of ReLU networks, we wish to place this type of approximation into
the usual setting of approximation theory, and thereby draw out the type of questions that should be answered.  As we have noted, approximation using the outputs of  neural networks with a fixed architecture   is a form of nonlinear approximation known as manifold approximation.  Given a target function $f$ in a Banach space $ X$, the approximation is given by $A_n(f):=M_n(a_n(f))$, where the two maps
$$
a=a_n:X\to\R^n, \quad M=M_n:\R^n\to X,
$$
select the  $n$ parameters of the network and output the approximation, respectively.

Of course, there are many methods of approximation.   The question we  address in this section is how could we possibly determine if approximation by NNs is in some quantifiable sense superior to other more
traditional methods of approximation.  Also, what are the inherent limits on the capacity of NNs to approximate, once  the number $n$ of parameters allocated to the approximation is fixed?
To answer such questions, we introduce   various traditional ways to compare approximation methods and  say with certainty whether an approximation method is optimal
among all methods of approximation,  or perhaps among all approximation methods with a specified structure.   How NNs do under such methods of comparison is not the subject of this section.  That topic is dealt
with in later sections of this paper.

To begin the discussion, we take the view that an {\it approximation method} is a sequence
\be
\nonumber
\{0\}=:\Sigma_0 \subset \Sigma_1\subset \cdots \Sigma_n\subset \cdots \subset X
\ee
of  nested   sets to be used in approximating functions $f$ from the Banach space $X$  in the norm $\|\cdot\|_X$.  Here $n$, in some sense,  measures the complexity of $\Sigma_n$.  The typical spaces $X$ used in practice are the spaces $L_p(\Omega).$ However, at this point, we let $X$ be any Banach space of functions on $\Omega$ with a norm $\|\cdot\|_X$.

The various methods of approximation are divided into two general categories: linear   and nonlinear.  A method is said to be linear if, for each $n$,  the set  $\Sigma_n$ is a linear space of dimension $n$, that is, $\Sigma_n$ is the linear span of $n$ elements from $X$.
The standard examples are spaces of polynomials, splines, and wavelets.   Note that the term linear
does not refer to how the approximation  depends on  $f\in X$.  It only refers to the structure of each $\Sigma_n$, $n\ge 0$.
All other methods of approximation are referred to as nonlinear.   For nonlinear methods,  a linear combination of elements from $\Sigma_n$ may not lie in $\Sigma_n$.   
There are three prominent examples of nonlinear approximation we wish to mention.    
  
  In the first,  $\Sigma_n$ consists of  piecewise polynomials (of a fixed and generally small degree $r$)  on a domain $\Omega\subset \R^d$. Let us  denote by $\cP_r$ the linear space of polynomials of degree $r$.  Here,  we can use any notion of degree in $d$ variables,  such as coordinate degree or total degree.  Given $n$, an element $S\in \Sigma_n$ is obtained by partitioning the domain $\Omega$ into $n$ disjoint cells  $\cC_j\subset \Omega$, $j=1,\dots,n$, and   assigning a  polynomial $P_j\in \cP_r$ to each cell.  Thus, we have
\be
\nonumber
S=\sum_{j=1}^nP_j\chi_{\cC_j},
\ee
where $\chi_{\cC_j}$ is the characteristic function of the cell $\cC_j$. The partitions are not fixed but allowed to vary   within a class of partitions that can be described by $n$ parameters. We have already seen an example of this
in the case of free-knot CPwL functions of one variable, in which case the partition was allowed to consist of any $n$ intevals.  In the multivariate case, the allowable partitions are more structured and usually generated adaptively.
The rough idea of this form of approximation is to use small cells where the target function is rough and large cells where the function is smooth.  

From our description of the sets $\Upsilon^{W,L}(\Relu;d,1)$, we see that NN approximation fits into the above framework of  piecewise polynomial 
approximation in the sense that  each element in one of these sets is a CPwL function on a polytope partition, see \S \ref{S:Relunetworks}.  However,  several notable distinctions arise.
First of all, we have many fewer restrictions on the partitions that arise when compared to other piecewise polynomial methods of approximation such
as FEMs, adaptive methods, free-knot splines, etc.  Another important point is that $\Upsilon^{W,L}(\Relu;d,1)$ is not the collection of all
CPwL functions subordinate to a fixed class of partitions.  Here, choosing the parameters of the network specifies in tandem the partition and
the CPwL.  One view of how this is done is nicely explained  in \cite{baraniuk}.

Another widely used example of nonlinear approximation is    {\it $n$-term approximation}.  Let $B:=\{\phi_j,\, j\geq 1\}$ be an unconditional  basis for $X$.  The set $\Sigma_n:=\Sigma_n(B)$ in this case consists of all functions $S\in X$ which are a linear combination of at most $n$ of  these basis elements.  Thus, each $S\in \Sigma_n(B)$ takes the form
\be
\nonumber
S=\sum_{k\in\Lambda} \alpha_k\phi_k,\quad \#(\Lambda)=n,\quad \alpha_k\in \R,
\ee
where $\Lambda\subset \N$ is a subset of $n$ indices. The index set $\Lambda$ is allowed to change at each occurrence with the only restriction being  that $\#(\Lambda)=n$. One can  generalize this setting  if $B$ is replaced by a frame or, more generally, a dictionary.  Typical examples used in numerical analysis
and signal/image processing are dictionaries of wavelets, curvelets, ridge functions, shearlets, and other families of waveforms.  In this generality, $n$-term approximation is not numerically implementable because the dictionary is infinite.  To circumvent this in practice, one uses a large but finite dictionary that is sufficiently rich for the problem at hand.

Neural network approximation fits most naturally into   a third type of nonlinear approximation known as {\it manifold approximation}.  In manifold approximation, the elements of $S\in\Sigma_n$ take the form  
  \be
 \nonumber
    S=M_n(y),\quad y\in\R^n,
    \ee
 where $M_n:\R^n\to X$, that is,  $\Sigma_n=\{M_n(y):\ y\in\R^n\}$.  
 As noted earlier, a numerical implementation of manifold approximation is made by 
 specifying a mapping 
 $a_n:K\to \R^n$    which,  when presented with $f\in K$, describes the parameters of the point on the manifold used to approximate $f$. 
 
 Given an approximation method  $\Sigma:=(\Sigma_n)_{n\geq 0}$ and  $f\in X$,  we let
\be
\nonumber
E_n(f)_X:=E(f,\Sigma_n)_X:=E_n(f,\Sigma)_X:= \inf_{g\in\Sigma_n} \|f-g\|_X,
\ee
denote the {\it error of approximation of $f$ by elements from $\Sigma_n$}.  Note that $E_n(f)_X$ gives the smallest error we can achieve using  $\Sigma_n$ to do the approximation, but it does not address the question of how to find such a best or near best approximation.  This is an important issue, especially for NN approximation,  that we address later in this section.

An often quoted property of NNs is their {\it universality}, which means that $E_n(f)_X\to 0$ as $n\to\infty$, for all $f\in X$.  This is a 
property possessed by all approximation methods used
in numerical analysis.  Universality is  not at all special and certainly cannot be used to explain the success of NNs.

 \subsection{Approximation of model classes}
 \label{SS:modelclasses}
 
 We do not measure the performance of an approximation method on a single function $f$  but rather on a class $K\subset X$ of functions. In this case, we have the {\it class error}
 \begin{equation}
 \label{approxerrorK}
 \nonumber
 E_n(K)_X:=E(K,\Sigma_n)_X:=E_n(K,\Sigma)_X:=\sup_{f\in K} E_n(f,\Sigma)_X, \quad n\geq 0.    
 \end{equation}
 Here, $K$  incorporates the knowledge we have about the function or potential functions  $f$ that we are trying to capture.   For example, when numerically solving a PDE, $K$ is typically provided by a regularity theorem for the PDE. In the case of signal processing, $K$ summarizes what is known or assumed about the underlying signal, such as bandlimits or sparsity.

  Note that $E_n(K)_X$ represents the worst case error.  It is also possible to measure error in some
 averaged  sense.  This would be meaningful, for example, when the set  $K$ is given by a stochastic process with some underlying probability measure.  For now, we discuss only the worst case error.
 
 A set $K$ on which we wish to measure the performance of an approximation method is called a {\it model class}. We always assume that $K$ is a compact subset of $X$. If the approximation process is universal,  then
 $E_n(K)_X\to 0$ as $n\to \infty$ for every model class $K$.   How fast it tends to zero represents how good the sets $(\Sigma_n)_{n\geq 0}$ are for approximating the elements of $K$.
 
 If we are presented with approximation processes given by $\Sigma=(\Sigma_n)_{n\geq 0}$ and $\Sigma'=(\Sigma'_n)_{n\geq 0}$ respectively,  then given a model class $K$, we can compare the performance of these methods on $K$ by checking the decay
 of $E_n(K,\Sigma)_X$ and $E_n(K,\Sigma')_X$  
 as $n\to \infty$.  If the decay rate of $E_n(K,\Sigma)_X$   is faster than that of $E_n(K,\Sigma')_X$ as $n\to\infty$, we are tempted to say that  $\Sigma$ is superior to  $\Sigma'$
 at least on this model class. However, the question of the computability  of the approximant is an important issue and has to be taken into consideration.

 To drive home this latter point, the following example is   germane.  Given  a compact set $K\subset X$ and  $\e>0$,
 let $S_1=S_1(\e)$ be a finite subset of $K$ such that ${\rm dist}(f,S_1)_X\le \e$ for all $f\in K$. For example, $S_1$ could be the set of centers of an $\e$ covering of $K$. 
Going  further, we can find a one dimensional manifold $\Sigma_1$ that is parameterized by $t\in[0,1]$ and passes through each point in $S_1$ as $t$ runs through $[0,1]$, and thus $E(K,\Sigma_1)_X\leq \varepsilon$. The point of this simple observation is to emphasize that we must place some further restrictions on what we allow as an approximation method
$(\Sigma_n)_{n\ge 0}$ so that we can have a meaningful theory.  What such restrictions should look like and what are  their implications  is the subject we address next.

 \subsection{Widths for measuring approximation error}
 \label{SS:widths}
 
  The concept of widths was introduced to quantify the best possible performance of approximation methods on a given model class $K$.  The best known  width is  the Kolmogorov width, which was introduced to quantify the best possible approximation when using linear spaces.   If $X_n$ is a linear subspace of $X$ of dimension $n$,
 then its performance in approximating the elements of the model class $K$ is given by the error $E(K,X_n)_X$ defined in \eqref{approxerrorK}. 
  If we fix the value of $n\ge 0$, the Kolmogorov $n$-width of $K$
 is defined as
 \be
 \label{Kwidth}
 d_0(K)_X:=\sup_{f\in K} \|f\|_X, \quad d_n(K)_X:=\inf_{\dim(Y)=n} E(K,Y)_X, \quad n\geq 1,
 \ee
 where the infimum  is taken over all linear spaces $Y\subset X$ of dimension $n$. An $n$ dimensional  space    which achieves the infimum in \eqref{Kwidth} is called a Kolmogorov space for $K$ if it exists.   

The Kolmogorov $n$-width of a model class $K$  tells us the optimal performance possible  for approximating  $K$ using linear spaces of dimension $n$ for the approximation. It does not tell us how to select
 a (near) optimal space $Y$ of dimension $n$ for this purpose nor how to find a good/best approximation from $Y$ once it is chosen.  In recent years, discrete optimization methods have been discovered for finding optimal subspaces.  They go by the name of greedy algorithms, see \cite{BMPPT}, \cite{BCDDPW}, 
 \cite{DPW}.   If $X$ is a Hilbert space and $Y$ is a finite dimensional subspace, then we can always find the best approximation from $Y$ to a given $f\in X$ by orthogonal projections.  This becomes a problem when $X$ is a general Banach space because linear projections onto a general $n$ dimensional space $Y$ may have large norm when $n$ is large.  Although a famous
 theorem of Kadec-Snobar says that there is always a projection with norm at most $\sqrt{n}$, finding such a projection is a numerical challenge.
 Also, projecting onto such a linear space does not give the best approximation from the space because 
 the norm of the projection is large.
 
For classical model classes such as the finite ball in  smoothness spaces like the Lipschitz, Sobolev, or Besov spaces,  the Kolmogorov widths are known asymptotically when $X$ is an $L_p$ space.  Furthermore, it is often known  that  specific linear spaces of dimension $n$ such  as polynomials,
 splines on uniform partition, etc.,   achieve this optimal asymptotic performance (at least within reasonable constants).  This can then be used to show that for such $K$, certain numerical methods, such as spectral methods or FEMs are also asymptotically optimal among all possible choices of numerical methods built on using linear spaces of dimension $n$ for the approximation.
 
Let us note that in the definition of Kolmogorov widths  we are not requiring that   the mapping which sends $f\in K$ into the approximation to $f$  is a linear map.    There is a concept of {\it linear width} which requires the linearity of the approximation map.  Namely, given $n\ge 0$ and a model class $K\subset X$, its {\it linear width} $d_n^L(K)_X$ is defined as
    \be
    \label{linwidth}
    d_0^L(K)_X:=\sup_{f\in K} \|f\|_X, \quad d_n^L(K)_X:=\inf_{L\in \cL_n} \sup_{f\in K} \|f-L(f)\|_X,\quad n\geq 1,
    \ee
where the infimum is taken over the class $\cL_n$ of all linear maps from $X$ into itself with rank at most $n$.  The asymptotic decay of linear widths for classical smoothness classes are also known. We refer the reader to the books \cite{Pinkusbook},
\cite{Lorenzbook} for the fundamental results for Kolmogorov and linear widths. When $X$ is not a Hilbert space,  the linear width  of $K$ can decay worse than the Kolmogorov width.
 
\begin{figure}[h]
  \centering
\includegraphics[scale=.5]{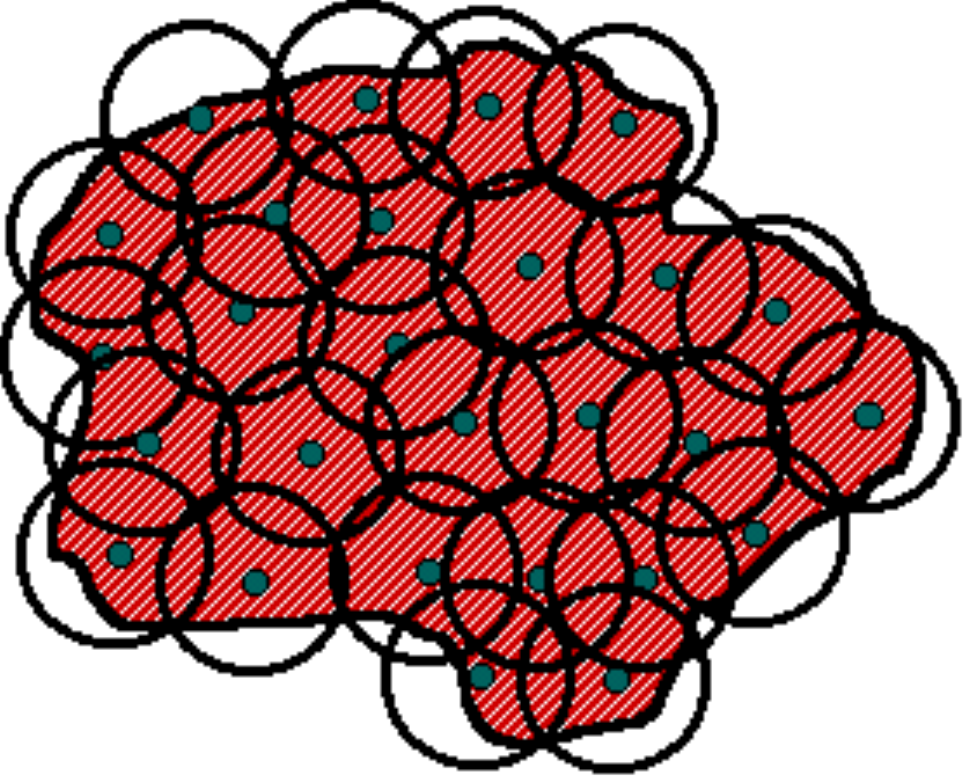}
\caption{Kolmogorov covering of $K$.}
\label{Figentropy}
\end{figure}

 Now, we want to make a very   important point.    There is a general lower bound on the decay of   Kolmogorov widths that was given by Carl in \cite{C}.
 This lower bound can be very useful in showing that a linear
 method of approximation is nearly optimal. To state this lower bound, we need to introduce the Kolmogorov  {\it entropy} of a compact set $K$. Given $\e>0$, compactness says that $K$ can be covered by a finite number of balls of radius $\e$, see Figure \ref{Figentropy}.  We define the covering number $N_\e(K)_X$ to be the  smallest number of balls of radius $\e$ that cover $K$, and we define  the entropy $H_\e (K)_X$ of $K$ to be the logarithm of this number
  \be
 \nonumber
  H_\e (K)_X:= \log_2(N_\e(K)_X).
  \ee
The entropy of $K$ measures  how compact the set $K$ is and is often used to give lower bounds on how well we can approximate the elements in $K$ and also how well we can learn an element from $K$ given data observations.  The  Kolmogorov entropy of a compact set is an important quantity for measuring optimality, not only in approximation theory and numerical analysis, but also in statistical estimation and encoding of signals and images.

  To formulate the lower bounds for widths in terms of entropy, we introduce the related  concept of entropy numbers.   Given $n\ge 0$, we  define the {\it entropy number} $\e_n(K)_X$ to be the infimum of all $\e>0$ for which $2^n$ balls of radius $\e$ cover $K$, that is,
  \be
  \nonumber
  \e_n(K)_X:=\inf\{\e: N_\e(K)_X\le 2^n\}.
  \ee
  The decay rate of entropy numbers for  all classical smoothness spaces in $L_p(\Omega)$ are known. 
  
  Carl proved that for each $r>0$, there is a constant $C_r$, depending only on $r$, such that 
 \be
 \label{carl}
 \Lambda:=\sup_{m\ge 0} (m+1)^rd_m(K)_X<\infty\quad \Rightarrow\quad \e_n(K)_X\le C_r\Lambda(n+1)^{-r}, \quad n\geq 0.
 \ee
 Thus,  for polynomial decay rates for approximation by $n$ dimensional linear spaces, this decay rate  cannot be better than the decay rate for  the entropy numbers of $K$.  Let us note that for many  standard model classes $K$, such as finite balls in Sobolev and Besov spaces, the decay rate of $d_n(K)_X$ is
 much worse than $\e_n(K)_X$.   A version of Carl's inequality holds for other decay rates, even exponential, and can be found in \cite{CDPW}.

   \subsection{Nonlinear widths}
  \label{SS:nlwidths}
  Since NN approximation is a nonlinear method of approximation,  the Kolmogorov widths are not an appropriate measure of performance.  Many different notions of nonlinear widths, see the discussion in \cite{DKLT}, have been  introduced to match the various  forms of nonlinear approximation used in numerical computations.  We shall discuss only nonlinear widths that match the form of approximation provided by NNs.
  
   Recall that NN approximation is a form of manifold approximation, where the approximation set $\Sigma_n$ consists of the outputs of a neural network with $n$ parameters.  Thus, $\Sigma_n=M_n(\R^n)$, with 
   $M_n$ being the mapping that describes how the output function is constructed once the parameters  and architecture of the NN  are set.  Note that   $\Sigma_n$ is a nonlinear set in the sense that the sum of  two elements from $\Sigma_n$  is generally not
   in $\Sigma_n$.  
   
    There are by now numerous papers that discuss the approximation by NNs.  They typically provide estimates for $E(K,\Sigma_n)_X$ for certain model classes $K$.  We will discuss such  estimates subsequently in \S \ref{S:Shallow} and \S\ref{S:deep}.   We have cautioned that such results must be taken with a grain of salt since they do not
   typically discuss how the approximation would be found or numerically constructed.  Our point of view is that it is not just an issue of how well $\Sigma_n$ approximates $K$, although this is indeed an interesting question, but
   also how a good approximation would  be found.  In other words, the parameter selection mapping $a_n$ is equally important.

When presented with an $f\in K$, one chooses the parameters of the NN to be used to construct the approximation  to $f$.     Typical algorithms in learning base this selection of parameters on some form of optimization, executed through gradient descent methods.  For our analysis, we denote this selection procedure by the mapping $a_n:K\to \R^n$.  So the approximation procedure is given by $A_n(f):=M_n(a_n(f))$.  If we wish to establish some form of  optimality of NNs,  we should compare NN approximation with other approximation  methods of this form.

Given any pair of mappings (not necessarily using NNs)  $a:X \to \R^n$ and $M:\R^n\to X$, we define
  {\it the error for approximating $f\in X$}   by
 \be 
 \nonumber
 E_{a,M}(f)_X:=\|f-M(a(f))\|_X,
 \ee
 and the approximation error on a model class $K\subset X$ by
 \be 
 \nonumber
 E_{a,M}(K)_X:=\sup_{f\in K} E_{a,M}(f).
 \ee
 For any such $a,M$ we have
 \be
 \label{compare11}
 E(K,\Sigma_n)_X\le E_{a,M}(K)_X,\quad \hbox{where}\quad \Sigma_n:=M(\R^n),
 \ee
 and we have equality when we choose $a(f)$ so that $M(a(f))$ is a best approximation to $f$ (assuming such a best approximation exists) from $\Sigma_n$.
   
A first possibility for defining optimal performance of such methods of manifold approximation on a model class $K$ would be to simply  find the minimum of $E_{a,M}(K)_X$ over all such pairs of mappings.
 However, we have already pointed out that this minimum would always be zero (even when $n=1$)  because of the existence of  space 
  filling  manifolds.  On the other hand, these space filling manifolds are   useless in numerical analysis.      
  Consider, for example, a one parameter space filling manifold. By necessity, a small perturbation of the parameter will generally  result in a large change in the output, which  makes parameter selection
  for fitting $f$ impossible.  The natural question that arises is  what restrictions need to be imposed on the mappings $a,M$ so that we have a theory which corresponds to reasonable numerical methods. We discuss this next.

\subsection{Restrictions on $a,M$ in manifold approximation}
\label{SS:restrictions} The first suggestion, given in \cite{DHM}, for the possible restrictions to place on the mappings $a,M$, was to require that they be continuous.
 This led to the following definition of  {\it manifold} widths $\delta_n(K)_X$,
 \be
 \nonumber
 \delta_n(K)_X:=\inf_{a_n,M_n} E_{a_n,M_n}(K)_X,
 \ee
 with the infimum taken over all  maps $a_n:K \to \R^n$ and $M_n:\R^n \to X$, where $a_n$ is continuous on $K$ and $M_n$ is continuous on $\R^n$.   Manifold widths are closely connected to other definitions of nonlinear widths, see the discussion in  \cite{DKLT}.  
 
 It turns out that even with these very modest assumptions on the mappings $a,M$, one can prove lower bounds for $\delta_n(K)_X$ when $K$ is a unit ball of a classical smoothness space, e.g. Besov, Sobolev, Lipschitz,  and these lower bounds show that manifold approximation is no better than other methods of nonlinear approximation such as $n$-term wavelet approximation or adaptive finite element approximation  for these model classes.  For example, if we approximate in $L_p(\Omega)$, with $\Omega=[0,1]^d$, and $K$ is a unit ball of any Besov space that embeds compcatly into $L_p(\Omega)$, then it was show in  \cite{DKLT} that
 \be 
 \label{Sobwidth}
 \delta_n(K)_{L_p(\Omega)}\ge Cn^{-s/d}, \quad n\ge 1.
 \ee 
   This should not be used
 to deduce that  manifold approximation, in general, and NN approximation, in particular,  offer nothing new in terms of their ability  to approximate. It may be that their power to approximate lies in their ability  to handle non-traditional model classes.  Nevertheless, this should make us proceed with caution.

  A stark criticism of manifold widths is that its requirement
 of continuity of the mappings is  too minimal and does not correspond to the notions of numerical stability used in practice.  In other words,  manifold approximation based on just assuming that $a,M$ are continuous may also not  be implementable
 in a numerical setting.  We next discuss what may be more viable restrictions on $a,M$ that match numerical practice.

\subsection{Stable manifold widths}
\label{SS:stablewidths}

 A major issue in the implementation of a method of approximation is its stability, that is, its sensitivity to computational error or noisy inputs.    The stability we want can be summarized in the following two properties:
 \vskip .1in
 \noindent
 {\bf (S1)}  When we input $f$ into the algorithm, we often  input a noisy discretization of $f$, which can be viewed as a perturbation of $f$.  So, we would like to have the property that when $\|f-g\|_X$ is small,
  the algorithm outputs $M(a(g))$  which is close to $M(a(f))$.   A standard quantification of this is to require that the mapping $A:=M\circ a$ is a Lipschitz mapping from $K$ to $X$.  Notice that in this formulation
 the perturbation $g$ should also be in $K$.
 
 \vskip .1in
 \noindent
 {\bf (S2)}  In the numerical implementation of the algorithm, the parameters $a(f)$ are not  computed exactly, and so we would like  that when $a,b\in \R^n$ are close to one another, then $M(a)$ and $M(b)$ are likewise close.
 Again, the usual quantification of this observation  is to impose that  $M:\R^n\to X$ is a Lipschitz map.  This property requires the specification of a norm on $\R^n$ which is controlling the size of the perturbation of $a$.
 \vskip .1in

  The simplest way to guarantee {\bf (S1)-(S2)} is to require that both  mappings $a,M$ are Lipschitz, which means that  there  is a norm $\|\cdot\|_Y$ on $\R^n$  and a  number $\gamma\ge 1$, such that 
 \be
 \label{Lipschitz1}
 \|a(f)-a(g)\|_Y \le \gamma \|f-g\|_X,\quad f,g\in K,
 \ee
 \be
 \label{Lipschitz2}
 \|M(x)-M(x')\|_X\le \gamma \|x-x'\|_Y, \quad x,x'\in\R^n.
 \ee
 If $a,M$ satisfy \eref{Lipschitz1}-\eref{Lipschitz2},  then obviously  {\bf (S1)} and {\bf (S2)} hold, where the Lipschitz constant in  {\bf (S1)} is $\gamma^2$.

Imposing Lipschitz stability on  $a,M$ leads to  the following  definition of {\it stable manifold} widths   %
 \be
 \nonumber
 \delta_n^*(K)_X:=\delta_{n,\gamma}^*(K)_X:= \inf_{a,M} E_{a,M}(K)_X,
 \ee
 where now the infimum is taken  over all maps  $a:K\to \R^n$ and $M:\R^n\to X$ that are Lipschitz continuous with constant $\gamma$.

\subsection{Bounds for stable manifold widths}
\label{SS:stablerate}

Both upper and lower bounds for   stable manifold widths of a compact set $K$ are given in \cite{CDPW}.  These bounds are tight in the case when the approximation takes place
in a Hilbert space.  Approximation in a Hilbert space is often used in  applications of NNs.     

Lower bounds for the decay of stable manifold widths in a general Banach space $X$  are given by the following Carl's type inequality, see \eqref{carl}, which compares $\delta_n^*(K)_X$ with the entropy numbers $\e_n(K)_X$. Specifically, for any $r>0$, we have
  \be
 \label{Carl}
 \e_n(K)_X\le C(r,\gamma)  (n+1)^{-r}\sup_{m\ge 0}(m+1)^r\delta^*_{m,\gamma}(K)_X, \quad n\ge 0.
 \ee
This shows that whenever the stable manifold widths $\delta_n^*(K)_X$ of a model class $K$ tend to zero like ${\cal O}(n^{-r})$, $n\to\infty$, then the entropy numbers of $K$ must have the same or faster rate of decay. 
Similar bounds are known when the decay rate $n^{-r}$, $n\to\infty$, is replaced by other decays, see \cite{CDPW}.  In this sense, the stable manifold widths  $\delta^*_{n,\gamma}(K)_X$ cannot tend to zero faster than the entropy numbers of $K$.

 The inequalities \eqref{Carl} give a bound for  how well manifold approximation can perform on a model class $K$ once Lipschitz stability of the maps $a,M$ is imposed.  One might speculate, however, that in general $\e_n(K)_X$ may go to zero   faster than $\delta_n^*(K)_X$.   This is not the case  when $X=H$ is a Hilbert space, since in that case  for any compact set $K\subset H$, we have the estimate
 \be
 \label{compare1}
 \delta_{26n,2}^*(K)_H\le 3 \e_n(K)_H,\quad n\ge 1,
 \ee
proved in \cite{CDPW}. 
This is a very useful information since it is often relatively easy to compute the entropy
numbers  of a model class $K$. In addition, it is also a very useful result for our, yet to come, analysis of NN approximation. 

The upper bound \eref{compare1} is proved through three fundamental steps.  The first is to select $2^n$ points $S_n:=\{f_i\}_{i=1}^n$ from $H$ such that the balls of radius $\e:=\e_n(K)_H$ centered at these points cover $K$.
The next step is to use the Johnson-Lindenstrauss dimension reduction lemma to find a (linear) mapping $a:S_n\to\R^{26n}$,  for which 
\be
\nonumber
\frac{1}{2}\|f_i-f_j\|_H\leq \|a(f_i)-a(f_j)\|_{\ell_2(\R^{26n})}\leq  \|f_i-f_j\|_H, \quad i,j=1,\ldots,2^n.
\ee
 According to the Kirszbraun extension theorem, see Theorem 1.12
from \cite{BL}, the mapping $a$ can be 
extended from $S_n$
to the whole $H$, preserving the Lipschitz constant 1.
The last step is to define 
$M$ on $a(f_j)$, $j=1,\ldots,2^n$, as
\be
\nonumber
M(a(f_j))=f_j,\quad j=1,\ldots,2^n.
\ee
Clearly 
$$
\|M(a(f_i))-M(a(f_j))\|_H=\|f_i-f_j\|_H\leq 2\|a(f_i)-a(f_j)\|_{{\ell_2(\R^{26n})}},
$$
and therefore $M$ is a Lipschitz map with a Lipschitz constant $2$ when restricted to the finite set $a(S_n)$. Again, according to the Kirszbraun extension theorem, we can extend $M$ to a Lipschitz map on the whole $\R^{26n}$ with the same constant $2$.

It is now easy to see that the approximation operator   $A:=M\circ a$ gives the desired approximation performance since,
with a suitable choice of $j$, we have
\begin{eqnarray}
\nonumber
\|f-A(f)\|_H &\leq &\|f-f_j\|+\|M(a(f_j))-M(a(f))\|_H\\ \nonumber
&\leq& \varepsilon_n(K)_H+2\|a(f)-a(f_j)\|_{\ell_2(\R^{26n})}\\
\nonumber
&\leq& \e_n(K)_H+2\|f-f_j\|_H\leq 3\varepsilon_n(K)_H.
\end{eqnarray}
Therefore, we have proved \eqref{compare1}.  Let us remark however that $A$ is not very constructive and that 
it is generally difficult to create Lipschitz mappings $a,M$ that achieve the optimal performance in stable nonlinear widths.

\subsection{Weaker  measures of stability}
\label{SS:other}

 It may be argued that requiring Lipschitz stability is too strong of a requirement. Recall that Lipschitz stability is just a sufficient condition to guarantee the stability properties {\bf (S1)-(S2)} that we want. In this direction, we mention that   {\bf (S1)-(S2)} will hold if $a,M$ satisfy the following weaker properties with $\|\cdot\|_Y$ some fixed norm on $\R^n$:

\vskip .1in
\noindent
{\bf (SP1)} The mapping $A:=M\circ a$, $A:K\to X$ is Lipschitz.  We can even weaken this further to requiring only $\|A(f)-A(g)\|_X\le C\|f-g\|_X^\alpha$, $f,g\in K$, for some $\alpha\in (0,1]$.  This is known as Lip $\alpha$ stability.

\vskip .1in
\noindent
{\bf (SP2)} The mapping $M:\R^n\to X$ is Lipschitz, or more generally, Lip $\alpha$ with respect to   $\|\cdot\|_Y$ on $\R^n$.
\vskip .1in
\noindent
While not directly needed for stability, the following property will play a role in our further discussions.
\vskip .1in
\noindent
{\bf (SP3)}  The mapping 
$a:K\to\R^n$ is bounded with respect to  $\|\cdot\|_Y$ on $\R^n$.  This property limits the search over parameter space.
\vskip .1in

It turns out that if the mappings $a,M$ satisfy the weaker assumptions {\bf (SP1)-(SP3)}, then one can still prove a version of Carl's inequality, and thus,  we still have limitations on the performance
of these approximation methods in terms of entropy number lower bounds. Let us briefly indicate how lower bounds for performance are proved when the mappings $a_n,M_n$ satisfy {\bf (SP1)-(SP3)}. For notational simplicity only,  we take $\alpha=1$, the Lipschitz constant of both $M_n$ and $M_n\circ a_n$ to be 
$\gamma$, and the image of $K$ under $a_n$ to be contained in the unit ball of $\R^n$ with respect to  $\|\cdot\|_Y$.

We fix $\e>0$ and 
let $a_n:K\rightarrow \R^n,\, M_n:\R^n\rightarrow X$, satisfy {\bf (SP1)-(SP3)} and approximate the elements of $K$ with the accuracy
\be
\label{errorbound1}
E_{a_n,M_n}(K)_X\le \e/3,\quad 
\hbox{for some $\e>0$}.
\ee 
We now show that \eref{errorbound1}  implies a bound on the entropy of $K$.
Let us denote by $\mathrm{Pack}_\e:=\set{f_1,\ldots,f_{P_\e(K)}}$ a maximal $\e$-packing of $K$, that is, a collection of 
points $\{f_i\}\in K$,
 with $\min_{i\neq j}\|f_i-f_j\|_X>\e$, whose size is maximal among all such collections.  
Now,  
define
\[y_i := a_n(f_i),\quad g_i := \lr{M_n\circ a_n}(f_i)=M_n(y_i),\qquad i=1,\ldots, P_\e(K).\]  It follows that for $i,j=1,\ldots, P_\e(K)$,
\be 
\nonumber
\|g_i-g_j\|_X\ge \|f_i-f_j\|_X-\|f_i-g_i\|_X-
\|f_j-g_j\|_X> \e/3,\quad i\neq j,
\ee 
where we used \eref{errorbound1}.  Since $M_n$ is $\gamma$ Lipschitz, we obtain
 \[\norm{y_i-y_j}_{Y} \geq \frac{1}{\gamma}\norm{M_n(y_i)-M_n(y_j)}_X=\frac{1}{\gamma}\norm{g_i-g_j}_X > \frac{\e}{3\gamma},\quad i\neq j.\]
 In other words, since 
 $\|y_i\|_Y=\|a_n(f_i)\|_Y\leq 1$, the collection $\set{y_1,\ldots, y_{P_{\e}(K)}}$ is an $\frac{\e}{3\gamma}$-packing of the unit ball in $\R^n.$ Well known volumetric considerations show that a maximal such packing can have at most $(1+6\gamma\e^{-1})^n$ elements,
and therefore
$P_{\e}(K)\leq \lr{1+6\gamma\e^{-1}}^n. 
$

Now, since the balls of radius $\e$ centered at the $f_i$ are a covering of $K$, we have that
\be 
\label{boundcovnumber}
N_{\e}(K)\leq P_{\e}(K)\leq \lr{1+6\gamma\e^{-1}}^n = 2^{n\log_2\lr{1+6\gamma\e^{-1}}}.
\ee 
For example, the above derivation shows that whenever there are mappings $a_n,M_n$ satisfying {\bf (SP1)-(SP3)}, then we have the Carl inequality 
\be 
\label{Carl2} 
E_{a_n,M_n}(K)_X\le Cn^{-r},\ n\ge 1\implies \e_n(K)_X\le C' n^{-r}[\log_2 n]^r,\ n\ge 1.
\ee 
Indeed, we take $\e= 3Cn^{-r}$ and use \eref{boundcovnumber} to find 
$\e_{cn\log_2n}\le Cn^{-r}$ which gives \eref{Carl2}.

\subsection{Optimal performance for  classical model classes described by smoothness}
\label{SS:classicalwidths}
 Although the definition of manifold widths places very mild conditions on the mappings $a,M$, it still turns out that these conditions are sufficiently strong to restrict how fast $\delta_n(K)_X$ tends to zero for model classes
 built on classical notions of smoothness
 described by smoothness conditions such as Sobolev or Besov regularity.  For example, if $B_q^s(L_\tau(\Omega))$, with $\Omega=[0,1]^d$, is any Besov space that lies above the Sobbolev  embedding line for $L_p(\Omega)$, then it is proven in \cite{DKLT} that
\be
\nonumber
\delta_n(U(B_q^s(L_\tau(\Omega))))_{L_p(\Omega)} \asymp \e_n(U(B_q^s(L_\tau(\Omega))))_{L_p(\Omega)} \asymp n^{-s/d},\quad n>0,
\ee
with the constants in this equivalence depending only on $d$.

It turns out that the decay rate ${\cal O}(n^{-s/d})$ can be obtained by many   methods of nonlinear approximation such as adaptive finite elements or $n$-term wavelet approximation.
The main message for us is that even with this mild condition of imposing only continuity on the maps $a,M$,  we cannot do better than the rate ${\cal O}(n^{-s/d})$ for these classical smoothness classes when using
manifold approximation.  In particular, this holds for NN approximation with the restriction of continuity on  the mappings $a,M$ associated to  the NNs. 

\subsection{VC dimension also limits approximation rates for model classes}
\label{SS:VClimits}
The results we have given above provide  lower bounds on how well a model class $K$ can be approximated by a 
stable manifold approximation.   If we remove the requirement of stability, it is still possible to give lower bounds on approximation rates
for model classes if the approximation method $(\Sigma_n)_{n\ge 0}$ is made up of sets $\Sigma_n$ which have limited VC dimension.  For such results, one needs some additional assumptions on  the model class $K$.  We describe results of this type in this section.

  Suppose $K$ is a model class in $L_p(\Omega)$ with $1\le p\le\infty$.   A common technique in proving  lower bounds on the Kolmogorov entropy or  widths of $K$ is to exhibit a function $\phi\in L_p(\Omega)$ with compact support 
  for which the normalized dilate
  \be
  \label{nd}
  \Phi(x ):= A\phi(\lambda x),\quad x\in \Omega,
  \ee
 is in $K$, provided $A$ and $\lambda$ are chosen appropriately.  The function $\Phi$ is called a {\it bump} function. By choosing $\lambda$ large,
  one concentrates the support of $\Phi$ but of course this is at the expense of making $A$ small in order to guarantee that the resulting $\phi$ is in $K$.
 The small support of $\Phi$ guarantees  that the shifted functions  $\Phi_i(\cdot)=\Phi(\cdot-x^{(i)})$, $i=1,\dots,N$, are also in $K$ and these functions have disjoint supports, provided $N$ is not too large and
  the $x^{(i)}$'s are suitably spaced out in $\Omega$.  It then follows that for any assignment of signs $\Lambda:=(\e_1,\e_2,\ldots,\e_N)$, $\e_i= \pm 1, \ i=1,\dots,N$, the function
  \be
  \label{bumps1}
  f_\Lambda:=B\sum_{i=1}^N\e_i\Phi_i,
  \ee
  is also in $K$ for a proper choice of $B$.  One then uses the rich family of functions $f_\Lambda$ as $\Lambda$ runs over the $2^N$ sign patterns to show that the Kolmogorov
  entropy of $K$ must be suitably large.

  This strategy can be used to bound from below how well a model class can be approximated by sets with limited VC dimension.
  For illustration,  we  consider  the simplest example where $K=U(C^r(\Omega))$, with $r$ being a positive integer, and measure approximation  error in the norm $\|\cdot\|_{C(\Omega)}$.   If   we  approximate the functions in $K$    by using a set $\cF$ 
  with $VC(\cF)\leq m$,
then we claim that there is a  constant $C=C(r,d)>0$ such that 
 \be
\label{never}
\delta:={\rm dist}(K,\cF)_{C(\Omega)}\ge C m^{-r/d}.
\ee
We prove this claim
 in the case $\Omega=[-1,1]^d$. Consider  a non-negative bump function   $\phi\in C^r(\R^d)$ which vanishes outside $[-1/2,1/2]^d$ and has norm  
 $\|\phi\|_{C(\Omega)}=\phi(0)$, see Figure \ref{Figbump}. The dilated function $\Phi$ of \eqref{nd}  is in 
$K$ if 
we choose $A$ so that $A\phi(0)+A|\phi|_{C^r(\Omega)}\lambda^{r}=1$.  The support of $\Phi$    is  contained in a $d$ dimensional cube centered at $0$ with side-length $\lambda^{-1}$.  So, if  {$N=\lfloor \lambda^{d}\rfloor$, we can make the $\Phi_i$'s appearing in \eqref{bumps1} to have disjoint support
by taking  the $x^{(i)}$  suitably separated. Moreover,  if $B=1$,
each of the $f_\Lambda\in K$.

\begin{figure}[h]
  \centering
\includegraphics[scale=.75]{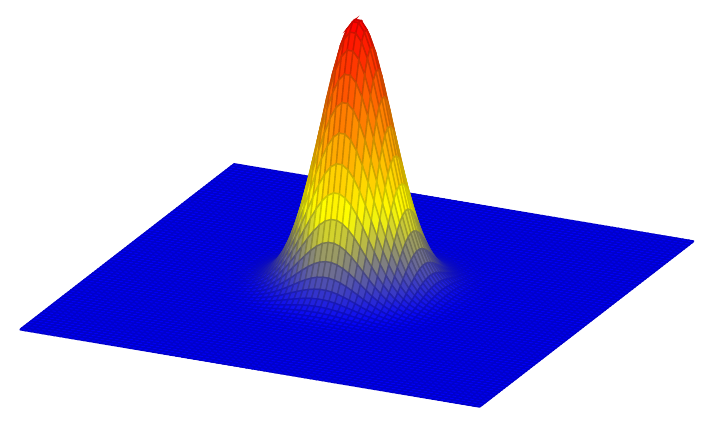}
\caption{A bump function in 3D}
\label{Figbump}
\end{figure}

  Now, to prove \eqref{never}, we take $\lambda= \lceil (m+1)^{1/d}\rceil $ and
obtain $N\ge m+1$ functions $\Phi_i(\cdot):=\Phi(\cdot-x^{(i)})\in K$, $i=1,\dots,N$,  with disjoint supports. Then, for each choice of sign patterns the function $f_\Lambda$ from \eqref{bumps1} is in $K$ and  
$f_\Lambda (x^{(i)})= A\phi(0)\e_i$, $i=1,\dots,N$.   Now $f_\Lambda$ is approximated by an $S_\Lambda\in\cF$ to accuracy
$\delta$.
  If $\delta$ were smaller than $A\phi(0)$,  then the function $S_\Lambda$ would carry the sign pattern of the $\e_i$ at each $x^{(i)}$.
    Hence, the points $x^{(i)}$, $i=1,\dots,N$,  would be shattered by $\cF$.
Since by assumption $VC(\cF)\leq m$, this is not possible, and we must  have $\delta \ge A\phi(0)$. 
Since we have that  $\phi(0)A=\phi(0)(\phi(0)+\lambda^r|\phi|_{C^r(\Omega)})^{-1}\geq Cm^{-r/d}$,
this proves \eqref{never}.

This argument can also be used to prove  that there is an absolute constant $C>0$ depending only on $s$, such that
for 
$K=U(B_q^s(L_\infty(\Omega)))$, with  $s>0$, $0<q\le \infty$,  
we have
\be
\label{VCBesov}
{\rm dist}(K,\cF)_{C(\Omega)}\ge C m^{-s/d}, 
\ee
whenever the VC dimension of $\cF$ is at most $m$.
We leave the proof to the reader.
 
Let us now specialize to the case where $\cF$ is the output of a ReLU network.   With an eye towards our bounds on VC dimension for
the spaces $\Upsilon^{W,1}(\Relu;d,1)$, see Lemma \ref{L:pseudo1}, and $\Upsilon^{W_0,L}(\Relu;d,1)$, see Theorem \ref{L:VCdeep}, we obtain  the following lower bounds for NN approximation for
$W\geq 1$ and $d\geq 2$,
\be
\label{lbW}
 {\rm dist}
 (U(B_q^s(L_\infty(\Omega))),\Upsilon^{W,1}(\Relu;d,1))_{C(\Omega)}\ge C(s,d)[W\cdot\log_2W]^{-s/d},
\ee
and 
\be
\label{lbL}
 {\rm dist}(
 U(B_q^s(L_\infty(\Omega))),\Upsilon^{W_0,L}(\Relu;d,1))_{C(\Omega)}\ge C(s,d)L^{-2s/d}.
\ee
 The logarithm in \eref{lbW} can be removed when $d=1$.

Notice that in the case of \eqref{lbW}, the lower bound can be stated as $C(s,d)[n\log_2n]^{-s/d}$, where $n=n(W,1)$ is the number of parameters used to describe
$\Upsilon^{W,1}$.  Thus, in this case, save for the logarithm, we cannot achieve any better approximation rates than that obtained by traditional linear methods of 
approximation.  We discuss later in \S \ref{S:Shallow}   what rates have been proved in the literature for one layer networks.

In the case of \eqref{lbL}, the lower bound is of the form 
$C(s,d)n^{-2s/d}$, where $n=n(W_0,L)$ is the number of parameters used to describe
the space $\Upsilon^{W_0,L}$.  The factor $2$ in the exponent leaves open the possibility of much improved approximation rates (when compared
with classical methods) when using deep networks.  We shall show in \S \ref{SS:super}   that these rates of approximation are attained.

We close this section by mentioning that the use of VC dimension to bound approximation rates from below seems to be
restricted to the case when approximation error is measured in the norm $\|\cdot\|_{C(\Omega)}$.  This makes one wonder if there is
a concept analogous to  VC dimension suitable for $L_p$ approximation when $p\neq \infty$.

\subsection{Another measure of optimal performance: approximation classes}
\label{SS:aclasses}

There is another important way to measure the performance of an approximation method $\Sigma=(\Sigma_n)_{n\ge 0}$ by 
looking at the set of all functions which have a given approximation rate
as $n\to\infty$.  Let $\lambda=(\lambda_n)_{n\ge 0}$ be a sequence of positive real numbers which decrease monotonically to zero.   We define
\be
\label{aclassdef}
\cA(\lambda):= \cA(\lambda,\Sigma) :=\{f\in X\,:\,\exists \Lambda >0\,\,E(f,\Sigma_n)_X\le \Lambda\lambda_n,\,\,\forall n\geq 0\},
\ee
and further define $\|f\|_{\cA(\lambda)}$ as the smallest number $\Lambda$ for which \eqref{aclassdef} holds.
The larger this set is, the better the approximation method $\Sigma$ is.  

The case when $\lambda_n:=(n+1)^{-r}$ is the most often studied since it corresponds to the rates most often encountered in numerical scenarios.  In this case,
$\cA(\lambda)$ is usually denoted by 
\be 
\nonumber
\cA^r=\cA^r(\Sigma)=\cA^r(\Sigma,X).
\ee 
A major chapter in approximation theory is to characterize the approximation classes $\cA^r$ for a given approximation method.  The main theorems of approximation theory
characterize $\cA^r$ for polynomial and spline approximation.  Such characterizations are also known for some methods of nonlinear approximation.  

As we shall see, we are far from understanding the approximation classes $\cA^r$ for NN approximation.  However, some useful 
results on the structure of these classes can be found in \cite{gribonval2019approximation}.

\vskip .2in

\section{Approximation using ReLU networks: overview}
\label{SS:ReLUefficiency}
As we have already noted, the collection   $\Upsilon^{W,L}=\Upsilon^{W,L}(\Relu;d,1)$ of  outputs of a ReLU network  with width $W$,  depth $L$, and input dimension $d$ is a nonlinear set of CPwL functions determined by $n(W,L)$ parameters.  Our interest in
the next few sections is to summarize the approximation power  of $\Upsilon^{W,L}$.  In the process of analyzing this, we shall not    address the question of whether there is a practical stable algorithm to produce  the
approximation, an issue we will discuss in  \S \ref{S:StableRelu}.

One of the impediments
to giving a coherent presentation of the approximation properties of the outputs of neural  networks, as the number of parameters increases, is the great variety of possible architectures of the networks.  Namely, when examining the approximation efficiency, we can fix $W$ and let $L$ change,  or fix $L$ and let $W$ change, or let  both change simultaneously. We can also vary the architecture by allowing full connectivity or sparse connectivity  between layers. We may also impose further structure on the weight matrices, leading, for example, to  convolution networks.  Moreover, we can as well  consider a variety of activation functions $\sigma$.

While each such setting is of interest,  we primarily concentrate    on  two cases  of ReLU networks.  The first is the case that most
closely matches classical approximation, the set $\Upsilon^{W,1}$ as $W\to\infty$. We shall see that  even this case is not completely understood.   At the other extreme is the case when we take the width $W$ to be some fixed constant $W_0$ and let  $L\to\infty$.  This is  a  most illuminating setting  in that
we shall see a  dramatic gain in approximation efficiency when the depth $L$ is allowed to grow.  This is commonly referred to as {\it the power of depth}. 

To provide a unified notational platform, we use  $\Sigma_n$ for the set $\Upsilon^{W,L}$ under consideration, where $n$ is equivalent to the number of parameters being used.  For example, we can take $\Sigma_n=\Upsilon^{n,1}$ or $\Sigma_n=\Upsilon^{W_0,n}$ since both of these sets depend on a number of parameters proportional to $n$.  Our goal is to understand how the family $\Sigma:=(\Sigma_n)_{n\ge 0}$  performs as an approximation tool.

In what follows in this section, we consider the set $\Upsilon^{W,L}$ restricted to the domain $\Omega:=[0,1]^d$. Recall that each function $S\in \Upsilon^{W,L}$ is the output of a neural network with at most $ n(W,L)= (d+1)W+W(W+1)(L-1)+ (W+1)$ parameters. We consider the error of approximation to be measured in an $L_p(\Omega)$ norm, $1\le p\le\infty$. 
Therefore, for $f\in L_p(\Omega)$, we are interested in the error of approximation
\be
\nonumber
E(f,\Sigma_n)_{L_p(\Omega)}:= E_n(f,\Sigma)_{L_p(\Omega)}:=\inf_{S\in \Sigma_n}\|f-S\|_{L_p(\Omega)},
\ee
when $\Sigma_n$ is one of the nonlinear sets $\Upsilon^{W,L}$ and $n\asymp n(W,L)$.  In the case $p=\infty$, we   assume that $f$ is continuous and the error is measured in the $\|\cdot\|_{C(\Omega)}$ norm, and so the results hold uniformly in $x\in\Omega$.

Note  that using $L_p(\Omega)$, $1\le p\le \infty$,  norms to measure error does not match the usual measures of  performance of classification algorithms, where 
the main
criteria is probability or expectation of misclassification, see \cite{Lugosi}.  This is an important distinction that we unfortunately will not address because of a lack of definitive results.
It may be that this distinction is in fact  behind the success of NNs in the learning environment.

The results we prove can be extended to approximation in 
$L_p(\Omega)$ for $0<p<1$, but this requires some technical effort we want to avoid.
We concentrate on the three most important cases $p=\infty$ (the case of uniform approximation), the case $p=2$ which is prevalent in stochastic estimates, and the case $p=1$ which monitors average  
error.   We always take 
the $L_p(\Omega)$  spaces with Lebesgue measure.   Let us also remark that the results we derive hold equally well
for general Lipschitz domains taken in place of $\Omega=[0,1]^d$. If we fix the value of $p$, the results we seek are of the    following two types.

\vskip .1in
\noindent
{\bf Model class peformance:} For a model class $K\subset L_p(\Omega)$, we have earlier defined
\be
\nonumber
E(K,\Sigma_n)_p:= E_n(K,\Sigma)_{L_p(\Omega)}:= \sup_{f\in K} E(f,\Sigma_n)_p.
\ee
Our interest is to describe the decay of this error (with estimates from above and below) as $n\to\infty$.
\vskip .1in
There are two types of model classes $K$  that are of interest.  The first are classical smoothness classes such as the unit ball of a Lipschitz, H\"older Sobolev, or Besov spaces, see \S\ref{S:modelclasses}.   In this way, we can compare
the approximation properties of NNs with more standard methods of approximation and see whether NNs offer better performance on these classical model classes.

A second type  of results of  interest  is to uncover new model
classes $K$ for which NNs perform well and classical methods of approximation do not.  Such new model classes would help clarify exactly when NN approximation is beneficial.
 Motivation for these  new model classes should come from  the intended application of NN approximation.
Such results might explain why NNs perform well in these applications.

\vskip .1in
\noindent
{\bf Characterization of Approximation Classes:} 
A second category of results that is  of interest would be to understand the approximation classes $\cA^r(\Sigma,L_p(\Omega))$ for NN approximation.    Recall that  these classes,  see \S\ref{SS:aclasses}  for their definition, consist of all functions $f$ whose approximation error satisfies 
\be
\label{AC1}
E(f,\Sigma_n)_{L_p(\Omega)}\le M (n+1)^{-r},\quad n\ge 0,
\ee
with the smallest $M$ defining 
$\|f\|_{\cA^r}$.

We would like to know which functions are in $\cA^r$. While a  precise characterization of these classes is beyond our current understanding of NN approximation,  the results that follow give sufficient conditions for a function $f$ to be in such a class.    In contrast, for many types of classical approximation, both linear and nonlinear, there are characterizations of their corresponding  approximation classes.  Such characterizations 
require what are called {\it inverse theorems} in approximation theory. An inverse theorem is a statement that whenever $f\in \cA^r$, we can prove that
$f$ is in a certain Banach space $Y_r$.

Consider, for example,  the case of approximation in $L_p(\Omega)$.   An inverse theorem is proved by showing an inequality of the form
\be 
\nonumber
|S|_{Y_r}\le C(n+1)^r\|S\|_{L_p(\Omega)},\quad S\in\Sigma_n, \quad n\ge 0.
\ee 
For example, if we consider approximation by   trigonometric polynomials of degree $n$ in one variable, in the metric $L_p([-\pi,\pi])$,
  one inequality of this type is the famous Bernstein inequality for trigonometric polynomials
\be 
\nonumber
\|T'\|_{L_p(\Omega)}\le n\|T\|_{L_p(\Omega)},
\ee 
which holds for any trigonometric polynomial of degree $n$.  So $r=1$ in this example, and $Y_1=W^1(L_p([-\pi,\pi]))$.

Such inverse theorems are not known for NN approximation save for the case of $\Sigma=(\Upsilon^{n,1} (\sigma;1,1))_{n\geq 0}$
for certain activation functions $\sigma$, including ReLU. Thus, there is quite a large gap in our understanding of NN approximation as compared to these more classical methods.   It is of major interest to establish  inverse inequalities for the elements in $\Sigma_n$ when $\Sigma_n$ is a set of outputs of a NN.

\section{Approximation   using single layer ReLU networks}
 \label{S:Shallow}

 In this section, we study approximation on the domain $\Omega:=[0,1]^d$ by the family $\Sigma:=(\Sigma_n)_{n\ge 0}$, where $\Sigma_0:=\{0\}$ and $\Sigma_n:=\Upsilon^{n,1}(\Relu;d,1)$, $n=1,2\dots$, with input dimension $d\geq 1$.  These sets are rarely used in numerical settings since deeper networks are preferred. However, for theoretical reasons, it is important to understand their approximation properties in order to see the advantages of the deeper networks studied later in this paper. Much of the activity on NN approximation has been directed at understanding the approximation properties of these single hidden layer networks.  Surprisingly, we shall see that  most  fundamental questions about approximation using $\Sigma$ are  not yet answered.

 We have discussed in \S\ref{S:d>1-L=1} the structure of $\Sigma_n$.  Each function $S\in \Sigma_n$ is a CPwL function in $d$ variables $x=(x_1,\dots,x_d)$ of the form
\be
\label{pwpd}
S(x):= b_0+\sum_{j=1}^n a_j \eta_j(x),\quad \eta_j(x):=(w_j\cdot x+b_j)_+,
\, w_j\in \R^d, \,
b_j,a_j\in \R,
\ee
where $\eta_j$ is linear on the half space
  $H_j^+:=\{x:\, w_j\cdot x+b_j>0\}$  and is zero otherwise.   Note that $S\in \Sigma_n$ is a CPwL function subordinate to a hyperplane partition $\cP$ of $\R^d$ into cells which are convex  polytopes.  
  
 In spite of the simplicity of the representation \eref{pwpd}, the set $\Sigma_n$ is quite complicated save for the case $d=1$, see \S \ref{SSS:shallowd1}.  First of all,
 the  possible partitions $\cP$ that arise from hyperplane arrangements are complex in the sense that the cells are  not isotropic, the  number of cells can be quite large,  and there is not a  simple characterization of these partitions. This is compounded by the fact that, as we have previously discussed, not every  CPwL function subordinate to  a partition given by an arrangement of $n$ hyperplanes is in $\Sigma_n$.  For example, this set does not contain any compactly supported functions. This is in contrast to the typical applications of CPwL functions in numerical PDEs.  Thus, $\Sigma_n$ is a complex, but possibly rich nonlinear family.  We shall see that this complexity inhibits our understanding of its approximation properties.

  Keeping in mind the discussion in the previous section, there are three types of results that we would like to prove in order to understand the approximation power of $\Sigma:= (\Sigma_n)_{n\ge 0}$, measured in the  $\|\cdot\|_{L_p(\Omega)}$ norm, $1\le p\le \infty$.  
  \vskip .1in
  \noindent
  {\bf Problem 3: }
 {\it  Give matching  upper and lower bounds for  $E_n(K,\Sigma)_{L_p(\Omega)}$ when $K$ is one of  the classical model classes   such as unit balls of Lipschitz, H\"older, Sobolev, and Besov spaces.}
  \vskip .1in
  \noindent 
We shall see that, save for the case $d=1$,  this problem is far from being solved.

As we have previously stressed, the partitions generated by hyperplane arrangements are complex and not well understood, with  cells that are possibly highly anisotropic.  This suggests the possibility of being able to approximate functions which are not described by classical isotropic smoothness   and leads us to expect new model classes that are well approximated by $\Sigma$.
  \vskip .1in
  \noindent
  {\bf Problem 4:} {\it  Describe  new model classes $K$ of functions that are guaranteed to be well approximated by $\Sigma$.}
  \vskip .1in
  \noindent
  Some advances on {\bf Problem 4} have been made, centering on the so-called Barron classes that we discuss in \S\ref{SS:novel}.   
  
  Finally, the most ambitious approximation problem for $\Sigma=(\Sigma_n)_{n\ge 0}$ is
  the following.
  \vskip .1in
  \noindent
  {\bf Problem 5:} {\it  For each $r>0$ and $1\le p\le\infty$, characterize the approximation class $\cA^r(\Sigma,L_p(\Omega))$ consisting of all functions $f\in L_p(\Omega)$ for which
  $$
  E_n(f,\Sigma)_{L_p(\Omega)}={\cal O}((n+1)^{-r}),\quad n\ge 0.
  $$}
  \noindent 
  Nothing is known on this last problem when $d>1$, and we are skeptical that any definitive result is around the corner for the   case of general $d$.
  
  In order to orient us to the type of results we might strive to obtain on these  problems for  general $d$, we begin in the next section by discussing the case $d=1$, where we have the most
   extensive results and the best understanding of approximation from these spaces.
  
  \subsection {Approximation by single layer networks when $d=1$}
  \label{SS:d=1}
   We begin by discussing the case $d=1$ not only because it is the best understood, but also because it can orient the reader as to what we can possibly expect when engaging the case $d>1$.  Because approximation by $\Upsilon^{n,1}(\Relu;1,1)$ is essentially the same as free-knot linear spline approximation, results for the NN approximation are derived from the known results on free-knot splines. The latter  are  well explained in \cite{devore1998nonlinear} and the literature cited therein, and summarized below.

  \subsubsection{Approximation of classical model classes when $d=1$}
  \label{SSS:d1classical}
  Here, we measure approximation error in $L_p(\Omega)$ with $1\le p\le \infty$ and domain  $\Omega=[0,1]$. The classical model classes for $L_p(\Omega)$
  are  finite balls in the Lipschitz, H\"older, Sobolev, and Besov spaces.  The latter  spaces are the most flexible for measuring smoothness and the approximation
  properties, for all of the other smoothness classes can be derived from them.  So, we restrict our discussion to the model classes $K=U(B_q^s(L_\tau(\Omega)))$, $0<q,\tau\le\infty$,
  which have smoothness of order $s>0$.  These spaces were introduced and discussed in \S \ref{SS:Lipschitz}, where we have noted that these spaces are compactly embedded in $L_p(\Omega)$ when
  $s > 1/\tau-1/p$, i.e., when these spaces lie above the Sobolev embedding line, see Figure \ref{fig2.pdf}.  They are not embedded in $L_p(\Omega)$  if they lie below the embedding line.
  
  The following theorem summarizes the results known about approximating Besov classes in  the case $d=1$.
  \begin{theorem}
  \label{T:d1classical}
  Let  $K=U(B^s_q(L_\tau(\Omega)))$ be  the unit ball of the Besov space $B^s_q(L_\tau(\Omega))$.  If $0<s\le 2$ and this space lies above the Sobolev embedding line for $L_p(\Omega)$ then
  \be 
  \label{T:d1}
  E_n(K,\Sigma)_{p}\le C(s,p,\tau)(n+1)^{-s},\quad n\ge 0.
  \ee 
  \end{theorem}
  
  Let us elaborate a little on what this theorem is saying.  First, note that the sets $K$ for which we obtain the approximation rate ${\cal O}((n+1)^{-s})$
  allow the smoothness describing $K$ to be measured in $L_\tau(\Omega)$, where $\tau\neq p$.  When $\tau\ge p$, the result does not need to exploit
  the nonlinearity of $\Sigma_n$ in the sense that the approximation rate can be obtained already by using linear spaces  corresponding to fixing the breakpoints in $\Sigma_n$ to be equally spaced on $[0,1]$.  It is only when $\tau<p$ that we need to exploit nonlinearity.
  
  A couple of simple examples may be in order.  Consider approximation in $C(\Omega)$ and smoothness of order  $s=1$.   Obviously, the space Lip 1 is compactly embedded in $C(\Omega)$ and the approximation rate is ${\cal O}((n+1)^{-1})$,   $n\to\infty$, when $K=U({\rm Lip} 1)$.  Note that  Lip $1$ is not a Besov space but  is continuously  embedded in $B_\infty^1(L_\infty(\Omega))$ and the latter space is covered by the theorem. Hence Lip $1$ also is.
  We can obtain the approximation rate ${\cal O}((n+1)^{-1})$ by taking the breakpoints equally spaced and thereby using a linear subspace of $\Sigma_n$.
     The Sobolev space $W^1(L_1(\Omega))$ is also contained in $C(\Omega)$, but not compactly.  Nevertheless, its unit ball has the approximation rate 
     ${\cal O}((n+1)^{-1})$.  The Sobolev spaces $W^1(L_p(\Omega))$, $p>1$, have unit balls that are compact in $C(\Omega)$ and the theorem gives that they also  have the approximation rate ${\cal O}((n+1)^{-1})$, $n\to \infty$.
      Recall that for $f$ to be in  Lip 1 requires that it has bounded derivative
  $\|f'\|_{L_\infty(\Omega)}< \infty$, while  $f\in W^1(L_p(\Omega))$ only requires $f'\in L_p(\Omega)$.  For example, the function  $f(t)=t^\alpha$, $0<\alpha<1$,  is in $W^1(L_p(\Omega))$ if $p>1$ is small enough, but this function 
  is  not in Lip 1.  The way one gets good approximation of $t^\alpha$ by $\Sigma_n$ is to put half of the  breakpoints of the output $S\in \Sigma_n$ near $0$ and the remaining half equally spaced in $\Omega$.  Thus, for these Sobolev spaces one truly needs the nonlinearity of $\Sigma_n$.
  To achieve the optimal approximation rate, we need to choose the breakpoints to depend on $f$, and thus we cannot choose them in advance.
  
  Finally, let us remark why we have the restriction $s\le 2$.  We are approximating locally by linear functions. A function $f$ with smoothness of
  order $s>2$ would need to use locally polynomials of degree higher than one to improve its local error of approximation (think of Taylor
  expansions).  Hence, when $f$ has smoothness of order $s>2$, we do not improve on the rate ${\cal O}((n+1)^{-2})$, $n\to\infty$, which we already have for functions with smoothness of order $2$.

  \subsubsection{Approximation classes for $d=1$}
  One of the crowning achievements of nonlinear approximation at the end of the last century was the characterization of the approximation classes for several classical methods of nonlinear approximation, including free-knot spline, $n$-term wavelet, and adaptive piecewise polynomial approximation.
  The key to establishing these results was not only to give upper  bounds for the error in  approximating functions from  Besov spaces but also to prove certain inverse theorems that say if a function $f$ can be approximated with a certain rate ${\cal O}((n+1)^{-r})$, $n\to \infty$, then  $f$ must possess a certain Besov smoothness.  These inverse theorems should not be underestimated since they allow precise characterization of approximation classes. 
  
  In the case of approximation using CPwL functions, the inverse theorems were provided by the seminal theorems of Pencho Petrushev, see  \cite{PP-NL}. 
  The approximation space $\cA^r=\cA^r(\Sigma,L_p(\Omega))$ is precisely characterized, provided $0<r<2$ and $1\le p\le\infty$, with $C(\Omega)$ used
  in place of $L_\infty(\Omega)$ when $p=\infty$.  In this case, $\cA^r$ is a certain interpolation space, see \cite{devore1998nonlinear}. Since we do not want to
  go too deeply into interpolation space theory here, we simply mention that  $\cA^r$ is sandwiched between two  Besov spaces of smoothness order $r$. More precisely, if $0<r<2$, and $1\le p\le \infty$ are fixed, and $\tau^*:= (r+1/p)^{-1}$, then for all $0<q\le \infty$, we have
  \be
  \label{AS1}
  B_q^r(L_{\tau}(\Omega))\subset \cA^r\subset  B_{\infty} ^r(L_{\tau^*}(\Omega)), \quad {\rm whenever} \ \tau>\tau^*. 
  \ee

 Since this result may be difficult to digest  at first glance, we make some comments to explain what these embeddings say. 
 First, recall the relation of Besov spaces to the Sobolev embedding line, see Figure \ref{fig2.pdf}.  For a fixed value of $r$, all spaces
 $B_q^r(L_\tau(\Omega))$ appearing on the left side of the embedding \eref{AS1} are compactly embedded in the space $L_p(\Omega)$, where we are measuring error.
 The left embedding says that any function in one of these spaces is in $\cA^r$, and hence has approximation error decaying at the rate  ${\cal O}((n+1)^{-r})$.  Note
 that these spaces get larger as we approach the embedding line.  The right embedding says that we cannot allow $\tau$ to be smaller than  $\tau^*$; in fact if $\tau$ is smaller than $\tau^*$ we do not even embed into $L_p(\Omega)$.
 Besov spaces that appear on the embedding line itself may or may not be compactly embedded in $L_p(\Omega)$, depending on $q$.  They are compactly embedded
 if $q$ is small enough.  
 
 Finally, we remark that we have the characterization of $\cA^r$ only if $r< 2$ for the same reason we had the restriction $s\le 2$ when discussing approximation of classical model classes in the previous section.  Going a little further,  note that if $S_n\in\Sigma_n$, $n\ge 1$, then $S_n\in\cA^r$ for all $r>0$, but $S_n$ is not in any smoothness space of order $s>2.$
 Moreover,  any function $f=\sum_{k\ge 1} \alpha_k S_k$, $\alpha_k\in \R$,  will be in $\cA^r$, $0<r<\infty$, if $(\alpha_k)_{k\ge 1}$ tends to zero sufficiently fast.  Yet, $f$ will not have any classical smoothness of order $s>2$.  So $\cA^r$, $r>2$, cannot be characterized by classical smoothness such as membership in a Besov space.

  \subsection{Results for $d\ge 2$}
  \label{d2}
 Continuing with one hidden layer networks, let us now  consider the case $d\ge 2$.    The difficulty in  constructing effective approximations in this case is the fact that when $d\ge 2$, the set of NN outputs  $\Sigma_n=\Upsilon^{n,1}(\Relu;d,1)$ does not have locally supported nodal functions that are commonly used to build approximants. So,  approximation
  methods are built on global constructions.  It is not surprising therefore, that the strongest  results are known in the case where the approximation is measured in the $L_2(\Omega)$ norm,
  where orthogonality can be employed in the constructions.  We discuss the $L_2$ approximation first.
  
\subsubsection{Approximation in $L_2(\Omega)$} 
\label{SS:L=1,p=2}
For approximation in $X=L_2(\Omega)$, it is known that when $f\in W^s(L_2(\Omega))$, we have
\be 
\label{d2first}
E_n(f,\Sigma)_{L_2(\Omega)}\le C n^{-s/d} \|f\|_{W^s(L_2(\Omega))},\quad n\ge 1,
\ee
provided $s\le 2+(d-1)/2$.
The case $d=2$ is given in \cite{DOP}, and the general case is considered in \cite{PP-NN}.

We  give a very coarse description of the ideas behind proving \eref{d2first}. In this discussion, it is useful to work with functions defined on the unit Euclidean  ball $\Omega^*\subset\R^d$ rather than on the cube $[0,1]^d$.   One can move between these different domains via restriction and extension operators which are known to preserve Sobolev and Besov regularity.   When  $g$ is a univariate function, then $g(a\cdot x)$,
$a\in\R^d$ is a {\it ridge function} of $d$ variables  (sometimes called a planar wave).  If $X_n$ is a linear space of dimension $n$ of univariate functions and $\Lambda$ is a fixed subset of $m$  unit vectors in
$\R^d$, then the set of functions $Y_{n,m}:=\mathrm{span} \{g(a\cdot x):\ g\in X_n, \ a\in\Lambda\}$ is a linear space of dimension at most $mn$.

The core of the proof of  \eref{d2first} is to show that if $X_n$ is effective in approximating univariate functions in $L_2$, and if the set $\Lambda$ is \lq uniformly distributed'  on the unit sphere in $\R^d$ (the boundary of the unit ball of $\R^d$),
then $Y_{m,n}$ will provide an approximation to $W^s(L_2(\Omega^*))$ functions when choosing $m={\cal O}(n^{d-1})$.  We can take $X_n$ to be the space of univariate CPwL functions on  an equally spaced partition.  The resulting space $Y_{n,m}$ is contained in $\Upsilon^{mn,1}(\Relu;d,1)$, and thereby proves \eref{d2first}.    
The proofs of these results are quite elaborate and technical. 

If we wish to characterize the approximation performance of $\Sigma$ on the model class $K:=U(W^s(L_2(\Omega^*))$, then we would need to establish lower
bounds for the approximation error $E_n(K)_{L_2(\Omega^*)}$ that match those of \eref{d2first}.   Such bounds are plausible but seem not to be known.  However, there are lower bounds for approximating $K$ by general ridge functions given in \cite{Maiorov}, which give for our setting and $d\ge 2$, the lower bound
\be  
\nonumber
E_n(K)_{L_2(\Omega^*)} \ge Cn^{-s/(d-1)},\quad n\ge 1,
\ee 
where $C$ depends only on $d$.

While the results given above are less than satisfactory, because of the lack of matching upper and  lower bounds,  the situation becomes  even  worse when we seek results that show the benefits of the nonlinear 
structure of the sets $\Sigma_n$, $n\ge 1$.
As we know from the case $d=1$, nonlinear methods
of approximation should allow smoothness to be measured in the weaker $L_\tau(\Omega)$ norms while retaining the same approximation order. Namely,
the question is  what are the approximation rates when $K$ is the unit ball of a Besov space $B_q^s(L_\tau(\Omega))$ that is above the Sobolev embedding line
for $L_2(\Omega)$.  In contrast to the case $d=1$,  we do not know results that quantify the performance of $\Sigma$, for the Besov spaces 
that compactly embed into $L_2$.

\subsubsection {Approximation in $L_p$, $p\neq 2$}
\label{SS:L=1,pnot2} 
When we consider approximation in $L_p(\Omega)$, $p\neq 2$, we are only aware of results for $p=\infty$   given in  \cite{Bach}. These are only stated  for 
 the unit ball $K$ of Lip 1 with approximation error measured in the norm of $C(\Omega)$, and take the form
\be 
\nonumber
E_n(K,\Sigma)_{C(\Omega)}\le C \frac{\log_2 n}{n^{1/d}},\quad n\ge 1.
\ee 
Since we are now considering approximation in the space $C(\Omega)$, we can use the known upper bounds for the  VC dimension of $\Upsilon^{n,1}(\Relu;d,1)$
to derive lower bounds on
the approximation error for $K$.  If we apply  Lemma \ref{L:pseudo1}  and employ arguments similar to those we used to prove \eref{never}, we obtain
\be 
\nonumber
E_n(K,\Sigma)_{C(\Omega)}\ge C [n\log_2 n]^{-1/d},\quad n\ge 1.
\ee 
In other words, modulo logarithms, the approximation rate of $K$ is $n^{-1/d}$.  It is of interest to remove these log terms.

We can derive  bounds on the approximation rates  for the model classes $K_\alpha:= U({\rm Lip}\ \alpha)$, 
$0<\alpha<1$, from the known ${\rm Lip} \ 1$ bound by using interpolation theory, see {\bf Extend 1}  in \S \ref{SS:interpolatiooperators},
which gives the bound
\be 
\nonumber
E_n(K_\alpha,\Sigma)_{C(\Omega)}\le C\left [\frac{\log_2 n} { n^{1/d}}\right]^{\alpha}, \quad n\ge 1.
\ee

One expects that these results also extend to error estimates for approximation by $\Sigma_n$ of the unit balls of the smoothness spaces
$B^s_q(L_\infty(\Omega))$ for some  range of 
$s$ larger than one.  However, these do not seem to be  found in the literature.  Equally missing are results for approximation in $L_p(\Omega)$ when $p\neq 2,\infty$.  Moreover, none of the known results reflect the expected gain from the fact that $\Sigma$ is a nonlinear method of approximation.

\subsubsection{Novel model classes for single layer approximation}
\label{SS:novel}
 
As we have noted earlier, there is much interest in identifying new model classes for which NN approximation is particularly effective.  One celebrated model class of this type was introduced by Andrew Barron in  \cite{barron1994approximation}.   This model class and its corresponding approximation results are nicely explained in the exposition \cite{pinkus1999approximation}.  The most recent results on NN appproximation of this class of functions can be found in
\cite{SX} and \cite{KB}.
We    limit ourselves to describing  how these model classes  fit into  the themes  of this article.

Given any domain $\Omega\subset \R^d$, Barron introduced the model class $K=K_\Omega $ consisting of all functions $f\in L_2(\Omega)$ which have an extension to all of $\R^d$ (still denoted by $f$) whose Fourier transform $\hat f$ satisfies
\be 
\label{Barronclass}
\int_{\R^d}
\|\omega\|_{\ell_1(\R^d)}|\hat f(\omega)|\, d\omega\le 1.
\ee 
Notice that \eref{Barronclass} imposes additional conditions over just requiring that  $f$ is square integrable.  Namely, \eref{Barronclass}    requires the decay of $\hat f(\omega)$ as the frequency
$\omega$ gets large. It is easy to check that this is equivalent to requiring that $f$ has a gradient (in the weak sense) whose Fourier transform
is in $L_1$.

Barron initially showed that for any sigmoidal activation function $\sigma$ the approximation family
$\Sigma:= (\Upsilon^{n,1}(\sigma;d,1))_{n\ge 1}$ approximates the model class $K$  in the norm of $L_2(\Omega)$ with the following accuracy

\be 
\label{Barronrate}
E_n(K_\Omega, \Sigma)_{L_2(\Omega)} \le C_\Omega n^{-1/2},\quad n\ge 1.
\ee 
  This result was then shown to hold also  for ReLU activation by using the fact that
$(\Relu(t)-\Relu(t-1))$ is a sigmoidal function. 

Barron's result has spirited a lot of generalizations and applications, and even the introduction of new Banach spaces, see  \cite{EMW}.
Important generalizations of \eref{Barronrate} were given     in 
 \cite{makovoz},  where it was shown that the above result for the class $K_\Omega $ holds for approximation in $L_q$, $1\le q<\infty$, and moreover, the rate of approximation can be improved to ${\cal O}(n^{-1/2-
 1/(q^* d)})$,  where $q^*$ is smallest even integer $\geq q$.  Further improvements on approximation rates for Barron classes and their generalizations have been given through the years.  We refer the reader to \cite{SX} for the latest information.

   We will not dig  too deeply into the known approximation rates for Barron classes and their generalizations  here.  Rather, we  confine ourselves to  some comments to properly frame these results in the context of nonlinear approximation. Let $H$ be a Hilbert space.  We say that a collection $\cD:=\{\phi\}$ of functions from $H$ is a dictionary if each $\phi$ has norm one and whenever $\phi\in \cD$, then so is $-\phi$.  Given such a dictionary $\cD$, we consider the closed convex hull
${\rm co}(\cD)$ of $\cD$.  A fundamental result in approximation theory is that whenever $f\in {\rm co}(\cD)$, then there exists
$g=\sum_{k=1}^n c_k\phi_k$ with the $\phi_k\in \cD$, such that
\be 
\label{fundamentalhilbert}
\|f-g\|_{\mathcal H} \le C n^{-1/2}.
\ee 
There is a constructive method to find such a $g$, known as the {\it orthogonal greedy algorithm}, see \cite{DT}.

To derive \eref{Barronrate} from this, it is enough to show that $K_\Omega$ is contained in the convex hull of the dictionary of
all functions $c\sigma(w\cdot x+b)$ with $w\in \R^d$, $b\in \R$ and $c$ a suitable normalizing constant. The proof of this fact  can be found  in \cite{Barron1} and \cite{pinkus1999approximation}. The improvements of Makovoz rest on the fact that in the case of sigmoidal functions, the dictionary elements $\sigma(w\cdot x+b)$ are very  close to one another when the parameters $w$ and $b$ change slightly and so one can reduce the number of
terms needed in the approximation when seeking an error $\e$.

Notice that neither the constant $C_\Omega$ nor the form of the decay $n^{-1/2}$ in \eref{Barronrate}  depend on $d$. This should be compared with approximation for Sobolev classes where the rates decrease and the constant explodes in size.   However, this must be viewed in the light that  the condition for membership in $K$
gets much stronger as $d$ gets large.  This class is analogous to requiring that $f$ have a Fourier series (in $d$ variables) whose coefficients
are absolutely summable.  Another important point is that the proof of \eref{Barronrate} exploits nonlinear approximation since the $n$ terms from the dictionary $\cD$ used to approximate $f$ are chosen to depend on $f$.

\section{Approximation using  deep ReLU networks}
  \label{S:deep}

  We now study in detail the approximation by the family $\Sigma=(\Sigma_n)_{n\ge 0}$ of  deep networks, with
  $\Sigma_0:=\{0\}$ and $\Sigma_n:=\Upsilon^{W_0,n}(\Relu;d,1)$, $n\geq 1$, where $W_0$ is fixed depending on $d$. The three main conclusions we uncover, following the order of our exposition, are:
  \begin{itemize}
\item  When error is measured in an $L_p$ norm, $1\le p\le \infty$, deep NNs approximate 
functions in the classical model classes (such as Lipschitz, H\"older,
  Sobolev and Besov classes)
at least as well as all of the known methods of nonlinear approximation, see \S\ref{SS:BesovBsplines}.  
 \item  For all classical model classes, deep NN approximation gives error rates dramatically better than all other standard methods of nonlinear approximation, see \S\ref{SS:super}.
 \item  There are novel model classes, built on the ideas of self similarity, where NNs provide approximation rates not available by standard approximation methods, see \S \ref{SS:novel1}.
\end{itemize}

    \subsection{Results obtained from basic decompositions}
    \label{SS:basicdeep}
    
    In this section, we describe what is perhaps the most common method of obtaining estimates for deep NN approximation.
    It is based on two principles.  The first is to show that the target function has a   decomposition in terms of fundamental building blocks with a control on the coefficients in the decomposition. These building blocks could be wavelets or some of their mathematical cousins,  such as shearlets or ridgelets, or they could be global representations  such as power series or Fourier decompositions. 
    For functions $f$ in classical smoothness spaces, we often know the existence of  such decompositions   with  quantifiable bounds on the coefficients of $f$.
    The second step is then to show that each of these building blocks can be captured very efficiently (usually with exponential accuracy) by deep networks. 
    
    These two  principles can then be put together in order to give quantifiable performance for approximation using deep NNs.  This technique appears often in the literature.  A partial list of 
    prominent papers using this method are 
    \cite{yarotsky2017error},
     \cite{OSZ}, \cite{BGKP}, 
     \cite{GRK},
     \cite{grohs2019deep},
     \cite{petersen2018optimal},
     \cite{petersennotes},
     \cite{wang2018exponential},
     \cite{Shencomp},
     \cite{lu2020deep}.

   We formalize the above mentioned  procedure by considering any 
   Banach space $X$ and representing $f\in X$ as  $f=\sum_{k\ge 1}\alpha_kg_k$, where the $\alpha_k$'s are scalars,  $g_k\in X $, and $\|g_k\|_X=1$. Then, we can  bound the   error in approximating $f$ by   its partial sum  by
  \be 
  \nonumber
  \|f-\sum_{k\leq n} \alpha_k g_k\|_X\le \sum_{k>n} |\alpha_k|.
  \ee 
  We can exploit this simple observation in the context of neural networks as follows.  If $g_k\in\Upsilon^{W_0-(d+1),n}(\Relu;d,1)$, $1\le k\le n$, then
  \be 
  \label{lc1} 
  E_{n^2}(f,\Sigma)_X\le   \sum_{k>n}|\alpha_k|,
  \ee 
since the partial sum $\sum_{k=1}^n \alpha_kg_k\in\Sigma_{n^2}
=\Upsilon^{W_0,n^2}(\Relu;d,1)$, see  {\bf Addition by increasing depth} in \S \ref{SS:highlighted}.   
  
  The bound \eref{lc1} is quite crude and can be improved in many ways.  For example, we can give a better control on depth needed,  when each $g_k$ is a composition of the same univariate  function $T$.  We shall use this fact in what follows and so we formulate it in the following proposition.
  \begin{proposition}
\label{P:compsum}
 If  $T\in\Upsilon
 ^{W_0-1,L_0}(\Relu;1,1)$, 
 then any linear combination $S= \sum_{i=1}^m \alpha_i T^{\circ i}\in \Upsilon^{W_0,mL_0}(\Relu;1,1)$.   \end{proposition}
{\bf Proof:} Let $\cN$ be a neural network with width $W_0-1$ and depth $L_0$,  with input and output dimension one,  whose output function is $T$. We concatenate $\cN$ with itself $(m-1)$ times to obtain the network $\cN^*$ of width $W_0$ and depth $mL_0$.  Note that the $kL_0$-th layer of $\cN^*$ can output
$T^{\circ k}$.
We add one collation channel to $\cN^*$, whose nodes pass value zero until layer $(L_0+1)$, where its  node 
collects $\alpha_1T$. This value is then passed forward  until layer $2L_0+1$, where $\alpha_2T^{\circ 2}$ is added, so that $\alpha_1T+\alpha_2 T^{\circ 2}$  is now  held in the node of this channel for layers, $2L_0+1,\dots 3L_0$.
 We continue in this way.  Then,
 we output $S$ from  the $mL_0$-th layer.\hfill $\Box$

  In deriving an estimate like \eref{lc1},  it is not necessary to assume that the functions $g_k\in \Sigma_n$, $k=1,\ldots,n$, but merely that the $g_k$'s  are approximated
  sufficiently well by $\Sigma_n$, as we see in the next proposition.  
  \begin{proposition}
  \label{P:superposition}
  If $g_k$, $k=1,\ldots,n$,  can be approximated by outputs $\hat g_k$  from $\Upsilon^{W_0-(d+1),n}(\Relu;d,1)$ with error
  $\|g_k-\hat g_k\|_X\leq \e$, 
  then the function $\hat S:=\sum_{k=1}^n \alpha_k\hat g_k\in \Upsilon^{W_0,n^2}(\Relu;d,1)$, and 
  \be 
  \label{lc1new} 
  E_{n^2}(f,\Sigma)_X\le \|f-\hat S\|_X\leq  \e\sum_{k=1}^n|\alpha_k|+\sum_{k>n}|\alpha_k|.
  \ee 
  \end{proposition}
 
  \noindent
  {\bf Proof:} The error estimate \eref{lc1new} follows from the fact that 
  $$
  \|f-\sum_{k=1}^n \alpha_k\hat g_k\|_X\leq \|f-\sum_{k=1}^n \alpha_k g_k\|_X+\sum_{k=1}^n |\alpha_k|\| g_k-\hat g_k\|_X.
  $$
  The network that outputs $\hat S$
  is obtained the same way as described above.

 \hfill $\Box$

  In this section, we shall use the following theorem.
  \begin{theorem} 
  \label{T:translate-dilate}
  Let $\varphi\in\Upsilon^{W_0,L_0}(\Relu;d,1)$, $A_j$ be a  $d\times d$ matrix,  and  $b_j\in\R^d$,   $j=1,\dots,n$.  Then, the function
  \be 
  \label{dilates}
  S=\sum_{j=1}^n c_j\varphi(A_jx+b_j), \quad c_j\in\R,
  \ee 
  is in $\Upsilon^{d+1+W_0,nL_0}(\Relu;d,1)$.
  \end{theorem}
  
  \noindent
  {\bf Proof:} Let $\cN_0$ be the network which outputs $\varphi$,
   and let us denote by $A^*$ the 
  $W_0\times d$ matrix of 
   input weights of $\cN_0$, and by $b^*$ the biases of its first layer.

   We build a special network $\cN$ with width 
  $W=d+1+W_0$ and depth $L=nL_0$ to output $S$. Its first $d$ channels are source channels to push forward $x_1,\dots,x_d$.
  The next $W_0$ channels will be the channels of $\cN_0$, and the final channel will be a collation channel to form the sum defining $S$.
  
   The network $\cN$ consists of  $n$ copies of $\cN_0$ placed next to each other. We feed the source channels to the $j^{th}$ copy of $\cN_0$, $j=1,\ldots,n$. For this copy we use input matrix $A^*A_j$ and bias $(b^*+A^*b_j)$ for its first layer. The nodes of the collation channel forward zeroes up to layer $L_0+1$, where the output $c_1\varphi(A_1x+b_1)$ of the first copy of $\cN_0$ is entered  and then forwarded.
   The output of the $j^{th}$ copy is multiplied by $c_j$
   through modification of the output weights of $\cN_0$,  and forwarded to the
   $(jL_0+1)^{st}$ node of the collation channel if $j<n$, where it is added to the current sum in that channel and then the result is forwarded.  When $j=n$ the output of the $n$-th copy is  outputted together with the content of the collation channel to produce $S$. 
  \hfill $\Box$

 \subsection{Approximation of products}
 \label{SS:products}
  We turn next to showing how to approximate certain simple building blocks with exponential accuracy using deep ReLU networks.
   These building blocks include monomials, polynomials, tensor products, and B-splines.  An important tool in establishing  such results is
   to show how one can approximate products of functions, which is our next item of interest.  
   
 Let $H$ be the hat function introduced in \eref{hat}.    We begin with the well known formula \footnote{It is not clear who was the first to observe this formula,
 but it appears already in \cite{Hata}.}
\be
\label{power1}
t(1-t)=\sum_{k\ge 1} 4^{-k} H^{\circ k}(t),\quad t\in [0,1]. 
\ee
  We define
\be
\nonumber
S(t):=t^2\quad {\rm and} \quad S_n(t) :=t-\sum_{k= 1}^n 4^{-k} H^{\circ k}(t),\quad n\ge 1, \quad t\in [0,1].
\ee
Let us note that we can also represent $S_n$ by
\be
\label{power111}
 S_n(t) :=t^2+\sum_{k=n+1}^\infty  4^{-k} H^{\circ k}(t),\quad n\ge 1.
\ee
These two representations of $S_n$ show that  
\be
\label{op}
t^2\le S_n(t)\le t, \quad t\in [0,1],
\ee
and so $S_n:[0,1]\to [0,1]$.

We now prove the following univariate result.
\begin{proposition} 
\label{Recoversquare}
For each $n\ge 1$, the function $S_n\in\Upsilon^{4,n}(\Relu;1,1)$ and satisfies
\be 
\label{squareerror}
\|S-S_n\|_{C([0,1])}\le \frac{1}{3}\cdot 4^{-n},\quad n\ge 1,
\ee 
and 
\be 
\label{Lipproduct}
\|S'-S_n'\|_{L_\infty([0,1])}\le 2^{-n},\quad    n\ge 1.
\ee 
\end{proposition} 
\noindent
{\bf Proof:}  
Since $H\in \Upsilon^{2,1}(\Relu;1,1)$, 
in view of Proposition \ref{P:compsum}, there is a ReLU network
of width $3$ and depth $L=n$ that outputs $\sum_{k=1}^n4^{-k}H^{\circ k}$.
  If we add one more channel to push forward $t$, then we can also output $S_n(t):=t-\sum_{k=1}^n 4^{-k}H^{\circ k}(t)$. 

Since  $S(t)-S_n(t)=-\sum_{k=n+1}^\infty4^{-k}H^{\circ k}(t)$, the bound  \eref{squareerror} 
follows from
\be 
|S(t)-S_n(t)|\le \sum_{k=n+1}^\infty 4^{-k}\le \frac{1}{3}\cdot 4^{-n},\quad t\in [0,1],
\ee 
whereas  \eref{Lipproduct} follows from the fact  that each $H^{\circ k}$ has Lipschitz norm $2^k$.  \hfill $\Box$

Let us mention that there are many functions other than $t^2$  for which explicit formulas like \eref{power1} hold.  These will be discussed in \S \ref{SS:novel1}.
For now,   we   want to examine how we can capture higher order monomials from the above results.  First, we begin by showing  how we can implement multiplication  using deep ReLU networks. We start with  the simple formula
  \be 
  \label{product1}
  \Pi(x_1,x_2):=x_1x_2=
 2S\lr{\frac{x_1+x_2}{2}}- \frac{1}{2}\{ S(x_1)+S(x_2)\},
 \quad x_1,x_2\in [0,1].
  \ee 
   We can construct a neural network with input dimension $2$ which outputs the function 
   $\Pi(\cdot,\cdot)$ with high accuracy.

For $n\ge 1$, we define for $x_1,x_2\in [0,1]$ the function 
\be
\label{PP}
    \Pi_n(x_1,x_2):= 2S_n\lr{\frac{x_1+x_2}{2}}- \frac{1}{2}\{ S_n(x_1)+S_n(x_2)\},
\ee
and  prove the following properties of 
$\Pi_n$.
\begin{proposition} 
\label{Pin}
For each $n\ge 1$, $\Pi_n(x_1,x_2)\in [0,1]$ for $(x_1,x_2)\in [0,1]^2$.
\end{proposition} 
\noindent
{\bf Proof:} First, we show that $\Pi_n(x_1,x_2)\leq 1$ for $(x_1,x_2)\in [0,1]^2$. Indeed, this follows from \eref{op} since for $(x_1,x_2)\in [0,1]^2$,
\begin{eqnarray}
\nonumber
\Pi_n(x_1,x_2)&=&2S_n\lr{\frac{x_1+x_2}{2}}- \frac{1}{2}\{ S_n(x_1)+S_n(x_2)\}  
\nonumber\\
&\leq &(x_1+x_2)-\frac{1}{2}
(x_1^2+x_2^2) 
=\frac{1}{2}\left[
x_1(2-x_1)+x_2(2-x_2)\right]
\leq 1.
\nonumber
\end{eqnarray}
To show that $\Pi_n\geq 0$, we start with
\be
\label{ll1}
2\Pi_n(x_1,x_2)=
 x_1+x_2+\sum_{k= 1}^n 4^{-k} 
[H^{\circ k}(x_1)+H^{\circ k}(x_2) -4 H^{\circ k}\left(\frac{x_1+x_2}{2}\right)].
\ee 
We introduce the function 
\be 
\nonumber
\zeta(t):= 2\min\{|t-m|\,:\,m\in \Z\},\quad t\in\R .
\ee 
Then, for $t\in [0,1]$, we have
\be 
\nonumber
H(t)= \zeta(t)\quad {\rm and} \quad H^{\circ k}(t)=\zeta(2^{k-1}t), \quad k\ge 2.
\ee 
Since $\zeta$ is subadditive, i.e., $\zeta(t+t')\leq \zeta(t)+\zeta(t')$,  we have
\be
\label{rs11}
H^{\circ k}\left(\frac{x_1+x_2}{2}\right)\leq 
H^{\circ k}\left(\frac{x_1}{2}\right)+H^{\circ k}\left(\frac{x_2}{2}\right) =H^{\circ (k-1)}(x_1)+ H^{\circ (k-1)}(x_2).
\ee 
We now replace each term  $H^{\circ k}\left(\frac{x_1+x_2}{2}\right)$ appearing  in   \eref{ll1} by the right side of \eref{rs11}.  The result is a telescoping sum.  Since $H(t/2)=t$, $t\in [0,1]$ this telescoping sum gives    
$$
2\Pi_n(x_1,x_2)\ge 4^{-n}[H^{\circ n}(x_1)+H^{\circ n}(x_2)]\ge 0,
$$
as desired \hfill $\Box$

 Next, we observe that $\Pi_n$ approximates $\Pi$ with exponential accuracy.
\begin{proposition}
   \label{P:mult} 
    For each $n\ge 1$, the function 
    $\Pi_n\in \Upsilon^{5,3n}(\Relu;2,1)$ and satisfies
    the inequalities
   \be 
   \nonumber
   \|\Pi-\Pi_n\|_{C([0,1]^2)}\le  
    4^{-n},
   \ee 
  and
   \be 
   \label{Pmult11}
   \|\partial_i\Pi-\partial_i\Pi_n\|_{L_\infty([0,1]^2)}\le 2\cdot 2^{-n},\quad \hbox{where}\quad \partial_i:=\partial_{x_i}, \quad i=1,2.
   \ee 
  \end{proposition}
   \vskip .1in
   \noindent
   {\bf Proof:} Let $\cN$ be the network of width 
   $W=4$ and depth $n$ 
   which   outputs $S_n$, see Proposition \ref{Recoversquare}. We now construct
   a network $\cN'$  which inputs $(x_1,x_2)$  and outputs $\Pi_n$. First, we add a 
   source channel to $\cN$ to push forward $x_2$
   ($\cN$ already has a source channel to push $x_1$). Then, we place 3 copies of this extended network next to each other.  We  output the three terms from \eref{product1} in the collation channel of $\cN$, and produce  
   $\Pi_n(x_1,x_2)$. 
   The new network has width $W=5$ and depth  
  $L=3n$. From \eref{squareerror}, we have $\|\Pi-\Pi_n\|_{C([0,1]^2)}\leq 4^{-n}$. 
   
  Finally, we check \eref{Pmult11} for $i=1$. The case $i=2$ is the same.  We have
   $\partial_{1}\Pi(x_1,x_2) =x_2$, and modulo a set of  measure zero,
\begin{eqnarray}
   \partial_1 
   \Pi_n(x_1,x_2)
   &=&\ 
   S_n'\lr{\frac{x_1+x_2}{2}} -  \frac{1}{2} S_n'(x_1)\nonumber \\
   &=& x_1+x_2+\e_1-\frac{1}{2}(2x_1+\e_2)
   =x_2+\e_1-\e_2 \nonumber\\
   &=&\partial_1\Pi(x_1,x_2)
   +\e_1-\e_2 \le 2^{-n+1},
   \nonumber
\end{eqnarray}
   where $|\e_1|,|\e_2|\le 2^{-n}$ because of \eref{Lipproduct}. The proof is completed.
   \hfill $\Box$}
 
  In general, we can 
    approximate any product
  \be 
 \nonumber
  \Pi^k(x_1,\dots,x_k):= x_1x_2\cdots x_k, \quad x_1,\dots,x_k\in[0,1],
   \ee 
up to  exponential accuracy,  using outputs of ReLU neural networks.   
We write
   $$\Pi^{k+1}(x_1,\dots,x_{k+1})= \Pi(x_{k+1},\Pi^k(x_1,\dots,x_{k})),
   $$
   denote $\Pi_n^2:=\Pi_n$, see \eref{PP}, and
   recursively define  
   $$\Pi_n^{k+1}(x_1,\dots,x_{k+1})= \Pi_n(x_{k+1},\Pi_n^k(x_1,\dots,x_{k})),\quad k=2,3,\dots.
   $$
It follows by induction, using Proposition \ref{Pin}, that $\Pi^{k}(x_1,\ldots,x_k)\in [0,1]$, and therefore $\Pi^{k+1}$ is well defined.
   Then, the following theorem holds. 
   \begin{theorem} 
   \label{T:product}
   For each $k\ge 2$, the  function $\Pi_n^k\in\Upsilon^{3+k, 3(k-1)n}(\Relu;k,1)$ and satisfies
   \be 
   \label{product21}
   |\Pi^k(x_1,\dots,x_k)-\Pi^k_n(x_1,\dots,x_k)|\le  C_k \cdot 4^{-n},\quad  x_1,\dots,x_k\in [0,1],
   \ee 
  where  for $k\ge 2$,
  $$C_k\le (k-1)\alpha_n^{k-2},\quad
  \alpha_n:= 1+2^{-n+1}.
  $$
  In particular, $C_k\le ek$, as long as $n\geq 1+\log_2 k$.
    \end{theorem} 
   \noindent
   {\bf Proof:}
 For $k\geq 3$, we construct a   network of width $(k+2)$  which takes the inputs $x_1,\dots,x_{k-1}$ and outputs $\Pi^{k-1}_n$ (when $k=3$, this is the network for
 $\Pi_n^{2}$).  Its first $3n(k-2)$ layers are the same as the network that inputs 
   $x_1,\dots,x_{k-1}$ and  outputs $\Pi_n^{k-1}$, except that we add an additional channel to push forward $x_{k}$.  We then follow this with the network for $\Pi_n^2$ using as inputs $x_{k}$ and $\Pi^{k-1}_n(x_1,\dots,x_{k-1})$.   This network will have width $W=k+3$
   and depth $L=3(k-1)n$ as desired.
   
    Next, we fix $n$ and prove  \eref{product21}    by induction on $k$.  The case $k=2$ is covered by Proposition \ref{P:mult} with 
    $C_2=1$.  To advance the induction hypothesis, we  assume that we have proven the result for some 
    $(k-1)$ with constant $C_{k-1}$,
   We write  (with the obvious abbreviation of notation)
  $$
   \Pi^{k}-\Pi^{k}_n= \Pi(x_{k},\Pi^{k-1})-\Pi_n(x_{k},\Pi^{k-1})+\Pi_n(x_{k},\Pi^{k-1})-\Pi_n(x_{k},\Pi_n^{k-1}),
   $$
and use Proposition \ref{P:mult}   to obtain
      \be 
      \label{Pmult22}
       \|\Pi^{k}-\Pi^{k}_n\|_{C([0,1]^{k})}\le 
       4^{-n}+ \|\partial_2\Pi_n\|_{L_\infty(([0,1]^2)} \|\Pi^{k-1}-\Pi_n^{k-1}\|_{C([0,1]^{k-1})}.
      \ee
      We now use \eref{Pmult11}, to conclude that $\|\partial_2\Pi_n\|_{L_\infty(([0,1]^2)}\le 1+2\cdot 2^{-n}$. Inserting this into
      \eref{Pmult22} gives
      $$
       \|\Pi^{k}-\Pi^{k}_n\|_{C([0,1]^{k})}\le 
       4^{-n}+ (1+2^{-n+1})C_{k-1}4^{-n} =(1+\alpha_n C_{k-1})4^{-n}.
       $$
    The recurrence formula  $C_{k}=1+\alpha_nC_{k-1}$, $k\ge 3$, with initial value $C_2=1$,  has the solution 
    $$
    C_{k} = \sum_{j=0}^{k-2}\alpha_n^j\le (k-1)\alpha_n^{k-2}. 
    $$ 
    Moreover, if $k\leq 2^{n-1}$, we have
    $$
    C_k\leq (k-1)\left(1+\frac{1}{2^{n-1}}\right)^{k-2}< k\left(1+\frac{1}{k}\right)^k<ek.
    $$
    This completes the proof of the theorem. \hfill $\Box$
    
       \begin{remark}
    \label{remR}
    If $a>1$, then a simple change of variables gives that  the function $\Pi^k$,  $k\geq 2$, now considered as a function in $C([0,a]^k)$, is approximated by $\tilde \Pi_n^k(x_1,\dots,x_k):=a^k \Pi_n^k(x_1/a,\dots,x_k/a)$  with accuracy
     $$
    \|\Pi^k-\tilde \Pi_n^k\|_{C([0,a]^k)}\leq C_ka^k\cdot 4^{-n},  
    $$
    with $C_k$ as in Theorem \ref{T:product}.  Moreover $\tilde \Pi_n^k\in \Upsilon^{3+k,3(k-1)n}(\Relu;k,1)$.
    
    \end{remark}

\vspace{.5cm}

   \subsection{Approximation of polynomials}
   \label{SS:polynomials}
   In the following, we show that polynomials can be well approximated by the outputs of  deep ReLU networks.  We begin with monomials.
   \vskip .1in
   \noindent
   {\bf Approximation of monomials:}    For  $\nu\in \N^d$, let $\phi_\nu(x):=x^\nu$, 
   $x\in [0,1]^d$.  We use the standard notation $|\nu|=\sum_{j=1}^d\nu_j$ for the length of $\nu$.
   Theorem \ref{T:product} shows that any  monomial $\phi_\nu$, with $|\nu|=m$,  is well approximated by deep ReLU networks. Namely, for each $n\ge 1$, there is an $S_\nu
   \in \Upsilon^{3+d,3(m-1)n}(\Relu;d,1)$ such that
   \be 
   \label{monomial1}
   \|\phi_\nu-S_\nu\|_{C([0,1]^d)} \le  em \cdot 4^{-n},\quad 
   n\geq 1+\log_2m.
   \ee
   Note that here we can keep the width of the network bounded by $3+d$ rather than $3+m$ because the $x_1,\dots,x_d$ are repeated; we leave the details to the reader.
   
   \vskip .1in
   \noindent
   {\bf Approximation of polynomials:}   If $P(x)=\sum_{\nu\in\Lambda}c_\nu x^\nu$, $x\in [0,1]^d$, where all of the indices $\nu\in\Lambda$ satisfy $|\nu|\le m$, then we can approximate $P$ by the function $S:= \sum_{\nu\in\Lambda}c_\nu S_\nu$.  Note that $S$ is the output of the concatenation of the networks
   that output the $S_\nu$'s, see {\bf Addition by increasing depth}  in \S\ref{SS:highlighted}.
   Since all networks producing the $S_\nu$'s have already source channels that forward the values $x_1,\ldots,x_d$, we need to add only a collation channel to collect the terms in the sum defining $S$, and therefore $S\in \Upsilon^{4+d,3(m-1)n\#(\Lambda)}(\Relu;d,1)$.
   We have the following estimate for the approximation error
   \be 
   \label{errorpoly}
   \|P-S\|_{C([0,1]^d)}\le 
   em \cdot 4^{-n}  \sum_{\nu\in\Lambda}|c_\nu|,
   \quad n\geq 1+\log_2m,
   \ee 
   obtained from \eref{monomial1}. There are several savings that can be made in the size of the network in such constructions by balancing the size of $c_\nu$ with the size of the networks for the $S_\nu$ when given a desired target accuracy.
   
   Constructions of NN approximations to polynomial sums have been employed to prove results on approximating real analytic functions using deep
   ReLU networks.  We do not formulate those results here but rather refer the reader to the papers \cite{OSZ} and \cite{E1}
   for statements and proofs.
   
   \subsection{Approximation of tensor products}
   \label{SS:tensors}
   Tensor structures are a very effective method for   approximation   in high dimensions.  It is beyond the scope of this article to lay this subject out in its full detail. We simply wish to point out here that
   a rank one tensor product
   \be 
   \label{rankone}
   f(x_1,\dots,x_d)=f_1(x_1)\cdots f_d(x_d),\quad  x_1,\dots,x_d\in [0,1],
   \ee 
   is well approximated in $C(\Omega)$, $\Omega=[0,1]^d$,  whenever the univariate components $f_j$ are well approximated. The starting point for this is the following simple proposition.
   \begin{proposition}
   \label{P:tensor}
   If $g_j:[0,1]\to [0,1]$, $g_j\in\Upsilon^{W_0,L_0}(\Relu;1,1)$, $W_0\ge 3$, 
  for $j=1,\dots,d$, then the rank one tensor
   \be 
   \label{r11}
   g(x_1,\dots,x_d)=g_1(x_1)\cdots g_d(x_d),
   \ee 
   can be approximated by $S\in \Upsilon^{dW_0,L_0+3(d-1)n}(\Relu;d,1)$ to accuracy
   \be 
   \label{r11acc}
   \|g-S\|_{C(\Omega)}\le 
   ed\cdot 4^{-n}, \quad n\ge 1+\log_2 d.
   \ee 
   \end{proposition} 
   \noindent
   {\bf Proof:}  We can take $
   S:=\Pi_n^d(g_1,\dots,g_d)$.   We claim that $S$ is an element of  $\Upsilon^{dW_0,L_0+3(d-1)n}(\Relu;d,1)$.  Indeed,
      we stack  the networks producing the $g_j$'s, $j=1,\dots,d$   on the top of each other and end up with a network $\cN_1$ with width $dW_0$ and depth $L_0$. Then, we concatenate it with the network $\cN_2$  producing $\Pi_n(y_1,\ldots,y_d)$. The latter has   depth $3(d-1)n $ and width $3+d\le dW_0$. The concatenation is done by forwarding the output of the network producing $g_j$ as an input to the 
$j$-th channel of $\cN_2$
(recall that the first $d$ channels of $\cN_2$ are source channels). We end up with a network with the desired width and depth. 
   The inequality \eref{r11acc} follows from Theorem \ref{T:product}.
   \hfill $\Box$
   
   \begin{remark}
   \label{R:tensor}
   We make two remarks on the above proposition:
   \begin{itemize}
       \item
   If  $0\leq g_j(t)\leq M$  instead of $0\leq g_j(t)\leq 1$, then by using Remark 
   \ref{remR}, we obtain an $S$ in the same ReLU space but the accuracy of approximation is now lessened by the factor $M^{d}$.
   \item If the $g_j$'s are not in the designated ReLU space, but are rather only approximated by $\hat g_j:[0,1]\to[0,1]$, $j=1,\dots,d$, from the designated ReLU space to an accuracy $\e$, then the function $S:=\Pi_n^d(\hat g_1,\dots,\hat g_d)$ is in the designated ReLU space and we can write
   \begin{eqnarray}
   \label{tenserror}
   \|g-S\|_{C(\Omega)}&\le& \|g-\Pi(\hat g_1\cdots\hat g_d)\|_{C(\Omega)} + \|[\Pi -\Pi_n](\hat g_1\cdots\hat g_d)\|_{C(\Omega)}\nonumber\\
   &\le& d\e+ ed\cdot 4^{-n}, \quad n\geq 1+\log_2 d.
  \end{eqnarray}
   where the first term does not exceed the sum of the $d$ errors
   \be
   \label{derrors}
  \| \hat g_1\cdots \hat g_j\cdot g_{j+1}\cdots g_d - \hat g_1\cdots \hat g_{j+1}\cdot g_{j+2}\cdots g_d\|_{C(\Omega)}\le \e. 
  \ee
   \end{itemize}
   \end{remark}

   \subsection{Approximation of B-splines}
   \label{SS:Bsplines}
   In our   presentation of classical smoothness classes $K$  of functions given in \S\ref{S:modelclasses}, we have stressed   that the elements in $K$ have
   certain atomic decompositions and their membership in  $K$ is characterized by the decay of their coefficients in such representations.  Thus, if we can show that the atoms in such a decomposition are well approximated by NNs, then we can obtain bounds for
   NN approximation of $K$.  The aim of the present section is to show how this unfolds when  we choose B-splines as the atomic representation system.
   
    Let $N_r$ be the univariate B-spline defined in \eref{Bspline}.  Let us recall that $N_r$ is supported on $[0,r]$ and is normalized so that $\|N_r\|_{C(\R)}=1$.
   
   \begin{proposition}
   \label{P:approxN}
   Let $r\ge  2$ and consider the  B-spline 
   $$N(x):=N_r(x_1)\cdots N_r(x_d)$$
   of $d$ variables.  There is a function $\hat N\in \Upsilon^{W,L}(\Relu;d,1)$ with  width $W=6d$ and depth $L=Cn$, with $C$ depending only on $r$ and $d$,  which satisfies
   \be 
   \label{approxN1}
   \|N-\hat N\|_{C(\R^d)}\le C'(r,d) 4^{-n}, \quad n\ge 1,
   \ee 
   with the constant $C'(r,d)$ depending only on $r$ and $d$.  Moreover, the support of $\hat N$ is contained in that of $N$.
   \end{proposition}
   \noindent
   {\bf Proof:}   This is proved by approximating in succession the functions
   \begin{equation}
   \label{thefunctions}
    t^{r-1},\ \rho_{r-1}(t):=t^{r-1}_+,\ N_r(t),\ N(x_1,\dots,x_d)=N_r(x_1)\cdots N_r(x_d),
    \end{equation}
   where $N_r$ is the univariate B-spline.   Our results of the previous sections on approximating products were stated for approximation
   on $[0,1]^d$ and now we want approximation on $[0,r]^d$.  
   This is done by using Remark \ref{remR}
   and changes the estimates by a constant depending
   only on $r$ and $d$.  
   We assume such changes without further elaboration in what follows.  All constants $C$ appearing in the proof depend at most on $r$ and $d$.

    Because of \eref{monomial1}, we can approximate
   the function $t^{r-1}$ by an element of $\Upsilon^{4,3(r-2)n}(\Relu;1,1)$ with an  error that does not exceed $C4^{-n}$.  By adding an extra layer for $\Relu$,
   we can approximate the function 
   $\rho_{r-1}$  by an $S_r$ from $\Upsilon^{4,3(r-2)n+1}(\Relu;1,1)$ with accuracy 
   \be 
   \label{approxB}
   \|\rho_{r-1}-S_r\|_{C([0,r])}\le C4^{-n}.
   \ee 
   
   Next, we use $S_r$ to approximate the univariate B-spline 
   $N_r$ by replacing 
   $\rho_{r-1}(k-t)$ by $S_r(k-t)$ in formula \eref{Bspline}. The resulting function
   $T_r$, see Theorem \ref{T:translate-dilate}, is in   $\Upsilon^{6,3(r+1)(r-2)n+r+1}(\Relu;1,1)$ 
    (note that the network producing the latter set has 1 source channel for $t$ and one collation channel). 
   In addition, we have
   \be 
   \label{approxB1}
   \|N_r-T_r\|_{C([0,r])}\le 
 C 4^{-n},\quad \hbox{and}\quad 
 \|T_r\|_{C([0,r])}\le 1+C4^{-n}.
   \ee 
 We next consider $T^+_r:=\Relu(T_r)$ which also satisfies \eref{approxB1}, with the additional property $T^+_r\geq 0$.
 
Remark \ref{R:tensor} with $\e=C4^{-n}$ and $M\leq (1+C/4)$ gives that 
$N$ can be approximated by
$\tilde  N=\Pi^d_n(T^+_r,\ldots,T^+_r)$ with accuracy
\be 
   \label{lastB}
   \|N-\tilde  N\|_{C([0,r]^d)} \le C 4^{-n}, \quad n\geq 1+\log_2 d.
   \ee 
The approximant $\tilde N\in \Upsilon^{W,L}(\Relu;d,1)$ with  width
$W=6d$ and depth $L=3(r+1)(r-2)n+r+2+3(d-1)n$.
The network producing $\tilde N$ has 
$d$ source channels for each of the variables $x_i$, $i=1,\ldots,d$, and $d$ collation channels.  

Finally, we modify the function $\tilde N$ of \eref{lastB} so that it vanishes outside $[0,r]^d$.  This is done by what should by now be
    a familiar technique to the reader.  We construct a function $S$ with support $[0,r]^d$  and $S\geq N_r$, 
    using the method for the construction of nodal functions, see \eref{expressphi}.  
    More precisely, $S:=\Relu(\min\{\ell_1, \ldots,\ell_{2^d}\})\in \Upsilon^{d+1,2^d}(\Relu;d,1)$, see {\bf MM2} of \S\ref{SS:highlighted},
    where 
     $\ell_j$, $j=1,\ldots, 2^d$, are affine functions, each of which vanishes on one  of the $2^d$ facets of the cube $[0,r]^d$ and is above the graph of the B-spline $N_r$.    
    Then, the function  $\hat N: =\Relu(\min\{\tilde N,S\})$ agrees with $\tilde N$ when $\tilde N$ is non-negative and vanishes outside $[0,r]^d$. Therefore  it satisfies all the properties of the theorem. 
    Since both networks producing $\tilde N$ and $S$ have source and collation channels, it follows that 
    $\hat N\in \Upsilon^{6d,Cn}(\Relu;d,1)$ with $C$ 
    depending only on  $r$ and $d$.
    \hfill $\Box$

  \subsection{Approximation of Besov classes with deep ReLU networks}
  \label{SS:BesovBsplines}
  With the results of the previous section on B-spline approximation in hand,  we can now show that deep neural networks are at least as
  effective as standard nonlinear methods (modulo logarthms), such as adaptive FEMs or $n$-term wavelets, when approximating the classical smoothness spaces (Sobolev and Besov).  The ideas used in the presentation below are put forward in the references \cite{Nouy},
  \cite{BGKP},
  \cite{gribonval2019approximation}.

  We fix $\Omega=[0,1]^d$ and measure the approximation   error in some $L_p(\Omega)$ norm, $1\le p\le\infty$.
  Let $K=U(B_q^s(L_\tau(\Omega)))$ be the unit ball of a Besov space that lies above the Sobolev embedding line for $L_p(\Omega)$,
  and therefore is  a compact subset of $L_p(\Omega)$.
  Classical methods of nonlinear approximation
  show that $K$ can be approximated in $L_p(\Omega)$ to accuracy ${\cal O}(n^{-s/d})$, where $n$ is the number of parameters used in the approximation.
  We now show how to achieve these rates  when using $\Sigma_n:=\Upsilon^{W_0,n}(\Relu;d,1)$ with $W_0$ fixed and depending only on $s$ and $d$.
  Note that the number of parameters needed to describe the elements in $\Sigma_n$ is at most $C(s,d)n$. Here, and later in this section, the constant $C(s,d)$ changes at each occurrence.

  \begin{theorem}
  \label{T:Besov}
  Let $s>0$ and $\Omega=[0,1]^d$.  Suppose $K=U(B_q^s(L_\tau(\Omega)))$, $0<q,\tau\le\infty$, is the unit ball of a Besov space lying above the Sobolev embedding line
  for $L_p(\Omega)$ with 
  $1\le p \le \infty$, that is
  \be 
  \nonumber
  \delta :=s -\frac{d}{\tau} +\frac{d}{p}>0.
  \ee 
  Then, we have
  \be 
  \label{errorBesov}
  E(K,\Sigma_{n [\log_2 n]^{\beta}})_{L_p(\Omega)}\le C(s,d,\delta)n^{-s/d},\quad n\ge 1,
  \ee 
  where 
  $\Sigma_n:= 
  \Upsilon^{6d,Cn}(\Relu;d,1)$, $n\ge 1$, with $C=C(s,d,\delta)$ fixed, depending only on $s$, $d$ and $\delta$, and where
   $\beta:=\max\{1,\lceil 2d/(s-\delta)\rceil\}$. 
    \end{theorem}
  \noindent
  {\bf Proof:}  We only treat the case $1\le p<\infty$ and leave to the reader to make the necessary changes for $p=\infty$.
 We fix $p$ and $s>0$. We can assume  $q=\infty$ since this is the largest unit ball for the given $\tau$, and $\tau<p$.  We can further assume that $\delta>0$ is arbitrarily small since the Besov spaces of order $s$ get larger as we approach the Sobolev embedding line which corresponds to $\delta=0$.

 To prove the theorem, it is sufficient to prove that it holds for $n=2^L$ with $L$ a sufficiently large positive integer.
 We take $r=\lceil s\rceil +1$ and let $N$ denote the multivariate tensor product B-spline of order $r$.  We recall the notation $\cD$ for dyadic cubes, $\cD(\Omega)$ for dyadic cubes $I\in\cD$ such that $N_I$ is nonzero on $\Omega$,  $\cD_k(\Omega)$ for these cubes at dyadic level $k$ (they have measure $2^{-kd}$), and 
 $\cD_+(\Omega):=\cup_{k\geq 0}\cD_k(\Omega)$.

 From \eref{Brep} and \eref{Bnorm}, we know that any  $f\in K$ has the representation
 \be 
 \label{repK}
 f=\sum_{I\in\cD_+(\Omega)} c_I(f)N_I,
 \ee 
 with
 \be 
 \label{dyadicbound1}
  \sum_{I\in\cD_k(\Omega)} |c_I(f)|^\tau |I| \le 2^{-ks\tau}, \quad k=0,1,2,\dots 
 \ee 
Here, we use   \eref{Bnorm} for the definition of the norm in the  Besov space.
For each $j\in \Z$ and $k\ge 0$, we define
 \be 
 \label{defLambdajk}
 \Lambda(j,k):=\{I\in\cD_k(\Omega): 2^{-j}\le |c_I(f)|<2^{-j+1}\},
 \ee 
 and estimate its  cardinality   from \eref{dyadicbound1}. We derive  that
 \be 
  \label{critical}
 \sum_{j=-\infty}^\infty
  2^{-j\tau} \#(\Lambda(j,k))\le 2^{k(d-s\tau)} ,\quad k=0,1,\dots . 
  \ee 
   It follows from \eref{critical} that
if $\Lambda(j,k)\neq \emptyset$, then
 $2^{-j\tau}\leq 2^{k(d-s\tau )}$, and therefore
\be
\label{jk}
 j\geq \left(s-\frac{d}{\tau}\right)k=
 \left(\delta-\frac{d}{p}\right)k=:J_k. 
\ee 
 In other words,
 \be 
 \label{jk1}
 \Lambda(j,k)=\emptyset, \quad\hbox{when}\quad  j< J_k.
 \ee 
 We will now replace 
 some of the $N_I$'s from \eref{repK}
 by approximants $\hat N_I$ from $\Sigma_{m(I)}$,
 where the nonnegative integers  $m(I):=m(j,k)\in\{1,2,\ldots\}$ will be chosen the same for each $I\in \Lambda(j,k)$ (as we shall see  below). The $N_I$'s that are not approximated are associated with    $m(I)=0$.
 
 It follows from  \eref{approxN1}, that
 \be 
 \label{guaranteeN1}
 \|N_I-\hat N_{I} \|_{L_p(\Omega)}\le C|I|^{1/p}4^{-m(I)},
 \ee 
 where
    here and later in this proof all constants $C$ depend only on $s,d$ and $\delta$.
  According to Proposition \ref{P:approxN}, we can also assume that  $\hat N_I$ is zero outside the support of $N_I$.
   
Next, we define
the functions
 \be 
 \label{weak}
 \hat S:= \sum_{I\in \cD_+(\Omega), m(I)>0} c_I(f)\hat N_I, \quad  \hat S_k:= \sum_{I\in \cD_k(\Omega), m(I)>0} c_I(f)\hat N_I, \quad k\geq 0,
 \ee 
 and proceed to  show that $\hat S$ provides the needed approximation if we choose $m(I)$ appropriately.
 
 In preparation for the choice of the $m(I)$, we first estimate how well $\hat S$ approximates $f$.  
 If we denote by 
  $S_k:= \sum_{I\in\cD_k(\Omega)} c_I(f)N_I$, $k\geq 0$, using
 \eref{guaranteeN1} 
  and  the fact that
 $\|N_I\|_{L_p(\Omega)}\le C|I|^{1/p}$, we obtain,
 
 \begin{eqnarray} 
 \nonumber
 \|S_k-\hat S_k\|^p_{L_p(\Omega)} &\le& C^p\sum_{I\in\cD_k(\Omega)} |c_I(f)|^p|I|4^{-m(I)p}\nonumber \\
 &\le&  C^p2^{-kd}\sum_{j\geq J_k} 2^{-jp} \#(\Lambda(j,k))4^{-m(j,k)p} \nonumber \\
 \nonumber
 &=&   C^p 2^{-kd}\sum_{j\geq J_k}2^{-j\tau}\#(\Lambda(j,k))
 2^{-2m(j,k)p-jp+j\tau}.      
 \end{eqnarray}
 For the definition of $m(j,k)$, let us introduce the notation
 \be
 \label{epsilon}
 \e_t:=2\log_2(t+1), \quad t\geq 0.
 \ee
 For every $j\geq J_k$,  $k\geq 0$, 
 we choose 
 $m(j,k)$ to be the smallest non-negative integer such that 
 $$
 \e_kp+ Lsp/d\leq ks\tau+2m(j,k)+1]p+j(p-\tau). 
 $$
  This choice satisfies
\be 
 \label{Lk1}
 0\leq 2m(j,k)p\le [\e_kp +Lsp/d -jp  - (ks-j)\tau]_+.
 \ee 
 Then, we obtain
\begin{eqnarray}
\nonumber
 \|S_k-\hat S_k\|^p_{L_p(\Omega)}
 &\leq& C^p 2^{-kd}2^{-\e_kp-Lsp/d+ks\tau}\sum_{j\geq J_k}2^{-j\tau}\#(\Lambda(j,k))\\ \nonumber
 &\leq& C^p2^{-Lsp/d}(k+1)^{-2p},
\end{eqnarray}
 where we used \eref{critical} for the last inequality.
It then follows from \eref{repK} that
 $$
 \|f-\hat S\|_{L_p(\Omega)}\leq
 \sum_{k=0}^\infty \|S_k-\hat S_k\|_{L_p(\Omega)}\leq
 C2^{-Ls/d}\sum_{k=0}^\infty
 (k+1)^{-2}=Cn^{-s/d}.
 $$

We are left to show that $\hat S\in \Sigma_{L^{\beta}2^L}$, 
 which in turn proves the theorem. We know that $\hat N_I\in 
\Upsilon^{6d,Cm(I)}(\Relu;d,1)$.   Since the 
network producing $\hat N_I$ already has $d$ source channels and a collation channel,  our  {\bf Addition by increasing depth} of \S\ref{SS:highlighted} gives that 
  $\hat S\in\Upsilon^{6d,CA}(\Relu;d,1)$, where
  \be 
  \label{boundA}
  A:= \sum_{k=0}^\infty \sum_{j=J_k}^{J_k^+} m(j,k)\#(\Lambda(j,k)).
  \ee 
The index in the second sum in \eref{boundA} has upper bound
$J_k^+$, where $J_k^+$ is defined by the equation
  \be 
  \label{nzpairs}
     J_k^+(1-\tau/p)+  ks\tau/p = \e_k +Ls/d,
  \ee
because  $m(j,k)=0$ when $j\ge J_k^+$, see \eref{Lk1}.
Later, we shall use the fact that
  \be 
  \label{kplus}
    \tau J_k^+= \lambda (\e_k +Ls/d -  ks\tau/p ), \quad {\rm with}\quad \lambda:=(1/\tau-1/p)^{-1}= \frac{d}{s-\delta}.
  \ee  
  
 For $j\in [J_k,J_k^+]$, $m(j,k)$ takes its maximum value at $j=J_k$ which is
  \begin{eqnarray}
 \nonumber
  m(J_k,k)&\le& \frac{1}{2}(\e_k+Ls/d-ks\tau/p -J_k(1-\tau/p))\nonumber \\
 &=& \frac{1}{2}(2\log_2(k+1)-k\delta+Ls/d) \nonumber \\
 &\le & CL,
 \nonumber
 \end{eqnarray} 
 where we used the definition of $J_k$ and \eref{Lk1}.
 Therefore, we have the estimate
  \be 
  \label{bounda}
  A\le CL \sum_{k=0} ^\infty \sum_{j=J_k}^{J_k^+} \#(\Lambda(j,k))\le CL\sum_{k=0}^{L/d-1} 2^{kd}+ 
  CL\sum_{k=L/d}^\infty 2^{kd-ks\tau +J_k^+ \tau },
  \ee 
   where in the first sum we used the fact that
   $$
\sum_{j=J_k}^{J_k^+} \#(\Lambda(j,k))\leq
 C2^{kd},
$$
because $\Lambda(j,k)\subset \cD_k(\Omega)$,
  and the second sum used that  
\begin{eqnarray}
\nonumber
  \sum_{j=J_k}^{J_k^+} \#(\Lambda(j,k))&\leq&
  2^{J_k^+\tau}\sum_{j=J_k}^{J_k^+} 2^{-j\tau}\#(\Lambda(j,k))
  \leq
 2^{J_k^+\tau}\sum_{j=J_k}^\infty 2^{-j\tau}\#(\Lambda(j,k)) \nonumber \\
 &\leq& 2^{kd-ks\tau+J_k^+\tau}.
 \nonumber
\end{eqnarray}

Obviously, the first sum on the right does not exceed $CL2^L$, and so we concentrate on the second sum.
     We first want to see what the exponent is in that sum.  From \eref{kplus}, we have
  \be 
  \label{findJ}
  kd  -ks\tau+J_k^+\tau  = \lambda(Ls/d+\e_k) +k\{d -s\tau(1+\lambda/p)\}.
  \ee   
  Going further, we find
\begin{eqnarray}
  \label{further1}
  d -s\tau(1+\lambda/p) &=&
  d-s\tau\left(1+\frac{1}{p(\frac{1}{\tau}-\frac{1}{p})}\right)=
  d-s\tau\left(1+\frac{1}{\frac{p}{\tau}-1}\right) \nonumber\\
  &=&
  d-\frac{sp}{\frac{p}{\tau}-1}=d-\frac{s}{\frac{1}{\tau}-\frac{1}{p}}=
  \frac{d(\frac{1}{\tau}-\frac{1}{p})-s}
  {\frac{1}{\tau}-\frac{1}{p}} \nonumber \\
  &=&
  \frac{-\delta}{\frac{1}{\tau}-\frac{1}{p}}=-\delta\lambda,
  \nonumber
\end{eqnarray}
 and thus
 $$
 kd  -ks\tau+J_k^+\tau  = \lambda(Ls/d+\e_k- \delta k).
 $$
 We substitute the latter relation into  \eref{bounda} and obtain, after change of index $i=k-L/d$ and using \eref{kplus},
\begin{eqnarray}
\nonumber
A&\le&  
 CL2^L+CL\sum_{k=L/d}^\infty
 2^{\lambda(Ls/d+\e_k- \delta k)}\nonumber \\
 &=&
 CL2^L+\sum_{i=0}^\infty
 2^{\lambda(s-\delta)L/d}
 2^{2\lambda\log_2(i+L/d+1)-i\lambda\delta}\nonumber \\
 &=&
 CL2^L+2^L\sum_{i=0}^\infty
 (i+L/d+1)^{2\lambda}2^{-i\lambda\delta}<
 CL2^L+CL^{2\lambda}2^L<
 CL^\beta2^L. 
 \nonumber
\end{eqnarray}
  This gives the bound  we want and proves the theorem. \hfill $\Box$

  Before proceeding on, we make the following remarks concerning the above theorem and its proof.
  
  \begin{remark}
  \label{R:proofT}
  The above result is not quite as good as the results for approximating unit balls of Besov classes when using
  other methods of nonlinear approximation, see \cite{devore1998nonlinear}, because of the appearance of the logarithm.
  We should mention that in \cite{Nouy} the authors prove a similar  to the above theorem 
  result by using spline wavelets
  rather than B-splines as the main vehicle.
  In the next section, we show that when $p=\infty$ this logarithm does not appear and, in fact, we can prove much
  better rates of approximation. These can in turn be used to improve the above results for all Besov classes, as we shall discuss at the end
  of the next section.  The determination of the best approximation rates for $L_p$ approximation  when using the outputs of deep networks such
  as $\Sigma_n=\Upsilon^{W_0,n}(\Relu;d,1)$  remains unsettled. 
  \end{remark}

  \subsection{Super convergence for deep ReLU networks}
  \label{SS:super}
  In this section, we present some very intriguing results
  on approximation by NNs that show quite unexpected rates of approximation for certain classical model classes $K$ described by smoothness.
  The initial results were given in \cite{yarotsky2018optimal} for the model classes $U({\rm Lip} \ \alpha)$, $0<\alpha\le 1$ on $\Omega=[0,1]^d$, and  were  later extended to more general model classes $U(C^s(\Omega))$, $s>0$, in \cite{lu2020deep}.   Our aim in this section is to show how  these super rates are established in the simple case of univariate functions in Lip 1, and leave the reader to consult the above references for the treatment for functions of $d$ variables and
  higher smoothness.  At the end of this section, we place these results into perspective and formulate some related questions.

  \begin{theorem} 
  \label{T:Yarotsky}
  If  $K= U({\rm Lip}\ 1)$ on $\Omega:=[0,1]$ and  $\Sigma:=(\Sigma_n)_{n\ge 1}$ with 
  $\Sigma_n:=
  \Upsilon^{11,16n+2}(\Relu;1,1)$, $n\ge 1$,  then  we have
  \be 
  \label{Yarotsky}
  E(K,\Sigma_n)_{C(\Omega)}\le 6n^{-2},\quad n\ge 1.
  \ee 
  \end{theorem}
  \noindent
  {\bf Proof:}  We use the same notation as in \S\ref{SSS:bitextraction}.   Namely, we define $N:=n^2$,with $n\geq 4$ an even positive integer and set $t_i:=i/N$, $0\le i\le N$,
  and $\xi_j:= j/n$, $j=0,\dots,n$.  Given $f\in K$, as a first step, we take  $S_0$ to  be the CPwL function with breakpoints precisely  the $\xi_j$'s, $j=1,\dots,n-1$, which interpolates $f$ at the $\xi_j$, $j=0,\dots,n$.  Then, $S_0$  has the following three properties:
  \begin{itemize}
  \item [(i)] $\|f-S_0\|_{C(\Omega)}\le 1/n,\quad n\ge 4$.
  \item [(ii)] $\|S_0\|_{{\rm Lip} \ 1}\le 1,\quad n\ge 4$.
  \item [(iii)] $S_0\in
  \Upsilon^{3,n}(\Relu;1,1)$, see 
  Proposition \ref{int11}.
  \end{itemize}
  Now, consider the  function $R:=f-S_0$.  It vanishes at each of the $\xi_j$, $j=0,\dots,n$, and $R\in {\rm Lip}\  1 $ with $\|R\|_{{\rm Lip}\ 1}\le 2$.  We next show that there is a sequence of $\e_i\in\{-1,+1\}$, $i=0,\dots,N-1$,  such that  $(y_i)_{i=0}^N$, defined recursively by $y_0:=0$ and
  \be 
  \label{relatedsequence}
  y_{i+1}:= y_i+\e_i, \quad i=0,\dots,N-1,
  \ee 
  satisfy
  \be 
  \label{satisfies1}
  y_{jn}=0,\quad j=0,\dots,n \quad {\rm and} \quad |R(t_i)-\frac{2y_i}{N}| \le \frac{2}{N},\quad 0\le i\le N.
  \ee 
  
  Let us for the moment assume we have found such a sequence $(\e_i)_{i=0}^{N-1}$ and show how to complete the proof of the theorem.  We apply 
  Theorem \ref{T:  bitextraction} to the $(y_i)_{i=0}^N$ defined in \eref{relatedsequence} and obtain  a function $S_1\in \Upsilon^{11,15n+2}(\Relu;1,1)$ with the properties guaranteed by
  this theorem.   Next, we consider the function 
  $$
  S:=S_0+\frac{2}{N}S_1.
  $$
  Note that $S\in \Upsilon^{11,16n+2}(\Relu;1,1)$ because of {\bf Addition by increasing depth}
  of \S\ref{SS:highlighted}, 
  taking into account that the network from Theorem \ref{T:  bitextraction} producing $S_1$ already has a source and collation channel. Moreover,
  \be 
  \label{satisfies2}
  \|f-S\|_{C([0,1])}=  \|R-\frac{2}{N}S_1\|_{C([0,1])}\le 6/N.
  \ee 
  Indeed, if  $t\in[t_i,t_{i+1}]$, $i=0,\ldots,N-1$, then   $S_1(t_i)= y_i$, and  we have 
  
  $$
  |R(t)-\frac{2}{N}S_1(t)|\le |R(t)-R(t_i)|+|R(t_i)- \frac{2}{N}S_1(t_i)|+\frac{2}{N}|S_1(t)-S_1(t_i)|\le 6/N,
  $$
  because of the Lipschitz properties of $R$, the properties of $S_1$, and \eref{satisfies1}.  This in turn would prove the theorem.

  So we are left with finding a sequence $(\e_i)_{i=0}^{N-1}$ such that  \eref{satisfies1} is valid.
  It is enough to show how to  define this sequence for $i=0,\dots,n-1$ since   for $i=jn,\dots, (j+1)n-1$ it is defined similarly.
 We choose  the sequence $\e_0,\e_1,\dots$ and the corresponding $y_{j+1}:=y_j+\e_j$  and verify \eref{satisfies1} recursively.  We first choose $\e_0\in\{-1,1\}$ so that $2\e_0/N$ is closest to
  $R(t_1)$ for this choice of the two possible values $\pm 1$. Clearly, since $|R(t_1)|\leq 2/N$, for $y_1:=\e_0$ we have  the inequality $|R(t_1)-\frac{2y_1}{N}|\le 2/N$. In other words,
  we have verified \eref{satisfies1} for $j=1$.

  Assume now that
  $\e_0,\dots,\e_{j-1}$ have been chosen and the corresponding $y_1,\dots,y_j$ have been shown to satisfy \eref{satisfies1}. We now choose $\e_j$ so that $\frac{2y_{j+1}}{N}=\frac{2(y_j+\e_j)}{N}$ is closest to $R(t_{j+1})$. Since $R$ changes by at most $2/N$ in moving from $t_j$ to $t_{j+1}$, this choice will
  also satisfy \eref{satisfies1}.  So, we are left to verify that $y_n=0$. Since $n$ is even, $y_n=\e_0+\ldots+\e_{n-1}=2m$ for some integer $m$. In addition, we have
  $|2y_n/N-0|\le 2/N$, and therefore we must have $m=0$.   Thus, we showed the existence of a sequence $(\e_i)_{i=0}^{N-1}$ with the required properties. The   proof of the theorem is completed.
  \hfill $\Box$

  \subsubsection{Remarks on Theorem \ref{T:Yarotsky}}
  \label{RMSY}
  We make some remarks on this theorem in order to put into perspective what it is saying.  In this section, 
  we take  $\Sigma:=(\Sigma_n)$, where the sets $\Sigma_n$ used for approximation  are $\Sigma_n=\Upsilon^{W_0,Cn}(\Relu;d,1)$ with $W_0$ and $C$ fixed and depending only on $d$.

  At a first glance,
  this theorem is very surprising to numerical analysts and approximation theorists since it is giving a rate of approximation ${\cal O}(n^{-2})$, $n\ge 1$, which is
  twice what standard approximation methods based on $n$ parameters give.  This indicates that the nonlinear manifold  $\Sigma_n:=\Upsilon^{W_0, C n}$ has 
  certain space filling properties in $X=C(\Omega)$.  While this seems like a great advantage of this manifold, recall that there are always
  one parameter manifolds which are dense in $X$, albeit not as neatly described as $\Sigma_n$.  But then,  we must throw in some caution.
  The theorem says that given $f\in K$, there is a mapping $a:K\to \R^n$ which selects the parameters $a(f)$ of the approximant that
  produces this exceptional approximation performance.  From our remarks in \S\ref{S:optimalperformance} on manifold width, the mapping $a$ cannot be continuous (note that the mapping $M$ is always continuous as will be discussed in more detail in the next section).   This shows a lack of numerical stability in the approximation process which yields Theorem \ref{T:Yarotsky}.  This means that we can expect that it will be very difficult to
  numerically find the parameters that attain the super convergence rate via a search over parameter domain.  On the other hand, if we are willing to allow a long enough search time with an a posteriori error   estimator,  we might be able to find such parameters.

  In spite of the negative comments just put forward, the theorem is intriguing and brings up several questions that we now discuss.
  The first natural question is in what generality does this super convergence hold.  We have already mentioned that Yarotsky proved it for multivariate functions of $d$ variables.  He also proved a general result which gives that the theorem holds for Lip $\alpha$ spaces, $0<\alpha\le 1$.   A generalization of this theorem is provided in \cite{lu2020deep}.  It shows that the set $K= U({\rm Lip} 1)$ can
  be replaced by the unit ball 
  $K$ of $C^s(\Omega)$, $\Omega=[0,1]^d$, for any  $s>0$.  However, in the latter presentation there is a loss of logarithm in that the proven approximation  rate  is $E(K,\Sigma_n)_{C(\Omega)}\le (\frac{\log n}{n})^{2s/d}$, $n\ge 1$.
  
  Next, let us remark that the results of 
  \S\ref{SS:VClimits} and Theorem \ref{L:VCdeep} give that for the model classes $K=U(C^s(\Omega))$ we have the lower bound
  \be 
  \label{Ylower}
  E(K,\Sigma_n)_{C(\Omega)}\ge c_0 n^{-2s/d},\quad n\ge 1.
  \ee 
  So, at least for the Lipschitz spaces, we have matching upper and lower bounds, and therefore a satisfactory understanding of the approximation properties of deep NNs for these classes.  
  
  \subsection{Super convergence for approximation in $L_p$}
  \label{SS:superLp}
  The above results were limited to approximation in $C(\Omega)$,  $\Omega=[0,1]^d$, and the  Sobolev spaces  $W^s(L_\infty(\Omega))$.  What happens when the approximation takes place in $L_p(\Omega)$,
  $1\le p< \infty$,  and what happens for general  Besov spaces that compactly embed in $L_p$?  We show in this section  that we can obtain super convergence results in this case as well by using results from the theory of interpolation spaces.
  
  \begin{theorem} 
  \label{T:superLp}
  We consider approximation in $L_p(\Omega)$, $1\le p\le\infty$, with domain  $\Omega=[0,1]^d$. Let  $\Sigma:=(\Sigma_n)_{n\ge 1}$, where $\Sigma_n:=\Upsilon^{W_0,Cn}(\Relu;d,1)$, $n\ge 1$, $W_0$  sufficiently large depending only on $d$, $C=C(s,d,\tau,p)$.  If $K:=U(B_q^s(L_\tau(\Omega)))$ is the unit ball of a Besov space above the Sobolev embedding line, then
  \be 
  \label{superLp}
  E_n(K,\Sigma)_{L_p(\Omega)}\le C[\log n]^\beta n^{- \theta s/d}, \quad n\ge 1,
  \ee 
 for any   $1\leq\theta<2-\frac{\tau^\ast}{\tau}$, with $\tau^\ast:=(s/d+1/p)^{-1}$, $\tau>\tau^*$,  and  $\beta$ depending  only on $s,d$ and $\theta$.
  \end{theorem}
  \noindent
  {\bf Proof:}  This is proved by using the K-functionals of interpolation theory. To keep the presentation simple and to just show the ideas of how
  this is done, we limit ourselves to proving one result of the above form when $d=1$ and $s=1$ with the approximation taking place in $L_\infty$.  Instead of Besov balls, we use   the unit balls 
  $K_\tau :=U(W^1(L_\tau(\Omega)))$, $1\le \tau\le \infty$ of the Sobolev spaces. After presenting this example, we give in Remark \ref{R:sLp}
 an outline of the proof of the general result stated in the theorem.

  We know the following two estimates
  \be
  \label{twoestimates}
   E_n(K_1,\Sigma)_{L_\infty(\Omega)}\le Cn^{-1},\quad E_n(K_\infty,\Sigma)_{L_\infty(\Omega)}\le Cn^{-2}, \quad n\ge 1,
  \ee 
where the first  was given in \S \ref{SSS:d1classical} (the approximant in this theorem can be viewed as element from  $\Upsilon^{3,n}(\Relu;1,1)$) and the second  is the super convergence result of Yarotsky, see Theorem \ref{T:Yarotsky}, with    $\Sigma_n:=\Upsilon^{11,16n+1}(\Relu;d,1)$.
  From interpolation between the pair $W^1(L_1(\Omega))$ and $W^1(L_\infty(\Omega))$ (this is where K-functionals are used, see \cite{DS}), we know that whenever $f\in K_\tau$, for any $t>0$ there is a function $g\in W^1(L_\infty(\Omega))$ such that
  \be 
  \label{interp11}
  \|f-g\|_{W^1(L_1(\Omega))} +t\|g\|_{W^1(L_\infty(\Omega))} \le Mt^{1-1/\tau},
  \ee 
  with $M$ an absolute constant.  We take $t=1/n$ in going further.   Now, let $S$ approximate $(f-g)$ in $L_\infty(\Omega)$ with the accuracy of the first statement in \eref{twoestimates},
  and let $T$ approximate $g$ with the acccuracy of the second statement.  Then $S+T\in\Upsilon^{11,17n+1}$ and
  \begin{eqnarray} 
  \label{approxinterp}
  \|f-(S+T)\|_{L_\infty(\Omega)}&\le &\|f-g-S\|_{L_\infty(\Omega)} + \|g-T\|_{L_\infty(\Omega)}\nonumber \\
  &\le& C\{n^{-1}\|f-g\|_{W^1(L_1(\Omega))}+ n^{-2}\|g\|_{W^1(L_\infty(\Omega))}\}\nonumber\\
  &\le& C n^{-1} n^{-1+1/\tau}=Cn^{-2+1/\tau}.
    \end{eqnarray} 
  In this case $\tau^*=1$ , so this  is the desired inequality.
Moreover, since 
$\|f-(S+T)\|_{L_p(\Omega)}\le \|f-(S+T)\|_{L_\infty(\Omega)}$, we also have 
$$
E_n(K_\tau,\Sigma)_{L_p(\Omega)}\leq Cn^{-2+1/\tau}, \quad 1\leq p \leq \infty.
$$
\hfill $\Box$
 
\begin{figure}[h]
  \centering
\includegraphics[scale=.5]{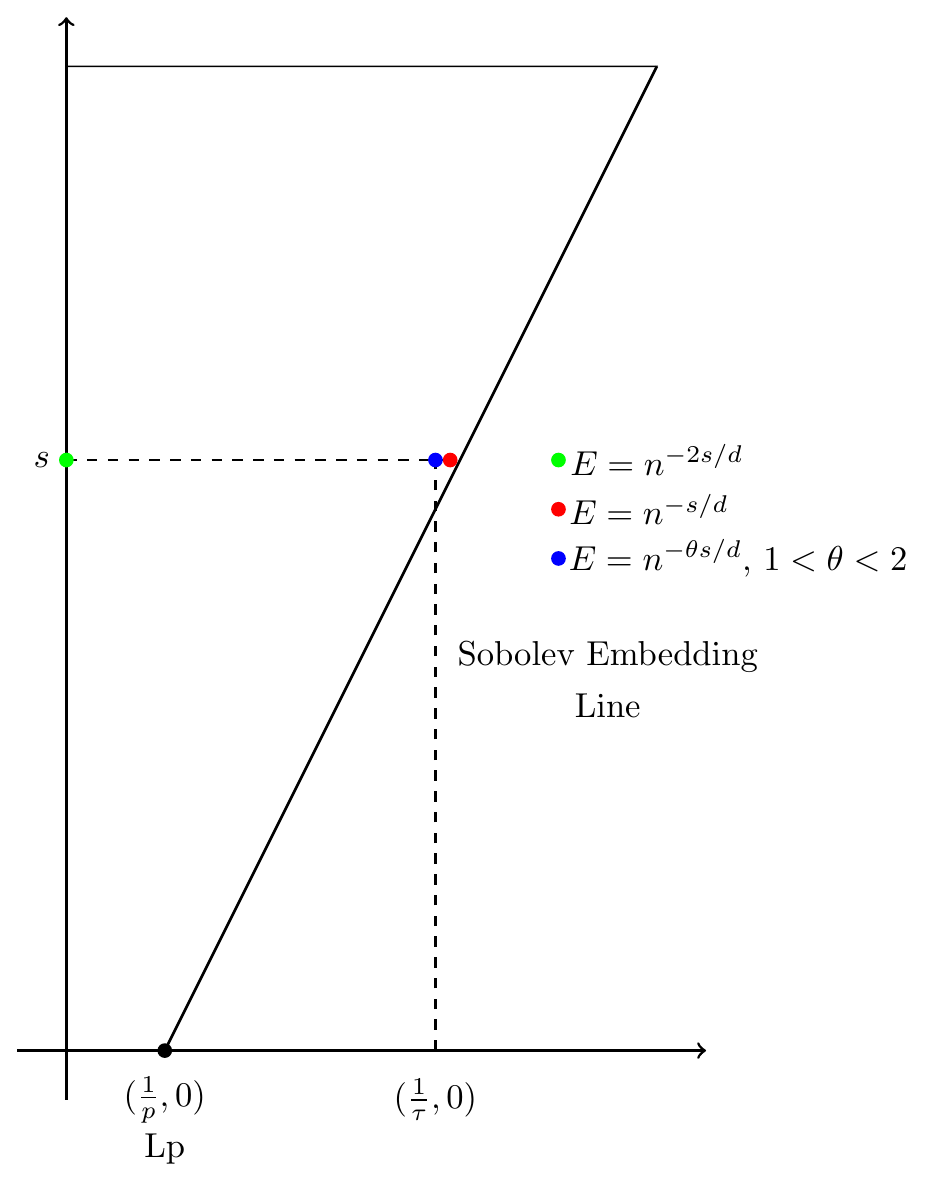}
\caption{Why we get general super convergence by using interpolation theory.
All error rates $E$ are modulo powers of logarithms when  $s>1$.}
\label{Figsuper}
\end{figure}

 \begin{remark}
 \label{R:sLp} 
 
 We outline the changes necessary to  prove the general case in the statement of the theorem. Now, we want to measure approximation error in $L_p(\Omega)$, $1\le p <\infty$, not just $C(\Omega)$.  
 Of course, the error of approximation in $L_p(\Omega)$ of a function $f$  is smaller than that in $C(\Omega)$.   
We use analogues of  
 \eref{twoestimates}    for approximation in $L_p(\Omega)$ and   two Besov balls.  The first is  $K_0=U(Z_0)$, $Z_0=B^s_{\infty}(L_{\tau_0}(\Omega))$ where we use Theorem \ref{T:Besov} to
 get the approximation rate $[\log_2 n]^{\beta_0} n^{-s/d}$.  Here, we can choose $\tau_0>\tau^*$ so that we are as close to
 the Sobolev embedding line as we want (but not on it).  The second inequality is the super convergence result for $K_1 =U(Z_1)$, $Z_1=C^s(\Omega)$.  For this, we use the generalization of Theorem \ref{T:Yarotsky}, as given in \cite{lu2020deep}, which gives the super  approximation rate $[\log_2 n]^{\beta_1}n^{-2s/d}$. We now interpolate between $Z_0$ and $Z_1$ to obtain the theorem for approximation in the fixed $L_p(\Omega)$ space. The reason we have the given restriction on $\theta$  is because we cannot take $Z_0$ directly on the Sobolev embedding line.  Figure \ref{Figsuper} may be useful for the reader to understand this theorem.
  \end{remark}
  
  \subsection{A summary of known approximation rates for classical smoothness spaces}
  \label{SS: summaryBesov}
  Let us summarize what we know about the  optimal approximation rates when approximating functions from Besov (and Sobolev) model classes using
  the outputs $\Sigma_n$ of deep neural networks, $\Sigma_n:=\Upsilon^{W_0,Cn}(\Relu;d,1)$, where $W_0$ is fixed, large enough, and  depending only on $d$, and $C$ depends on $d$ and the model class.  Given a value of $p$ with $1\le p\le \infty$, recall that 
  \be 
  \label{Enagain}
  E_n(K,\Sigma)_{L_p(\Omega)}:= \sup_{f\in K}\dist (f,\Sigma_n)_{L_p(\Omega)},\quad \Omega:=[0,1]^d.
  \ee 
  We want to address what we know regarding the following problem.

  \vskip .1in
  {\bf Problem 6:}  For  each model class $K$ which is the unit ball of a  Besov space $B_q^s(L_\tau (\Omega))$ which lies above the Sobolev embedding line for $L_p(\Omega)$, determine asymptotically matching upper and lower bounds for $E_n(K)_{L_p(\Omega)}$, $n\ge 1$.
  \vskip .1in
  Even for the most favorable case $p=\infty$, we only have a satisfactory answer to this question when $0<s\le 1$, in which case  the optimal rate is   $n^{-2s/d}$, $n\ge 1$.
  The above results on super convergence provide the upper bounds. The lower bounds follow from the   derivation of lower bounds on approximation rates using VC  dimension, given in \S\ref{SS:VClimits}.  Going further with the case $p=\infty$, the above results only provide a complete description of approximation rates when $s\le 1$ because of the  the appearance of a logarithm in the extension of Yarotsky's results given in \cite{lu2020deep}.
 
 When we move to the case $p<\infty$, the situation is even less clear.  First, Theorem \ref{T:superLp} does give a super rate.  However,
 we have no corresponding lower bounds that come close to matching this rate because we can not use VC dimension theory for $L_p$ approximation.
 In summary,  for all Besov spaces that compactly embed into $L_p(\Omega)$, we obtain
 error bounds for approximation in $L_p(\Omega)$ strictly better than classical methods.  What is missing vis a vis {\bf Problem 6} is what are the best bounds
 and how do we prove lower bounds for  approximation rates in $L_p(\Omega)$, $p\neq \infty$.

 \subsection{Novel model classes}
 \label{SS:novel1}
 
 While the performance of NN approximation on the classical smoothness spaces is an intriguing question that deserves a full and complete answer,
 we must stress the fact that such an answer will not provide an explanation for the success and popularity of NNs in their current domains of 
 application, especially in deep learning.  Indeed, the problems addressed via deep learning typically have the feature that the functions to
 be captured are very high dimensional, that is,  the input dimension $d$ is very large.  Since all of the classical model classes built on smoothness have large entropy and suffer the curse of dimensionality
 as $d$ gets large, they are not appropriate model classes for such learning problems.  This amplifies the need to uncover new model classes that
 do not suffer the curse of dimensionality, that are well approximated by outputs of NNs, and are a good match for the  targeted application.  We must say that little is formally  known in terms of rigorously defining new model classes in high dimension, showing that they have reasonable entropy bounds, and then analyzing their approximation properties by NNs.   However, several ideas have emerged as to how such model classes may be
 defined.  We mention some of those ideas here with the intention  to outline a  road map of how to  possibly proceed with defining model classes in high dimensions.

 \subsubsection{Comments on the curse of dimensionality}
 \label{SSS:curse}
 First, let us say a few words about the curse of dimensionality. One frequently hears the claim that a certain numerical method \lq breaks the   curse of dimensionality\rq.  There are two components to such a statement.  The first is that the numerical problem under study is such that it can be solved
 in high dimensions without suffering adversely from dimensionality.  The second is that a particular numerical method has been found that actually does the job.
 
 In the setting of numerical methods for  function approximation, the first statement has to do with the model class assumption on $f$, or the model class information that can be derived about $f$ from the context of the problem.  For example, when solving a PDE numerically, the model class information is usually given by a regularity theorem for the solution to the PDE.
  In other words, it is 
 the model class $K$ that  determines whether or not the problem is solvable by a numerical method that avoids the curse of dimensionality. 
 
 Heuristically, it is thought that the crucial  factor on whether or not a given model class $K$ suffers from the curse of dimensionality is its Kolmogorov entropy in  the metric where the error is to be measured, see \S \ref{SS:widths} for the definition of this entropy and the entropy numbers $\e_n(K)_X$.  There is not always a clear cut mathematical proof that entropy is indeed the deciding factor. This  lack of clarity stems from our vagueness in describing what is an allowable numerical method.  This returns  us back to the use of space filling manifolds in approximation.  We have already noted that such manifolds
 have the capacity to approximate arbitrarily well.  But are they a fair method of approximation?  Implementing such a manifold numerically as an approximation tool  requires an inordinate amount of computation.  So really, the computational time to implement the numerical method is an issue.  This is well known in the
 numerical analysis community,  but seems to be not treated sufficiently well in the learning community.  The latter would involve statements
 about how many steps of a descent algorithm are necessary to guarantee a prescribed accuracy.

We have touched on this subject in \S \ref{SS:stablerate}, where we have introduced  stable methods of approximation.  The introduction of stability was made precisely to quantify when a numerical method could be implemented within a reasonable computational budget. Under the imposition of stability in manifold approximation, we have shown that indeed the entropy of $K$ governs optimal approximation rates.

   Regarding the second factor, the question is whether we can put forward a concrete numerical scheme which can approximate the target function with a computational budget which does
 not grow inordinately with the dimensionality $d$.  In this sense, it is not only an issue of how well we can approximate a given $f$ using a specific
 tool $\Sigma:= (\Sigma_n)_{n\ge 1}$, but whether we can find an approximant within a reasonable  computational budget.

\subsubsection{Model classes in high dimension}
\label{SSS:hdmodel}
With these remarks in hand, our quest is to find appropriate model classes for high dimensional functions which have reasonable entropy when $d$ is large and yet match intended
 applications.  In this context, it is allowable for the entropy of the model class to grow polynomialy with $d$, but not exponentially.

 The search for appropriate high dimensional model classes has carried on independently of deep learning or NN approximation, since  
  it has always been  a driving issue whenever we are dealing with  high dimensional approximation.
 We next mention some of the ideas that have emerged over the recent decades on how to possibly define high dimensional model classes and how these ideas intersect with NN approximation.
 \vskip .1in
 \noindent
 {\bf Model classes built on sparsity:}  The idea  of using  sparsity to describe high dimensional model classes appeared  largely in the context of signal/image processing.
 The simplest example is the following.  Assume $\{\phi_j\}_{j\ge 1}$, with $\|\phi_j\|_X=1$, is an unconditional basis in a  Banach space $X$ of
  functions of $d$ variables. So, every 
  $f\in X$ has a unique representation
  \be 
  \label{comprepf}
  f=\sum_{j=1}^\infty \lambda_j(f)\phi_j,
  \ee 
where $\lambda_j$ are linear functionals on $X$ and the convergence in \eref{comprepf} is absolute. Here, the reader may assume that $X$ is an $L_p$ space to fix ideas.  The space $X$ defines the norm where we will measure
  performance (error of approximation).   Given any $q\le 1$, let $K_q$ consist of all functions $f\in X$ such that
  \be 
  \label{compK}
 f=\sum_{j=1}^\infty \lambda_j(f)\phi_j,\quad \sum_{j=1}^\infty |\lambda_j(f)|^q\le 1.
 \ee 
 If one wishes to approximate functions from $K_q$, the most natural candidate is $n$-term approximation using the basis $(\phi_j)_{j\ge 1}$. Let $\Sigma_n$ be the (nonlinear) set consisting of all functions $S=\sum_{j\in\Lambda} a_j\phi_j$, $\#(\Lambda)\leq n$.  It is a simple exercise to show that
 \be 
 \label{ntermerror}
 E(K_q,\Sigma_n)_X\le C_qn^{-1/q+1},\quad n\ge 1.
 \ee 
 Note that the Besov model classes take a form similar to \eref{compK}  because of their characterization by atomic decompositions using splines or wavelets, see \S \ref{SS:atomic}.  
 There are numerous generalizations of this notion of sparsity. For example, one can replace the unconditional basis by a more general set $\cD$ of functions,  which form a frame or  a dictionary.  
 
 Even though they give approximation rates that do not depend on the number of variables $d$, model classes built on sparsity are not necessarily immune to the curse of dimensionality because the basis or dictionary is infinite.
 To avoid this, one has to impose other conditions on the sequence of coefficients $(\lambda_j(f))_{j\geq 1}$ that allows one to truncate the sum to a finite set of indices when seeking an $n$-term approximation.
 This is often imposed by putting mild decay assumptions on these coefficients. The other central issue is whether the model class built on sparsity matches the intended application.  That is, there should be some justification that
 the sparsity class is a natural assumption in the application area.

 We have already seen an example of using sparsity in terms of a dictionary in discussing NN approximation when we introduced  the  Barron  class.  The Barron class appears as a natural model class when using shallow neural networks as an approximation tool.
  The neat  thing about the Barron class is that its definition was not made in terms of a dictionary but rather classical notions such as Fourier transforms. Generalizations of Barron classes to deeper networks is given in \cite{EMW}.  Then it was shown to be a sparsity class for a suitable dictionary of waveforms.
  
  \vskip .1in
  \noindent{\bf Model classes built on composition:}  Since  NNs are built 
  on the composition of 
  functions, it is natural to
  try to define model classes based on such compositions.  The basic idea is that the model class should consist of functions $f$ 
  with the representation $f=g_1\circ g_2\circ \cdots \circ g_m$,
  where  $g_k$, $k=1,\dots,m$, are
  simple component functions.  This approach is studied, for example, in
  \cite{MP,Shencomp,Schmidt}.
  
  The key question in such an approach is what assumptions should be placed on the component functions.  One expects to build the model
  class in a hierarchical fashion by showing that when $g_1$ and $g_2$ are well approximated then so is their composition. Let us consider for a moment the simple setting of approximating in the univariate uniform norm $\|\cdot\|_{C(\Omega)}$, $\Omega=[0,1]$.  Given $g_1,g_2$  and approximants $\hat g_1$ and $\hat g_2$,
  the simplest inequality for how well $\hat g_1\circ \hat g_2$ approximates $g_1\circ g_2$ is
  \begin{eqnarray} 
  \nonumber
  \|g_1\circ g_2-\hat g_1\circ \hat g_2\|_{C(\Omega)}&\le & \|g_1\circ g_2-\hat g_1\circ g_2\|_{C(\Omega)}+ \|\hat g_1\circ g_2-\hat g_1\circ \hat g_2\|_{C(\Omega)}\nonumber \\
 &\le&
  \|g_1-\hat g_1\|_{C(\Omega) }+\|\hat g_1\|_{{\rm Lip}\ 1}\|g_2-\hat g_2\|_{C(\Omega)},
  \end{eqnarray}
 which points to the observation that  formulations of such model classes will probably  involve mixed norms.

  \vskip .1in
  \noindent
  {\bf Model classes built on self similarity:}  Let us continue with the last example of the composition $g_1\circ g_2$.
  If $g_2$ is a CPwL function
  (as is the case of outputs of  ReLU NNs), then as the input variable $t$ traverses $[0,1]$, the composition traces out scaled copies of $g_1$ or parts of it. For example, if $g_2$ is the saw tooth function $H^{\circ L}$ of Figure \ref{F16}, then we trace out multiple copies of $g_1$.  The composition is, therefore,   a self similar function.
  This self similarity is prevalent in  outputs of deep NNs and  has been used
  to show that certain functions such as the Weierstrass nowhere differentiable function are well approximated by  outputs of deep NNs.  There are
  even classes of functions, generated by dynamical systems, which are  efficiently approximated by outputs of deep NNs.  So, it is natural to try and build model classes using  self similarity or fractal like structures,  and then  show that its  members are well approximated by deep NNs.  Examples of such 
  univariate function classes are given in \cite{daubechies2019nonlinear}, including the so-called Tagaki class. In higher dimension, it is shown in \cite{DSD} that
  the characteristic functions $\chi_S$ of certain fractal sets are also efficiently approximated by the outputs of deep networks.
  This may relate to the success of deep learning in classification problems.

\vskip .1in
\noindent
{\bf Model classes built on dimension reduction:}  A common high dimensional model class with reasonable entropy is the set of functions with anisotropic smoothness. These  functions   depend non democratically on their variables, that is, certain variables are more important than others, see e.g. \cite{DPW2}. This is a dominant theme in numerical methods for PDEs, where notions of hyperbolic smoothness classes and numerical methods built on sparse grids or tensor structures arise.
 
Another prominent example is a  model class viewed as low dimensional manifolds in a high dimensional ambient space. Since our approximation tool is itself a parameterised manifold, these model classes seem like  a good fit for NN approximation.  This is related to the viewpoint that  the 
NN output is
an adaptive   partition/filter design as expressed in \cite{baraniuk}.

  \section{Stable approximation}
  \label{S:StableRelu}
  Up to this point, we have mainly been interested in how well we can approximate a target function $f$ by the elements of the sets
  $\Upsilon^{W,L}(\Relu;d,1)$.  The results that we have obtained do not usually provide an actual procedure that could be implemented 
  numerically.  In this section, we discuss in more detail issues surrounding the construction of numerical approximation procedures and whether we can guarantee their stability. This section builds on the general discussion in \S\ref{S:optimalperformance} which the reader needs  to keep in view.
   
Here, we   measure error in the norm of $X=L_p(\Omega)$, $\Omega=[0,1]^d$, $1\le p\le \infty$.  We let
  $\Sigma_n$, $n\ge 1$, be   the ReLU sets $\Upsilon^{W,L}(\Relu;d,1)$ with the number of parameters $n(W,L)\asymp n$.  As usual, the two main examples that we have in mind are when $W=n$ and $L=1$ and secondly, when $W=W_0$ is fixed (depending on $d$) and $L=n$. Let $K$ be a model class in the chosen $L_p(\Omega)$.   
  
  We have mentioned  before that  any approximation method is described by two mappings  
  $$
  a_n:K\to\R^n;\quad M_n:\R^n\to \Sigma_n,
  $$
  where $a_n$ chooses the parameters of the network for a given
  $f\in K$, and $M_n$ describes how the neural network takes a vector $y\in\R^n$ of parameters and assigns the output $M_n(y)\in\Sigma_n$.
  Thus, the approximation to $f$ is the function $A_n(f)=M_n(a_n(f))$.  Notice that once we have decided to use NN with specific architecture for the method of approximation, the mapping $M_n$ is fixed and we do not get to choose it.
  
  We now wish to understand   two main issues:
  \begin{itemize}
  \item {\bf Stability Issue 1: }  How does 
  imposing  stability restrictions on the mappings $a_n$ and $M_n$ affect the approximation rates we can obtain?
 \item {\bf Stability Issue 2:}  How can we construct stable numerical algorithms for approximation?
 \end{itemize}
  
 We have already discussed {\bf Stability Issue 1} in some detail, see \S \ref{SS:stablewidths}.  We have seen that imposing stability limits
 the achievable approximation rates for NN approxcimation of a model class $K$ in the sense that the decay rate cannot be better than the entropy numbers $\e_n(K)_X$ of $K$.  Of course, this does not say we can necessarily achieve (with NNs) an approximation rate equivalent to $\e_n(K)_X$.
 However, this does give a benchmark for optimal performance.  This leads us to the following problem.
 \vskip .1in
 \noindent
 {\bf Problem 7:} What are  the optimal stable approximation rates for classical model classes such as Sobolev and Besov balls when using
 $\Sigma:=(\Sigma_n)_{n\ge 1}$ with $\Sigma_n:=\Upsilon^{W_0,Cn}(\Relu;d,1)$ as the approximation tool?  In other words, we want matching upper and lower bounds for stable approximation of these model classes.   A more modest question would be to replace stability by simply asking for continuity of these mappings.
 \vskip .1in
  
 Consider, for example, approximation in  $L_p(\Omega)$ with  $\Omega=[0,1]^d$ of the Besov balls $B_q^s(L_\tau(\Omega))$ that embed into $L_p(\Omega)$.    
The entropy of such a ball is known and gives the lower bounds ${\cal O}(n^{-s/d})$ for the best approximation rate by
stable method of approximation. However, we have not provided stable 
 mappings for NNs that achieve this approximation rate.  A similar situation holds when we assume only continuity of these maps.

  \subsection{Stability of $M_n$}
  \label{SS:stabilityM}
  As we have noted, when using NN approximation, the mapping $M_n$ is determined  by the architecture of the NN.  In this section, we discuss the stability of this mapping.  
  We always take $M_n$ to be the natural mapping
  which identifies the output $S\in \Upsilon^{W,L}(\Relu;d,1)$ with  the parameters that are the entries of the matrices and bias vectors of the NN, see \S\ref{SS:FCN}.
  We   identify these parameters with a point in $\R^n$ in such a way that the parameters at 
  layer $\ell$ appear before those for the next layer and the ordering  for each hidden layer is done in  the same way.
  
  It is easy to see that the  mapping $M_n:\R^n\to C(\Omega)$,
  $\Omega=[0,1]^d$, is continuous.  In fact, as we shall now show, it is a Lipschitz map on any bounded set of
  parameters, that is, $M_n$ is locally Lipschitz. To describe this, we need to specify a norm to be used for $\R^n$.  
  We take this norm to be the $\ell_\infty(\R^n)$ norm, that is,   $\|y\|_{\ell_\infty^n}:=\max_{1\le i\le n}|y_i|$.  This choice is not optimal for obtaining the best constants in estimates but it will simplify the exposition that follows.
  
  \begin{theorem}
  \label{T:MLip}
If $B$ is any finite ball in $\ell_\infty(\R^n)$, then 
$M_n:B\to C(\Omega)$ is a Lipschitz mapping, that is,
  \be 
  \label{MLip}
 \|M_n(y)-M_n(y')\|_{C(\Omega)}\le C\|y-y'\|_{\ell_\infty^n},\quad y,y'\in B,
  \ee 
  with the constant $C$ depending only on $B,W,L$, and $d$.
  \end{theorem}
  \noindent
  {\bf Sketch of Proof:}
 We will be a bit brutal and not search for the best constant in \eref{MLip}.  In what follows in this proof,   $C$ denotes a  constant depending only on $B,W,L$, and $d$,
and  may  change from line to line. For 
  $y,y'\in B$  we wish to bound $\|M(y)-M(y')\|_{C(\Omega)}$ by $\delta:=\|y-y'\|_{\ell_\infty^n}$. 
  For a vector valued continuous function $g$, we use $\|g\|$ to denote its $C(\Omega)$ norm, which is the maximum of the $C(\Omega)$ norm of its components.

  We denote by $\eta^{(j)}$,
  $j=1,\ldots,L$, the
  vector valued function of $x\in \R^d$ computed at level $j$ by the network with parameters $y$, and by $\eta'^{(j)}$ the corresponding vector of functions computed with parameters $y'$.   
 We can write
\begin{eqnarray} 
 \nonumber
  \eta^{(1)} = \Relu (A_0x+b^{(0)}),\quad  \eta'^{(1)} = \Relu (A'_0x+b'^{(0)}),
  \nonumber \\
  \eta^{(j+1)} = \Relu (A_j\eta^{(j)}+b^{(j)}),\quad  \eta'^{(j+1)} = \Relu (A'_j\eta'^{(j)}+b'^{(j)}),
  \nonumber
\end{eqnarray}
  where  $A_j$ is the matrix determined by $y$ to go from level $j$ to level $j+1$, and $b^{(j)}$ is the bias vector,
  $j=0,\ldots,L-1$.  Similarly,
  $A_j'$ and $b'^{(j)}$ correspond to the parameter $y'$. 
  
  Since $y,y'\in B$, all entries in the $A_j,A'_j,b^{(j)},b'^{(j)}$ are bounded.  Likewise, the matrix norms of  $A_j,A'_j$ as mappings  from $\ell_\infty$
  to $\ell_\infty$ are bounded. Also, we have,
  
  \begin{itemize}
  \item $\|A_0-A_0'\|_{\ell_\infty^d\to\ell_\infty^W}\le C\delta\quad$, 
  $\|A_j-A_j'\|_{\ell_\infty^W\to\ell_\infty^W}\le C\delta, \quad j=1,\dots,L-1$.
  \item $\|b_j-b_j'\|_{\ell_\infty^W}\le  \delta,\quad j=0,\dots,L.$
\end{itemize}
Using 
  \be 
  \nonumber
  \eta^{(j+1)}= \Relu[A_{j}(\eta ^{(j)}-\eta'^{(j)}) + A_j(\eta'^{(j)})+b_j], 
  \ee 
  and the fact that $\Relu(\cdot)$ is a Lip 1 function, one derives 
  \begin{eqnarray}
  \nonumber
  \|\eta^{(j+1)}-\eta'^{(j+1)}\| &\le& \|A_j\| \|\eta^{(j)}-\eta'^{j)}\|+\|A_j-A_j'\|\|\eta'^{(j)}\| + \|b_j-b_j'\|\nonumber\\
  &\le & C\{  \|\eta^{(j)}-\eta'^{j)}\|+\delta \|\eta'^{(j)}\|+\delta\}.
  \nonumber
  \end{eqnarray} 
  One then proves by induction that $\|\eta'^{(j)}\|\le C$, $j=0,1,\dots,L$, and that
  \be 
  \nonumber
  \|\eta^{(j)}-\eta'^{(j)}\|\le C\delta, \quad j=1,2,\dots,L.
  \ee
 The final step is that
 \be 
 \nonumber
 \|M_n(y)-M_n(y')\|_{C(\Omega)}\le (\|\eta^{(L)}-\eta'^{(L)}\|+\delta),
 \ee 
 which gives the theorem.\hfill $\Box$
 
  \begin{remark}
  \label{R:stableMn}
  A closer look at the above estimates shows that the  Lipschitz constant for $M_n$ can be controlled if we take $B$ as a small ball 
  around the origin.  The size of the ball is chosen so that each of the matrices $A_j,A'_j$ have small norm. To do this, the required  size of the ball
  gets smaller as $W$ gets larger. 
  \end{remark}
  
  \subsection{Stability of $a_n$}
  \label{SS:stablea}
  With the above analysis of $M_n$ in hand, we see that the stability of a
  NN approximation method rests on the properties of the parameter selection $a_n$.    It is of interest to understand whether the most common methods of parameter selection based on gradient descent provide any stability.  We discuss this issue in \S \ref{SSS:select}.  For now, we limit ourselves to recalling our discussion on how imposing stability limits
  approximation rates and when we know stable methods.  For the discussion that follows, we limit ourselves to approximation using
  $\Sigma_n=\Upsilon^{W_0,n}(\Relu;d,1)$ with $W_0$ fixed.  Many of the same issues we raise concerning stability appear in approximation using shallow networks.

  Let us first observe that the parameter selection procedures that generate the super rates of convergence for Besov and Sobolev classes cannot be continuous because of \eref{Sobwidth}.
  If we require that the mappings $a_n$ are only continuous and consider approximation in $L_p(\Omega)$, $\Omega=[0,1]^d$, then we can never  attain rates of approximation better than ${\cal O}(n^{-s/d})$ for the unit ball of   any Besov space $B_q^s(L_\tau(\Omega))$ that embeds compactly into $L_p(\Omega)$.
  The only cases where we know that we can actually attain this rate is when $\tau\ge p$. In these cases, there are linear spaces, such as FEM spaces, contained in $\Sigma_n$ that provide this rate and the  approximation can be done by a linear operator. So the following problem is not solved except for very special cases.
  \vskip .1in
  \noindent
  {\bf Problem 8:} Consider the approximation of the unit ball of a Besov space 
  $B_q^s(L_\tau(\Omega))$ compactly embedded in $L_p(\Omega)$ using the manifold
  $\Sigma_n$. Give matching upper and lower bounds for the approximation rate in the case $a_n$ and $M_n$ are Lipschitz mappings.  Similarly,
  determine upper and lower bounds when the parameter selection mapping $a_n$ is continuous.
  
  A question closely related to stability is whether one can approximate well under the very modest restriction that $a_n$ is bounded.   Recall that boundedness
  helps us with $M_n$ as well (see the above discussion).  The issue of what approximation rates are possible when  one  imposes boundedness on $a_n$ was studied in detail in
  \cite{BGKP}.  The motivation in that paper was different from ours in that they were interested in NN approximation from the viewpoint of encoding.
  However, there is an intimate connection with stability as we have just discussed.

  \section{Approximation from data}
  \label{S:data}
  Thus far, we have limited ourselves to understanding the approximation power of neural networks.  The approximation rates we have obtained
  assumed full access to the target function $f$.  This scenario does not match the typical application of NN approximation to the tasks of learning.
  In problems of learning,  the only information we have  is  data observations of $f$.  Such data observations alone do not allow
  any rigorous  quantitative guarantee of how well $f$ can be recovered, that is, how accurately  the behavior of $f$  at new points can be predicted.  What is needed for the latter is additional information about $f$, which we have referred to as model class information.  The  model class information is   an assumption
  about $f$ that is often not provable but based more on heuristics about the application area.

  Learning from data is a vast area of research that cannot be covered in any detail in this exposition.   So, we limit ourselves to pointing  out some
  aspects of this problem and how they interface with the theory of NN approximation that we have discussed so far. Obviously, any performance guarantees derived in the learning setting must necessarily be worse than those for approximation, where full information about $f$ is assumed.
  Thus, an important issue is to quantify this loss in performance.

  The most common setting for 
  the learning problem is a stochastic one, where it is assumed that the  data is given by random draws from an underlying probability
  distribution.   However, it is useful to consider the deterministic setting as well since it sheds some light on the stochastic formulation
  and the type of results that we can expect.
  
  \subsection{Deterministic learning; optimal recovery}
  \label{SS:optimalrecovery}
  In this section,  we wish to learn a function $f$ which is an element of a Banach space $X$. Our goal is to recover $f$ from some finite set of data observations.  We assume that the data observations are in the form of bounded
  linearly independent
  linear functionals applied to $f$.  Thus, our data takes the form
  \be 
  \label{Data}
  (\lambda_1(f), \ldots,\lambda_m(f))\in \R^m,
  \quad 
  \lambda_j\in X^*,\quad j=1,\ldots,m,
  \ee  
  where  $X^*$ is the dual space of $X$. 
  As we have pointed out numerous times, to give quantitative results on how well $f$ can be recovered requires 
  more information about $f$ which we call model class information, i.e., information of the form $f\in K$, where
  $K$ is a compact set in $X$.
  When we inject the model class assumption that $f\in K$, we have the question of how
  accurately we can recover $f$ from the two pieces of information, the data and the model class.
  We shall present the functional analytic view of this problem which is known as {\it optimal recovery}. It will 
  turn out that the optimal recovery problem is not always amenable to a simple numerical method for the recovery of
  $f$.  Nevertheless, this viewpoint will be useful in motivating  specific numerical methods and analyzing how
  well they do when compared with the optimal solution.

  \subsection{Optimal recovery in a Hilbert space}
  \label{SS:optimalrecH}

  We shall restrict our development here to the most popular setting where $X=H$ is a Hilbert space.  The reader interested in
  the more general Banach space setting can consult \cite{DPW}.  In the   Hilbert space setting, each of the functionals $\lambda_j$  has a representation 
  \be 
  \nonumber
  \lambda_j(f)=\langle f,\omega_j\rangle,\quad \omega_j\in H,\ j=1,\dots, m,
  \ee   
  which is referred to as the Riesz representation of $\lambda_j$.  The functions $\omega_j$ span an $m$  dimensional subspace
  \be 
  \nonumber
  W:={\rm span}\{\omega_j\}_{j=1}^m
  \ee  
  of $H$.
  We can assume without loss of generality that the $\omega_j$'s are an orthonormal system. From the given data, we can find the projection 
  $$
  w:=P_W f
  $$
  of $f$ onto $W$. We think of $w$ as the given data.  
  
  Now, let us assume in addition that $f$ is in a certain model class $K$, and ask what is the best approximation
  (with error measured in the norm of $H$) that we can give to $f$
  based on this information, i.e., the data and the model class information.
  One may think that the best we can do is to take $P_Wf$ as the approximation. However,  this is not the case since the information that $f\in K$ allows us to say something about the projection of $f$ onto the orthogonal complement $W^\perp$ of $W$.

  Indeed, the model class information will allow us to
  give a best  approximation to $f$ from the available information (model class and data $w$) as follows. 
  Let 
  $$
  K_w:=\{g\in K\,:\,P_Wg=w\}.
  $$
  Then, the membership of $f$ in $K_w$ is the totality of information we have about $f$.  The best approximation
  to $f$ is now given by the center of the set $K_w$.  Namely, let $B:=B(K_w)$ be the smallest ball in $H$ which contains $K_w$.  This ball is referred to as the {\it Chebyshev ball}, its center $b_w\in H$ is called the {\it Chebyshev center}, and its radius $R_w$ is the {\it Chebyshev radius}.
  The best approximation we can give to $f$ is to take $b_w$ as the approximation and the error that will ensue is $R_w$.
  The function $b_w$ is the {\it optimal recovery} and $R_w$ is its {\it error of optimal recovery}.

   Let us reflect a bit on the above optimal solution.  Every function in $K_w$ is a possibility for $f$.  From the information presented to us (model class plus data), we do not know which of these functions is the desired $f$.
  So, we do the best we can to approximate all of the possible $f$'s,  which turns out to be the Chebyshev center.
   Each 
   $g\in K_w$ (the possibilities for  
   approximants of $f$) is of the form $w+\eta$, where $\eta$ is in the null space $\cN=W^\perp$.  So, in essence, we are trying to find the $\eta \in W^\perp$  that we can add to $w$ so that the sum $w+\eta\in K$.
   \vskip .1in
   \begin{remark}
   \label{R:optrec}
   Notice that if we find any $\eta\in \cN$ such that $w+\eta$ is in $K$, then we have essentially solved the problem since $\hat f:=w+\eta$  approximates $f$ to accuracy at worst  $2R_w$.  Such an $\hat f$ is called a near best solution.
   \end{remark}
   \vskip .1in

  The above description of optimal recovery, despite being elegant and optimal, is not very useful in constructing a numerical procedure since the Chebyshev ball is difficult to find numerically.  Also in practice, we often are not sure what is the apropriate model class $K$ in a given setting. However, optimal recovery  is still a good guide for the development of  numerical procedures.

  There are two standard approaches to developing numerical algorithms for optimal recovery. The first one is to numerically generate a recovery  through least squares minimization with a constraint that enforces the model class assumption. We will not 
  engage this approach here but simply mention that several elegant  results show that for certain model classes these optimization problems have exact solution in NN spaces, especially the ones with a single hidden layer.  We refer the reader to
  \cite{unser}, \cite{PN},
  \cite{OWSS},
  \cite{SESS} 
  for the most recent results using this approach.

    The second approach, which is more  closely tied to approximation, is to replace $K$ by a simpler set $\hat K$ which is less complex than $K$, and yet accurate. One then solves the optimal recovery problem on the simpler surrogate model class $\hat K$. We discuss this approach in the following two sections.
    
    \subsection{Optimal recovery by linear space surrogates}
    \label{SS:ORlinear}
    The usual approach to finding a surrogate $\hat K$ for $K$ is to approximate $K$
  by a linear space   of dimension $n$, or more generally, a nonlinear manifold  $\Sigma_n$, with $n$ the number of parameters needed for its description.  If we know that $\Sigma_n$ approximates $K$ to accuracy $\e_n$ (here is where our error estimates for approximation are useful),  we then can replace $K$ by
  \be 
  \label{surrogateK}
  \hat K:=\{h\in H:\ \dist (h,\Sigma_n)_H\le \e_n \}.
  \ee 
 Clearly, $K\subset \hat K$.
  Usually, we also have some knowledge on the norm $\|f\|_H$ 
  for functions $f\in K$  and this can be used to trim the set $\hat K$ even further.
  
  Once a surrogate $\hat K$ has been chosen,  we solve the optimal recovery problem for $
  \hat K$ in place of $K$ by using Chebyshev balls for $\hat K_w$ as described above. As we shall now see, we can often
  solve the optimal recovery problem for the surrogate exactly by a numerical procedure.

  We assume for the time being that  $\hat K$ is given by \eref{surrogateK}    with $\Sigma_n$  a linear space of dimension $n<m$.  In this case, the problem is a much simpler recovery problem than the one  for $K$, and optimal recovery has an exact solution that we now describe, see \cite{BCDDPW1}. Let us define by 
  $H_w:=\{h\in H:\ P_Wh=w\}$, that is, $H_w$ is the set of all functions in $H$ which satisfy the data. Since $f\in K$, 
  $K\subset \hat K$, and $P_Wf=w$, we see that   $\hat K_w:=H_w\cap\hat{K}$ is non-empty. The center of the Chebyshev ball $B(\hat K_w)$ for $\hat K_w$ is the point $u^*(w)\in H_w$  which is closest to $\Sigma_n$, that is
  \be 
  \nonumber
   u^*(w):= \argmin_{h\in H_w} \dist(h,\Sigma_n)_H.
  \ee 
  The function $u^*(w)$ is found as follows.
  One solves the least squares problem
  \be
 \nonumber
  v^*(w):= \argmin_{v\in\Sigma_n} \|P_Wv-w\|_H,
  \ee 
  and then $u^*(w)= w+P_{W^\perp} v^*(w)$,  where $W^\perp$ is the orthogonal complement of $W$ in $H$ (the null space of $P_W$).  One can also compute the Chebyshev radius $\hat R_w$  of $B(\hat K_w)$ as
  \be  
  \nonumber
  \hat R(w)=\mu(\Sigma_n,W)_H(\e_n^2-\|u^*(w)-v^*(w)\|_H^2)^{1/2},
  \ee 
  where 
  \be 
 \nonumber
  \mu(\Sigma_n,W)_H:= \sup_{\eta\in W^\perp}\frac{\|\eta\|_H}{{\rm dist}(\eta,\Sigma_n)_H}.
  \ee

  Here are a few remarks to put the above results into context. 
  
  \begin{remark}
  \label{R:data2}
  The quantity $\mu(\Sigma_n,W)$  is the reciprocal of the cosine  of the angle between the
  two spaces $\Sigma_n$ and $W$.  It   reflects the quality of the data relative to $\Sigma_n$.
  This number will be large when the data is not well positioned relative to the linear space $\Sigma_n$.
  In particular,  it will always be infinite whenever the dimension $n$ of $\Sigma_n$ is larger than $m$.  This is because there will always be elements from $\Sigma_n$ in the null space of $P_W$ and hence there will be points in $\hat K_w$ that are arbitrarily far apart in this case.
  \end{remark}
  
  \begin{remark}
  \label{R:data3}  The above results give a bound for the performance of least squares, see \cite{BCDDPW1}.  Namely, 
  given data 
  $w_j^*=\lambda_j(f)$, $j=1,\ldots,m$, for some $f\in H$, let
  \be 
  \nonumber
  S^*:= \argmin_{S\in\Sigma_n }\sum_{j=1}^m[
  w^*_j-\lambda_j(S)]^2.
  \ee 
  Then, for any $f\in H$  which satisfies the data, we have
  \be 
  \nonumber
  \|f-S^*\|_H\le \mu(\Sigma_n,W)_H {\rm dist}(f,\Sigma_n)_H,
  \ee 
  and this bound cannot be improved in the sense that there are always $f\in H$ for which we have equality.
  \end{remark}
    
  The above analysis and remarks only apply to the case that $\Sigma_n$ is a linear space and $X=H$ is a Hilbert space.
  In the spirit of this paper, we would take $\Sigma_n$ to be $\Upsilon^{W,L}(\Relu;d,1)$, the outputs of a ReLU network which depends on roughly $n$ parameters.  We should choose the architecture to match $K$ as best as possible,  given the budget $n$ of parameters.
  
  Let us, for example, consider the case where $\Sigma_n:=\Upsilon^{W_0,n}$, $n\ge 1$,  with $W_0$ fixed, i.e., the case of a deep network with constant width, and continue to assume that
 $X=H$ is a Hilbert space.  We suppose that $\Sigma_n$ provides an approximation 
  with error
 \be 
 \nonumber
\dist (K,\Sigma_n)_H=\e_n.
 \ee 
 We view $\hat K:=\{h\in H: \ \dist(h,\Sigma_n)_H\le \e_n\}$ as a surrogate for $K$. 
 Note that  $K\subset \hat K$.
If $f,g\in \hat K_w:=\{h\in\hat K:  \,P_Wh=w\}$ then $\eta:=f-g\in W^\perp$, 
and 
\be 
\nonumber
\dist(\eta, \bar\Sigma_n)_H \le 2\e_n,
\ee 
where $  \bar \Sigma_n:= \Upsilon^{W_0+d+1,2n}(\Relu;d,1)$.
It follows that
 \be 
 \label{estChebradius}
 \|f-g\|_H \le  2\mu_n \e_n,\quad \hbox{where}\quad 
  \mu_n:= \sup_{\eta\in W^\perp} \frac{\|\eta\|_H}{{\rm dist}(\eta,\bar\Sigma_n)_H},\quad n\ge 1.
  \ee 
This tells us that the Chebyshev radius  $\hat R_w$ of  $\hat K_w$ (and thereby the Chebyshev radius $R_w$ of $K_w$)  satisfies
$$
R_w\leq  \hat R_w\le \mu_n\e_n.
$$
This is the same estimate as in the case when $\Sigma_n$ is a linear space, except that now we have to expand $\Sigma_n$ to $\bar \Sigma_n$ because of the nonlinearity of $\Sigma_n$.

We are left with finding an approximation to the Chebyshev center of $\hat K_w$ (and thereby $K_w$). For this we take
any $S^*\in\Sigma_n$ which satisfies
\be 
\label{lssigma1}
\|w-P_W S^*\|_H =\inf_{S\in\Sigma_n}\|w-P_W S\|_H\le \e_n,
\ee 
where the last inequality follows because 
$$
\|w-P_W S\|_H=\|P_Wf-P_WS\|_H\leq \|f-S\|_H,
$$
and we know $\dist(f, \Sigma_n)_H\leq \e_n$.
This is a least squares problem which does not necessarily have a unique solution.  However, we now show that any solution $S^*$ provides a good estimate for the Chebyshev center of $\hat K_w$.

Indeed, let us take any of its solutions $S^*\in\Sigma_n$ and 
consider 
$$
h^*:=w+P_{W^\perp}S^*\in H_w.
$$
With an eye towards \eref{lssigma1}, we see that 
$$\|h^*-S^*\|_H=\|w-P_WS^*\|_H\leq \e_n,$$
and thus $h^*\in \hat K_w$. Moreover, it follows from 
\eref{estChebradius} that for every $f\in \hat K_w$ we have
$$
\|f-h^*\|_H\leq 2\mu_n\e_n,
$$
and therefore, the ball of radius $2\mu_n\e_n$ with center $h^*$ contains $\hat K_w$. Thus, $h^*$ can be taken as an approximation
to the Chebyshev center of $\hat K_w$ (and thus to the Chebyshev ceneter of $K_w$).
A cruder, but less laborious  approximation to $f$ is provided by $S^*$, since
\be
\label{mn}
\|f-S^*\|_H\leq \|f-h^*\|_H+\|h^*-S^*\|_H
\leq 2\mu_n\e_n+\e_n, \quad f\in K_w,
\ee
and therefore
$$
\dist(K_w,S_n^*)_H\le (2\mu_n+1)\e_n.
$$

Inequality \eref{mn} can be reformulated in the 
following way. 
For any $f\in H$, the least squares solution $S_n^*$ for $w:=P_Wf$ provides an approximation to $f$ of accuracy
\be 
\label{estC}
\|f-S_n^*\|_H\leq (2\mu_n+1) \dist(f,\Sigma_n)_H,
\ee 
since the the above argument can be repeated with $\e_n=\dist (f,\Sigma_n)_H$.

Finally, note again that if $n\ge m$, then there will be elements of $\Sigma_n$ that interpolate the data and hence $\mu_n$ is infinite which renders the bound \eref{estC} useless.  Yet, this is the case of overparametrized learning which is often used in practice.  So, something must be added to least squares minimization in order to have viable results in the overparameterized case.  What this additional ingredient should be is  the subject of the next section.

  \section{Using Neural Networks for  data fitting}
  \label{S:overparameterized}

The typical setting  for supervised learning is to find an approximation of an unknown function $f$, given a training data set of its point values
\begin{equation}\label{E:data-def}
    \set{(x^{(i)}, f(x^{(i)})},\qquad x^{(i)}\in \R^d,\qquad f(x^{(i)})\in \R, \qquad i=1,\ldots, m.
\end{equation}
 We refer to the points $x^{(i)}$, $i=1,\dots,m$, as the {\it data sites}. Thus, the data observation functionals are point evaluations (delta functionals). In many applications, the dimension $d$ is very large.
For example, in classification problems for images, $d$ is the number of pixels in the images, typically somewhere in the range of $10^3$ to $10^6$, and for videos it is even higher.  The learning problem is then to numerically produce from this data  a function $\hat f$ that is in some sense a good  predictor of $f$ on
new unseen draws $x\in \R^d$.

In the preceding section, we described  a systematic approach to learning from data, called optimal recovery.  It  begins with two vital requirements: (i) a known model class $K$ to which  $f$ is assumed to belong, and (ii) a specific norm or metric in 
which the recovery of $f$ by $\hat f$ is measured.   In the optimal recovery formulation of the problem, a solid theory exists to describe the optimal solution via the Chebyshev ball.  The deficiency in this  approach is that the construction of numerical algorithms to generate a surrogate $\hat f$ may be a significant computational challenge.

Optimal recovery is not the viewpoint taken in the general literature on learning.  Rather, in the learning community,  the  data fitting task is formulated in a stochastic setting, where one assumes that the data comes from random draws 
of the data sites  $x^{(i)}$ with respect to a probability distribution, and the $f(x^{(i)})$'s are noisy observations of some unknown function $f$.  Performance is then evaluated  on new draws of data in the sense of probability or expectation of accuracy on these draws. This is commonly referred to as
{\it generalization error}.  Note that in this setting, there is no model class assumption on the function $f$ giving rise to the data, and so there can be no provable bound for the generalization error.  What is done in practice is to
give an empirical bound based on checking performance on a  lot of new (random) draws which are referred to as test data.

Traditionally, model class assumptions on the unknown function $f$  played a dominant role in the classical formulation and proof of a priori performance guarantees, see \cite{Lugosi}.
However, as noted in the previous paragraph, in the now dominant field of deep learning, where neural network approximation is  an important technique, one deviates from the classical setting of model class assumptions. Our goal in the sections that follow is    to understand what     role   approximation using neural networks plays in this new setting.

\subsection{Deep learning}
Deep learning is characterized by its ability to successfully treat very high dimensional problems, where one begins with inordinately large data sets and employs intensive computation for generating surrogates.  Its success  in handling high dimensional problems is  provided only by empirical
verification that the numerically created surrogate performs well on new draws of $x$.  A priori guarantees of performance are generally lacking.
In fact, performance is not typically  formulated under model class assumptions, which  in turn prevents such a priori analysis.
 The lack of a specific model class assumption is probably  due, at least in part,  to the high dimensionality $d$, since in this case it is often unclear what  appropriate  model classes should be.  Note, however, that since the data observations are point evaluations, a minimal assumption is that $f$ is in a Reproducing Kernel Hilbert Space (RKHS), although the specific RKHS  is not known or postulated.
 
 Another important feature of deep learning is its use of over parameterization in the search for a surrogate.  This runs in the face of classical learning which warns against overfitting the data because it leads to fitting the noise.
 
 \subsection{Possible model class assumption in high dimension}
 \label{SS:possiblemodel}

Before turning to the overparameterized setting, we wish to make a few remarks on possible viable model class assumptions that could be used towards providing a priori guarantees in deep learning.  One valid view point is that  the functions we are trying to recover do possess some
special properties; we just do not know what they are.

The fact that
neural networks are used quite successfully suggests that the functions  we are trying to learn are
well approximated by   neural networks.  If this is the case, then a natural model class assumption would be that $f$ is in an approximation
class $\cA^r((\Sigma_n)_{n\ge 1},X)$, which we recall consists of the functions $f$ for which
\be 
\label{aclassagain} {\rm dist}(f,\Sigma_n)_X\le Mn^{-r},\quad n\ge 1,
\ee 
where again 
there is the question what is the appropriate space $X$ in which  to measure error.  Here, $(\Sigma_n)_{n\geq 1}$ would be the family of spaces outputted by
the chosen NNs and $n$ would represent the number of their parameters. This underlines the importance of understanding the approximation performance
of neural networks in a rate/distortion sense, and, in particular, which functions are well approximated by neural network outputs.

 \subsection{Overparameterization}
 \label{SS: oevverp}
  We turn now to learning from data using overparameterized models.    When searching for an approximation  to $f$ from  a set
  of outputs  of a neural network with a given architecture, say $\Sigma_n=\Upsilon^{W,L}(\sigma; d;1)$, it is usually the case in practice that the number of trainable parameters, that is the number $n= n(W,L)$ of weights and biases, exceeds the number $m$ of data sites $x^{(k)}$,
\[
\#\text{data points} ~=~m ~\ll ~n~=~\#\text{parameters}. 
\]
In other words,  neural networks are usually \textit{overparameterized}.
 This means that  there are  generally infinitely many choices  of the parameter vector $\theta$ (of network weights and biases) so that 
 the 
 network 
 with these parameters  outputs a function $S(\cdot;\theta)$ that   fits (interpolates) the data, that is
\[
S(x^{(i)};\theta) = f(x^{(i)}),\qquad i=1,\ldots, m.
\]
 Characterizing exactly which interpolant is chosen by the numerical method is at the heart of learning via overparameterized neural networks.  
 In this section, we want to understand how this selection is done in practice and whether the selection  has an analytic interpretation.  In particular,
 there is the question of whether the numerical method itself is in a certain sense specifying a model class assumption.  If so, it would be important to unravel what this hidden assumption is.

\subsubsection{Selecting the interpolant by gradient descent}
\label{SSS:select}
The standard way of selecting an approximant $\hat f$  to $f$ 
in the practice of overparameterized deep learning using neural networks is to begin with a random starting guess $\theta^{(0)}$ for the parameters and  thereby specifying the first guess $S(\cdot;\theta_0)$ for a surrogate.  Successive approximations $S(\cdot,\theta^{(k)})$, for  $k=1,2,\dots$, are then generated  by applying a gradient descent (or stochastic gradient descent) to 
finding the minimum  of a loss $\cL$, which usually takes the form of an empirical risk, that is
\be 
\label{lossfunction}
\cL(\theta):=\sum_{i=1}^m (f(x^{(i))}-S(x^{(i)};\theta))^2.
\ee 
If the step sizes are appropriately chosen in the descent algorithm, then this procedure seems to work well in practice in that $S(\cdot,\theta^{(k)})$
 with $k$ large is an approximation to $f$ which generalizes well.
 Here, $\theta^{(k)}$ is the output parameter  of the gradient descent algorithm at the $k^{th}$ step.

 The above method for selecting a surrogate does not employ the  traditional remedy for working with approximation methods  that have the capacity to overfit the  data  which is to add a regularizer such as $\ell_1$ or $\ell_2$ penalty on the parameter vector in the iteration.   The effect of incorporating such regularizers is studied for example in   \cite{SESS,OWSS,PN}.  The main conclusion of the above papers is that one can view the addition of a constraint as a model class assumption.  However, it is important to note that adding such a regularizer is not usually done in NN practice. While a weak regularization is sometimes employed, empirical evidence  seems to indicate that it is  not necessary for good generalization performance, see \cite{zhang2016understanding}.

A number of attempts have been made to understand why descent algorithms, employed to train overparameterized neural networks, provide a surrogate that generalizes well, see \cite{jacot2018neural},
\cite{dziugaite2017computing},
\cite{arora2018stronger},
\cite{bartlett2017spectrally}. However, the resulting a priori performance guarantees are often vacuous in practice in the sense that the probability of misclassification of a new sample
is bounded from above
by a number greater than one. Of course, no such guarantee can hold in the absence of a  model class assumption on the underlying function $f$ which provided the data.

On the other hand, some heuristic explanations have been put forward to explain the success of this approach.
One of the most popular is that the descent algorithm itself provides a form of  \textit{implicit regularization} that biases learning towards selecting  parameter values $\theta^*$ that correspond in some sense to low complexity functions $S(\cdot ;\theta^*)$. The idea is that the starting guess $S(\cdot;\theta^{(0)})$ has relatively low complexity with high probability. Then, since the model is overparameterized, there are many values  of $\theta$ for which $S
(\cdot;\theta)$ interpolates the data. In particular, there is often such a value  $\theta^*$ near $\theta^{(0)}$. Since gradient descent is essentially a greedy local search, it is reasonable to expect that it will converge to  such a $\theta^*$ that is near $\theta^{(0)}$. 

These heuristics  would match a model class assumption that
$f$ itself is well approximated by the output of neural networks depending on relatively few parameters, that is, $f $ is in a model class $\cA^r$ with
a large value of $r$. Or, more generally, that $f$ is well-approximated by networks depending on many parameters, but with some additional constraints on the size or complexity or these parameters. The purpose of 
the next section is to provide some  support for this idea in the simple case of overparameterized regression.

\subsubsection{Gradient descent for linear regression}
\label{SSS:GDlinear}
As we have seen, the outputs of a neural network form  a complicated nonlinear family which is  difficult to analyze.   It could be therefore
useful to understand what the above numerical approach  based on gradient descent yields  in the simpler case of linear regression.  We briefly describe this in the present section.

We seek to model a data set
\[\{(x^{(i)},f(x^{(i)}))\},\qquad x^{(i)}\in \R^d,\, f(x^{(i)})\in\R,
\quad i=1,\ldots,m,\] 
 by using a function from a linear space  
\[
V_n = \mathrm{span}\set{\phi_j,\, j=1,\ldots, n}.
\]
The key assumption we make is that the model is overparameterized, meaning that $m < n$. If $A:=(a_{ij})$ is the $m\times n$ matrix with entries
\be 
\nonumber
a_{i,j}:=\phi_j(x^{(i)}), \quad i=1,\dots,m;\ j=1,\dots,n,
\ee 
the coefficients $\theta=(\theta_j)_{j=1}^n$ of any interpolant 
\begin{equation}
\label{intr}
S(\cdot,\theta)=\sum_{j=1}^n\theta_j \phi_j(\cdot)\in V_n
\end{equation}
to the data satisfy  the
underdetermined system of equations
\be 
\label{matrixproblem}
A\theta =y, \quad y:=(f(x^{(1)}),\ldots,f(x^{(m)}))\in\R^m.
\ee 
The standard way of choosing a solution to \eref{matrixproblem} is to choose the Moore-Penrose pseudoinverse
$\theta^*\in\R^n$, which we recall is the solution which has minimum $\ell_2(\R^n)$ norm.

Let $W^\perp$ denote the null space, which corresponds to all
$\theta\in\R^n$ which are solutions to \eref{matrixproblem} with the zero vector on the right side.  Further, we denote by $W$ the orthogonal  complement of
$W^\perp$ in $\R^n$. Note that $\theta^*\in W$ since $\theta^*$ is itself a solution to \eref{matrixproblem}.

\vskip .1in
\noindent
{\bf Claim:} Suppose that we apply the gradient descent algorithm,  with  appropriate step
size restrictions,  for finding the minimum of the loss function \eref{lossfunction}, where $S$ is the interpolant \eref{intr}. 
Then, this procedure determines parameter selections $\theta^{(k)}\in \R^n$, $k=1,2,\dots$,  which have  a limit
\be 
\label{limit}
\hat \theta=\lim_{k\to \infty}\theta^{(k)} = \theta^* +P_{W^\perp}\theta^{(0)},
\ee 
where $\theta^{(0)}$ is the initial guess.
\vskip .1in

We do not provide a full detailed proof of this claim, but make the following remarks, which will allow the reader to fill in the details.
The iterative procedure chooses  step sizes $\eta_k$ and defines an optimization trajectory as follows,
\[
\theta^{(k+1)}: = \theta^{(k)} -\eta_k \nabla \mathcal L(\theta^{(k)}),\quad k=0,1,\dots, 
\]
where $\nabla \cL$ is the gradient of the loss function $\cL$.  
Note that $\nabla \mathcal L(\theta)\cdot \theta'=0$ for any $\theta\in \R^n$ and
$\theta'\in W^\perp$ since the function
$h(t):=\mathcal L(\theta+t\theta')$ is a constant function of $t\in \R$. Thus, the vector $\nabla \mathcal L(\theta^{(k)})$, $k=0,1,\ldots$, does not have components in $W^\perp$.  It follows that,
$P_{W^\perp}(\theta^{(k+1)})=P_{W^\perp}(\theta^{(k)})$,\  $k=0,1,\dots$, which gives
\be 
\label{iteration}
\theta^{(k)}= P_W\theta^{(k)}+ P_{W^\perp}(\theta^{(0)}),\quad k=0,1,\dots,
\ee
and 
\be 
\label{lossiteration}
\cL(\theta^{(k)})=\cL(P_W\theta^{(k)}).
\ee 
The function $\cL$ is strictly  convex on $W$ with minimizer $\theta^*$. Since the iterations of gradient descent converge
under restriction on the step size provided by the eigenvalues of $A^TA$, we obtain the claim.

In summary, we find that optimization by gradient descent from a random initialization has at least two important effects. First, the choice of initialization determines the value of the component $P_{W^\perp}(\theta^{(0)})$ not ``seen'' by the data. Its norm is precisely the distance between the $\theta^*$ and $\hat \theta$, which suggests that it is important to properly initialize the optimization. Second, the gradient descent was greedy, leaving $P_{W^ \perp} (\theta^{(0)})$ unchanged during the optimization. This can be viewed as a form of implicit regularization, since at least it does not increase this component. In addition, it  implies an implicit model class assumption that the function $f$ underlying the data $\{(x^{(i)},f(x^{(i)}))\}$ is of low complexity which means that it is well approximated by $V_n$.

\subsubsection{Gradient descent selection for neural networks}
\label{SSS:GDselection}

The above discussion does not carry over directly to overparameterized data fitting with neural networks because the set of NN outputs is not a linear space. However, a recent line of work, see  \cite{jacot2018neural},
\cite{du2018gradient},
\cite{allen2019convergence},
\cite{du2019gradient},
\cite{liu2020toward}, has shown that for sufficiently \textit{wide} networks such considerations are still approximately valid. In short, as we sketch immediately below, a number of rigorous results show that, as $W\gives \infty$, gradient descent on the mean squared error loss $\mathcal L$ using neural networks $\Upsilon^{W,L}(\sigma;d,1)$ can be recast as overparameterized regression in a RKHS $H_{\sigma,L}$, determined by $\sigma$ and $L$. The reproducing kernel of $H_{\sigma,L}$ is called the neural tangent kernel and is fixed throughout training in the limit when $W\gives \infty.$

To explain this point, suppose we are given a dataset  as in \eqref{E:data-def}. Let us fix $L$ and solve the learning problem for this dataset using a class of neural networks $\Upsilon^{W,L}(\sigma;d,1)$ in which $W$ is large. Starting from a random guess $\theta^{(0)}$, the trajectory of the gradient descent on the loss $\mathcal L$, see  \eref{lossfunction},  for the network parameters is given by 
\begin{equation}
\label{E:GD}
\theta^{(t+1)}=\theta^{(t)}-\eta_t \nabla_\theta \mathcal L(\theta^{(t)}).    \end{equation}
 
Varying $W$ changes the number of components of $\theta$. It is convenient to introduce the functions
\[
v_i=v_i(\theta):=S(x^{(i)};\theta),\qquad i=1,\ldots, m, 
\]
which record the values of $S$ on the data set. A simple calculus exercise 
(Taylor's formula) shows that the trajectory of $v_i$ induced by \eqref{E:GD} is 
\begin{equation}\label{E:output-GD}
    v_i^{(t+1)} = v_i^{(t)}-\eta_t \sum_{j=1}^m K_{\theta^{(t)}}(x^{(j)},x^{(i)}) (v_j^{(t)}-y_j) + {\cal O}(\eta_t^2),
\end{equation}
where $v_i^{(t+1)}:=S(x^{(i)};\theta^{(t+1)})$, and 
$K_\theta$ is the so-called \textit{neural tangent kernel}
\[
 K_\theta(x^{(j)},x^{(i)}) := 2 \sum_{l=1}^n \frac{\partial S(x^{(j)};\theta)}{\partial\theta_l}\frac{\partial S(x^{(i)};\theta)}{\partial\theta_l}.
\]
Note that $K_{\theta^{(t)}}$ 
depends on the current setting $\theta^{(t)}$ of trainable parameters. However, it turns out that in the limit when $W$, and hence $n$, tends to infinity, $K_{\theta^{(t)}}$ is given for all $t$ by the average 
\[
K_{\sigma, L}(x,x'):=\E{K_{\theta^{(0)}}(x,x')}, \quad
x,x'\in \R^d,
\]
of $K_{\theta^{(0)}}$ over the randomness in $\theta^{(0)}$. The notation $K_{\sigma, L}$ is meant to emphasize that this limiting kernel depends on the network depth $L$ and the activation function $\sigma$, see \cite{jacot2018neural} and subsequent work. Thus, the training dynamics are summarized by \begin{equation}\label{E:wide-output-GD}
    v_i^{(t+1)} = v_i^{(t)}-\eta_t \sum_{j=1}^m K_{\sigma,L}(x^{(j)},x^{(i)}) (v_j^{(t)}-y_j).
\end{equation}
The term multiplied by $\eta_t$ on the right hand side is precisely the derivative with respect to $v_i$ of 
\[
\norm{v^{(t)}-y}_{K_{\sigma, L}}^2:=\sum_{j,i=1}^m K_{\sigma,L}(x^{(j)},x^{(i)})(v_j^{(t)}-y_j)(v_i^{(t)}-y_i),
\]
where $v^{(t)}-y:=(v_1^{(t)}-y_1,\ldots,v_m^{(t)}-y_m)$ and the norm is with respect to the RKHS structure determined by $K_{\sigma,L}$.

This derivation shows that in the case of small step sizes and large widths, using gradient descent on the loss function ${\cal L}$, see \eref{lossfunction}, for neural networks of fixed depth  
 is similar  to using
gradient descent for the least squares regression problem in the RKHS determined by $K_{\sigma,L}$.

 While the discussion above gives some view of what gradient descent minimization is doing,  a satisfactory understanding  of why overparameterized learning generalizes well remains elusive.  This is an  important but poorly understood topic with a rapidly growing literature, see  \cite{ghorbani2019linearized},
 \cite{bartlett2020benign,chizat2019lazy}.  
  
  \subsubsection{Stability of gradient descent}
  \label{SSS:GDstability}
  
  A natural   question  when  applying gradient descent to find  an approximant to the underlying function is its stability as a numerical algorithm.
  That is, if we slightly change the input data (the data sites and the values  assigned to  these points), how does this effect the output of
  the numerical algorithm.    In this section,  we ask some natural questions that would aid our understanding of stability.
  
  In our earlier treatment of stability, see \S \ref{SS:stablewidths}, we assumed full access to the target function $f$ in the formulation of what stability meant and what was
  an optimal performance of a stable recovery when  using nonlinear manifolds.  Recall that the optimal recovery rate on a model class $K$ was given by the  stable widths $\delta_n^*(K)_X$,  and these  were connected to the entropy of $K$.
 
 Let us denote by $D:=\{(x^{(i)},f(x^{(i)})\}$, $x^{(i)}\in \R^d$, $f(x^{(i)})\in \R$, 
 $i=1,\dots,m$,  the data provided to us.  So, $D$ is a collection of $m$ points in $\R^{d+1}$ and $D$ itself can be viewed as a point in $\R^{(d+1)m}$.  We let $\Sigma_n=\Upsilon^{W,L}:=\Upsilon^{W,L}(\Relu;d,1)$ be the output set  of the neural network architecture that has been chosen.  Here, $n$ is the total number of parameters used to describe the elements of $\Sigma_n$,
  that is, $n=n(W,L)$.  We denote by $a_n:D\mapsto \theta(D)$ the mapping of the data into the parameters   $\theta(D)$  chosen by the numerical algorithm which, for the time being, we assume is based on gradient descent.  Then,   $A_n(D)= M_n(a_n(D)) \in \Sigma_n$ is the output of the algorithm and the learned surrogate.
 \begin{itemize}
 \item {\bf Question 1:}  What  are the regularity properties of  $A_n$?  Is it continuous or perhaps even smoother? 
 
 Of course, the answer will depend on the step size restrictions imposed during the steps of gradient descent and, in addition, on the stopping criteria for the iterations.
 Recall that we know from \S  \ref{SS:stabilityM} that $M_n$ is locally Lipschitz, that is, on each bounded set $B$ in $\R^n$ it is Lipschitz with Lipschitz constant $\gamma_B$.  This leads us to ask the next question.
 
 \item
 {\bf Question 2:}  On which compact sets in $\R^n$ does $M_n$ have a reasonable Lipschitz constant?

 Some information about this question can be extracted from our discussion in \S\ref{SS:stabilityM}, but the analysis there was quite crude.
 Given an answer to {\bf Question 2}, we would like $a_n$ to map into such a ball which leads us to the next question.
 
 \item 
 {\bf Question 3:}  What can be said about the range of $a_n$ as it relates to the initial parameter guess and subsequent step size restrictions?
\end{itemize}
 
Our next questions center on whether $A_n(D)$ is a good surrogate. Although model classes do not appear in the construction of $A_n$,  there is a  belief  that $A_n(D)$ provides a good surrogate for the target function $f$ that gave rise to the data. If this is indeed the case, then
   this statement needs an analytic formulation.  One such possible answer  is that $A_n$ is good
  for a universal collection of model classes.  To try to formulate this, let us now introduce a model class $K$ into the picture, where
  $K\subset X$ is a compact subset of $X$.  We take the view
  that $K$ exists but is unknown to us.
  
  Given such a model class $K$, the datasets given to us are now of the form $D=D(f)$, $f\in K$, where $f(x^{(i)})$ are the observed values
  at the data sites $x^{(i)}$, $i=1,\dots,m$.  We can further add in variability of the data sites by introducing
  $\cX:=(x^{(i)})_{i=1}^m$.  In this way, we can view the data provided to depend on both the selection of sites and the $f\in K$, and write $D(\cX,f)$.
  One can then revisit {\bf Questions 1-3} in this setting.
  
  We can now view $A_n$ as a map $A_n:\cX\times K \to \Sigma_n$ and treat it as a random variable. This would allow us to measure its performance in expectation or with high probability.  At this point, there would be no need to require that the mapping $a_n$ be given by
  gradient descent but rather put gradient descent into competition with more general mappings.  This would lead to various notions of optimal performance similar to those, considered in Information Based Complexity, see e.g. \cite{TW}.   One of these is
  \be
  \label{mnPerformance}
  E_{m,n}(K)_X:=\inf_{\#(\cX)=m;A_n\in\cA}\sup_{f\in K}\|f-A_n(\cX,f)\|_X,
  \ee 
  where the infimum is taken over a class $\cA$ of algorithms $A_n$, perhaps imposing some stability on $A_n$.
  Another meaningful measure of optimality would involve expected performance over random draws $\cX$.
  
  Whatever measure of performance is chosen, one can introduce a corresponding concept of width.  Now, the width $\delta_{m,n}(K)_X$  for a model class $K$ would depend on both 
  $m$ and $n$, and the properties imposed on the algorithms in $\cA$ such as Lipschitz mappings.  
  With such a width in hand, one can now ask for lower and upper bounds for these widths.

 \vskip .1in
 \noindent
 {\bf Acknowledgment: }All three authors were supported by a MURI grant
 N00014-20-1-2787, administered through the Office of Naval Research.  RD and GP were supported by an NSF grant
DMS-1817603, 
BH was supported by an NSF grant DMS–1855684.

\bibliography{Refs}

\label{lastpage}

\end{document}